\documentclass[reqno,10pt]{amsart}
\usepackage{amscd,amssymb,verbatim}
\numberwithin{equation}{section} \theoremstyle{plain}

\newcommand\alp{\alpha}         
\newcommand\bet{\beta}
\newcommand\gam{\gamma}         \newcommand\Gam{\Gamma}
\newcommand\del{\delta}         \newcommand\Del{\Delta}
\newcommand\eps{\varepsilon}
\newcommand\zet{\zeta}
\newcommand\tet{\theta}

\newcommand\lam{\lambda}                \newcommand\Lam{\Lambda}
\newcommand\sig{\sigma}         

\newcommand\ome{\omega}         \newcommand\Ome{\Omega}

\newcommand\calB{{\mathcal{B}}}

\newcommand\calF{{\mathcal{F}}}
\newcommand\calG{{\mathcal{G}}}

\newcommand\calI{{\mathcal{I}}}

\newcommand\calK{{\mathcal{K}}}

\newcommand\calM{{\mathcal{M}}}
\newcommand\calN{{\mathcal{N}}}

\newcommand\calR{{\mathcal{R}}}

\newcommand\bfk{{\mathbf k}}

\newcommand\RR{\mathbb{R}}

\newcommand\ZZ{\mathbb{Z}}

\newcommand\CC{\mathbb{C}}

\newcommand\nek{,\ldots,}
\newcommand\sdp{\times \hskip -0.3em {\raise 0.3ex
\hbox{$\scriptscriptstyle |$}}} 

\newcommand\const{\operatorname{const}}

\newcommand\Det{\operatorname{Det}}

\newcommand\Hom{\operatorname {Hom}}
\newcommand\Id{\operatorname {Id}}

\newcommand\IM{\operatorname{Im}}

\newcommand\MOD{\operatorname{mod}}

\newcommand\Ker{\operatorname{Ker}}

\newcommand\rank{\operatorname{rank}}
\newcommand\rk{\operatorname{rk}}

\newcommand\RE{\operatorname{Re}}

\newcommand{\sign}{\operatorname{sign}}

\newcommand\spec{{\rm spec\,}}

\newcommand\tr{\operatorname{tr}}
\newcommand\Tr{\operatorname{Tr}}

\newcommand\os{{\bar{s}}}

\newcommand\oz{{\overline{z}}}

\newcommand\oeta{{\overline{\eta}}}

\newcommand\oxi{{\overline{\xi}}}

\newcommand\ozet{\overline{\zet}}

\newcommand\hata{{\hat{a}}}
\newcommand\hatA{{\widehat{A}}}

\newcommand\hatc{{\hat{c}}}
\newcommand\hatC{{\widehat{C}}}
\newcommand\hatB{{\widehat{B}}}

\newcommand\hath{{\hat{h}}}
\newcommand\hatH{{\widehat{H}}}

\newcommand\hatV{{\widehat{V}}}

\newcommand\hatGam{{\widehat{\Gamma}}}

\newcommand\hatOme{{\widehat{\Ome}}}

\newcommand\tilA{{\widetilde{A}}}
\newcommand\tila{{\tilde{a}}}
\newcommand\tilB{{\widetilde{B}}}

\newcommand\tilc{{\tilde{c}}}
\newcommand\tilC{{\widetilde{C}}}

\newcommand\tilh{{\tilde{h}}}

\newcommand\tilH{{\widetilde{H}}}

\newcommand\tiln{{\tilde{n}}}

\newcommand\tilDel{{\widetilde{\Delta}}}

\newcommand\tilGam{{\widetilde{\Gamma}}}

\renewcommand{\>}{\rangle}
\newcommand{\<}{\langle}

\theoremstyle{plain}
\newtheorem{Thm}[subsection]{Theorem}
\newtheorem{Cor}[subsection]{Corollary}
\newtheorem{Lem}[subsection]{Lemma}
\newtheorem{Prop}[subsection]{Proposition}
\newtheorem{Conjec}[subsection]{Conjecture}

\newtheorem{Def}[subsection]{Definition}

\theoremstyle{remark}

\newtheorem{Rem}[subsection]{Remark}

\def\TeXref#1{%
        \leavevmode\vadjust{\setbox0=\hbox{{\tt
                \  {\tiny \textrm #1}}}%
        \theight=\ht0
        \advance\theight by \lineskip
        \kern -\theight \vbox to
        \theight{\rightline{\rlap{\box0}}%
        \vss}%
        }}%
\newif\ifShowLabels
\ShowLabelstrue
\newdimen\theight
\def\TeXrefEq#1{%
        \leavevmode\vadjust{\setbox0=\hbox{{\tt
                \  {\tiny \textrm #1}}}%
        \theight=\ht1
        \advance\theight by \lineskip
        \kern -\theight \vbox to
        \theight{\rightline{\rlap{\box0}}%
        \vss}%
        }}%

\newcommand{\refs}[1]{Section ~\ref{S:#1}}
\newcommand{\refss}[1]{Subsection ~\ref{SS:#1}}

\newcommand{\reft}[1]{Theorem ~\ref{T:#1}}
\newcommand{\refl}[1]{Lemma ~\ref{L:#1}}
\newcommand{\refp}[1]{Proposition ~\ref{P:#1}}

\newcommand{\refd}[1]{Definition ~\ref{D:#1}}
\newcommand{\refr}[1]{Remark ~\ref{R:#1}}
\newcommand{\refe}[1]{\eqref{E:#1}}

\newenvironment{thm}[1]%
        { \begin{Thm} \label{T:#1}  \ifShowLabels \TeXref{T:#1} \fi }%
        { \end{Thm} }

\renewcommand{\th}[1]{\begin{thm}{#1}  }
\renewcommand{\eth}{\end{thm} }

\newenvironment{lemma}[1]%
        { \begin{Lem} \label{L:#1}  \ifShowLabels \TeXref{L:#1} \fi }%
        { \end{Lem} }

\newcommand{\lem}[1]{\begin{lemma}{#1} }
\newcommand{\elem}{\end{lemma}}

\newenvironment{propos}[1]%
        { \begin{Prop} \label{P:#1}  \ifShowLabels \TeXref{P:#1} \fi }%
        { \end{Prop} }

\newcommand{\prop}[1]{\begin{propos}{#1} }
\newcommand{\eprop}{\end{propos}}

\newenvironment{corol}[1]%
        { \begin{Cor} \label{C:#1}  \ifShowLabels \TeXref{C:#1} \fi }%
        { \end{Cor} }
\newcommand{\cor}[1]{\begin{corol}{#1}  }
\newcommand{\ecor}{\end{corol}}

\newenvironment{conjec}[1]%
        { \begin{Conjec} \label{Conj:#1}  \ifShowLabels \TeXref{C:#1} \fi }%
        { \end{Conjec} }
\newcommand{\conj}[1]{\begin{conjec}{#1}  }
\newcommand{\econj}{\end{conjec}}

\newenvironment{defeni}[1]%
        { \begin{Def} \label{D:#1}  \ifShowLabels \TeXref{D:#1} \fi }%
        { \end{Def} }
\newcommand{\defe}[1]{\begin{defeni}{#1}  }
\newcommand{\edefe}{\end{defeni}}

\newenvironment{remark}[1]%
        { \begin{Rem} \label{R:#1}  \ifShowLabels \TeXref{R:#1} \fi }%
        { \end{Rem} }
\newcommand{\rem}[1]{\begin{remark}{#1}}
\newcommand{\erem}{\end{remark}}

\newcommand{\eq}[1]%
        { \ifShowLabels \TeXrefEq{E:#1} \fi
           \begin{equation} \label{E:#1} }
\newcommand{\eeq}{\end{equation}}
\newcommand{\meq}[1]%
        { \ifShowLabels \TeXrefEq{E:#1} \fi
           \begin{multline} \label{E:#1} }
\newcommand{\emeq}{\end{multline}}

\newcommand{\prf}{ \begin{proof} }
\newcommand{\eprf}{ \end{proof} }
\newcommand{\Label}[1]{\label{#1}  \ifShowLabels \TeXref{#1} \fi }

\ShowLabelsfalse

\newcommand{\n}{\nabla}

\renewcommand{\tiln}{\tilde{\n}}

\renewcommand{\b}{\bullet}
\newcommand{\pa}{\text{\( \partial\)}}
\newcommand{\tilpa}{\text{\( \tilde{\partial}\)}}\newcommand{\hatpa}{\text{\( \hat{\partial}\)}}
\newcommand{\odd}{{\operatorname{odd}}}
\newcommand{\even}{{\operatorname{even}}}
\newcommand{\Fp}{\operatorname{F.p._{s=0}}}

\newcommand{\TRS}{T^{\operatorname{RS}}}

\newcommand{\Detgrtet}{\Det_{{\operatorname{gr},\tet}}'}
\newcommand{\Detgrnp}{\Det_{\operatorname{gr}}}
\newcommand{\Detgrtetnp}{\Det_{{\operatorname{gr},\tet}}}

\newcommand{\symb}{\sig}

\newcommand{\LD}{\operatorname{LDet}'}\newcommand{\LDnp}{\operatorname{LDet}}

\newcommand{\LDgrtet}{\operatorname{LDet}_{\operatorname{gr},\tet}'}\newcommand{\LDgrtetnp}{\operatorname{LDet}_{\operatorname{gr},\tet}}

\newcommand{\fdm}{\text{\( \|\cdot\|_{\Det(H^{\b}(\pa))}\)}}
\newcommand{\B}{\calB}\newcommand{\tB}{\tilde{\calB}}

\newcommand{\rat}{\rho_{\operatorname{an}}}
\newcommand{\RSn}[1]{\|#1\|_{\Det(H^\b(M,E))}^{\operatorname{RS}}}
\newcommand{\RSnp}[1]{\|#1\|_{\Det(H^\b(M,E'))}^{\operatorname{RS}}}
\newcommand{\RSnnp}[1]{\|#1\|_{\Det(H^\b(M,E\oplus{}E'))}^{\operatorname{RS}}}
\newcommand{\bigRSnnp}[1]{\big\|#1\big\|_{\Det(H^\b(M,E\oplus{}E'))}^{\operatorname{RS}}}

\newcommand{\bigRSnp}[1]{\big\|#1\big\|_{\Det(H^\b(M,E'))}^{\operatorname{RS}}}

\renewcommand{\hatc}{\hat{c}}
\newcommand{\tV}{V^{*_\tau}}

\newcommand{\hGam}{{\Gam^*}}
\renewcommand{\hatOme}{\check{\Ome}}

\begin{document}

\title{Refined Analytic Torsion as an Element of the Determinant Line}
\author[Maxim Braverman]{Maxim Braverman$^\dag$}
\address{Department of Mathematics\\
        Northeastern University   \\
        Boston, MA 02115 \\
        USA
         }
\email{maximbraverman@neu.edu}
\author[Thomas Kappeler]{Thomas Kappeler$^\ddag$}
\address{Institut fur Mathematik\\
         Universitat Z\"urich\\
         Winterthurerstrasse 190\\
         CH-8057 Z\"urich\\
         Switzerland
         }
\email{tk@math.unizh.ch}
\thanks{${}^\dag$Supported in part by the NSF grant DMS-0204421.\\
\indent${}^\ddag$Supported in part by the Swiss National Science foundation, the programme SPECT, and the European Community through the FP6
Marie Curie RTN ENIGMA (MRTN-CT-2004-5652)} 
\begin{abstract}
We construct a canonical element, called the refined analytic torsion, of the determinant line of the cohomology of a closed oriented odd-dimensional
manifold $M$ with coefficients in a flat complex vector bundle $E$. We compute the Ray-Singer norm of the refined analytic torsion. In particular, if
there exists a flat Hermitian metric on $E$, we show that this norm is equal to 1. We prove a duality theorem, establishing a relationship between
the refined analytic torsions corresponding to a flat connection and its dual.
\end{abstract}
\maketitle \tableofcontents

\section{Introduction}\Label{S:introduction}

The aim of this paper is to define the notion of the {\em refined analytic torsion} associated to an arbitrary vector bundle $E\to M$ over a closed
oriented manifold of {\em odd} dimension $d=2r-1$ and an arbitrary flat connection $\n$ on $E$. The refined analytic torsion $\rat=\rat(\n)$ is an
element of the determinant line $\Det\big(H^\b(M,E)\big)$ of the cohomology of the bundle $E$. If $\dim{M}\equiv 1\ (\MOD 4)$ or if the rank of the
bundle $E$ is divisible by 4, then $\rat(\n)$ is independent of any choices. If $\dim{M}\equiv 3\ (\MOD 4)$ then $\rat(\n)$ depends on the choice of
a compact oriented manifold $N$, whose oriented boundary is diffeomorphic to two disjoint copies of $M$, but only up to a factor $i^{k\cdot\rk{E}}$ \
($k\in \ZZ$).

The Ray-Singer norm of the refined analytic torsion is equal to $e^{\pi\IM\eta}$, where $\eta$ is the $\eta$-invariant of the Atiyah-Patodi-Singer
odd signature operator $\B$ associated to $\n$ and a Riemannian metric on $M$. In particular, if the connection $\n$ is Hermitian then $\eta$ is real
and the Ray-Singer norm of $\rat(\n)$ is equal to 1. This property justifies calling $\rat(\n)$ the refined analytic torsion. Indeed, the Ray-Singer
torsion can be viewed as an element of the determinant line $\Det\big(H^\b(M,E)\big)$ whose Ray-Singer norm is equal to one. Such an element is
defined up to a multiplication by $t\in\CC$ with $|t|=1$.  The refined analytic torsion $\rat(\n)$ is a canonical choice of such an element.

If the connection $\n$ is acyclic, i.e., if the cohomology $H^\b(M,E)=0$, then $\Det\big(H^\b(M,E)\big)$ is canonically isomorphic to $\CC$ and the
refined analytic torsion is a complex number. Under the additional assumption that $\B$ is bijective, this number was introduced and studied in
\cite{BrKappelerRATshort,BrKappelerRAT}. In this paper we drop this assumption and define the refined analytic torsion for arbitrary flat connections
and compare it to the Ray-Singer torsion. Recall that in \cite{BrKappelerRATshort,BrKappelerRAT}, we also computed the quotient of the refined
analytic torsion and Turaev's refinement of the combinatorial torsion (see also \cite{BurgheleaHaller_function} for related results). In our
subsequent paper \cite{BrKappelerRATdetline_hol}, we extend this computation to arbitrary flat connections.

Having applications in topology in mind Dan Burghelea suggested constructing a holomorphic function on the space of acyclic connections, whose
absolute value is (related to) the Ray-Singer torsion. Such functions were constructed by Burghelea and Haller in
\cite{BurgheleaHaller_Euler,BurgheleaHaller_function} and also in our papers \cite{BrKappelerRATshort,BrKappelerRAT}. The results of the present
paper fit nicely into this program. In fact, in \cite{BrKappelerRATdetline_hol}, we show that the refined analytic torsion is a holomorphic section
of the determinant line bundle over the space of flat connections (or, equivalently, over the space of representations of the fundamental group of
$M$). Thus we extend the construction of \cite{BrKappelerRATshort,BrKappelerRAT} to general, not necessarily acyclic, connections.

Let us now briefly describe our construction of the refined analytic torsion.

\subsection{The refined torsion of a finite dimensional complex}\Label{SS:reftorfd}
Let
\eq{Cpintrod}
    \begin{CD}
       (C^\b,\pa):\quad  0 \ \to C^0 @>{\pa}>> C^1 @>{\pa}>>\cdots @>{\pa}>> C^d \ \to \ 0
    \end{CD}
\end{equation}
be a complex {\em finite dimensional}\/ $\CC$-vector spaces of odd length $d=2r-1$. A {\em chirality operator}\/ $\Gam:C^\b\to C^\b$ is an involution
such that $\Gam(C^j)= C^{d-j}$, for all $j=0\nek d$. As a first step towards our construction of the refined analytic torsion we introduce and study
the {\em refined torsion} of the pair $(C^\b,\Gam)$. Consider the determinant line
\[
    \Det(C^\b) \ := \ \bigotimes_{j=0}^d\, \Det(C^j)^{(-1)^j},
\]
where $\Det(C^j)^{-1}:= \Hom(C^j,\CC)$ denotes the dual of $C^j$. For an element $c_j\in \Det(C^j)$ we denote by $c_j^{-1}$  the unique element in
$\Det(C^j)^{-1}$ satisfying $c_j^{-1}(c_j)=1$. We also denote by $\Gam{}c_j\in \Det(C^{d-j})$ the image of $c_j$ under the map $\Det(C^j)\to
\Det(C^{d-j})$ induced by $\Gam:C^j\to C^{d-j}$.

For each $j=0\nek r-1$, fix an element $c_j\in \Det(C^j)$ and consider the element
\eq{Gamcintrod}\notag
    c_{{}_\Gam} \ := \ (-1)^{\calR(C^\b)}\cdot
    c_0\otimes c_1^{-1}\otimes \cdots \otimes c_{r-1}^{(-1)^{r-1}}\otimes (\Gam c_{r-1})^{(-1)^r}\otimes (\Gam c_{r-2})^{(-1)^{r-1}}
    \otimes \cdots\otimes (\Gam c_0)^{-1}
\end{equation}
of $\Det(C^\b)$, where $(-1)^{\calR(C^\b)}$ is a normalization factor introduced in \refe{Gamc}. It is easy to see that $c_{{}_\Gam}$ is independent
of the choice of $c_0\nek c_{r-1}$.

Denote by $H^\b(\pa)$ the cohomology of the complex $(C^\b,\pa)$. In \refss{isomorphismn}, we construct a sign refined version of the standard
isomorphism $\Det(C^\b)\to \Det(H^\b(\pa))$, cf. \cite{Milnor66}. Our isomorphism $\phi_{C^\b}:\Det(C^\b)\to \Det(H^\b(\pa))$ is similar but not
equal to the one considered by Turaev \cite{Turaev86}.
\defe{refinedtorsionintrod}
The {\em refined torsion} of the pair $(C^\b,\Gam)$ is the element
\eq{refinedtorintrod}\notag
    \rho_{{}_\Gam} \ = \ \rho_{{}_{C^\b,\Gam}} \ := \  \phi_{C^\b}(c_{{}_\Gam})\ \in \ \Det(H^\b(\pa)).
\end{equation}
\edefe

\subsection{Calculation of the refined torsion of a finite dimensional complex}\Label{SS:calcreftorfd}
To compute the refined torsion we introduce the operator
\[
    \B \ := \ \Gam\,\pa\ + \pa\,\Gam.
\]
This operator is a finite dimensional analogue of the signature operator on an odd-dimensional manifold, see \cite[p.~44]{APS1}, \cite[p.~405]{APS2},
\cite[pp.~64--65]{Gilkey84}, and section~7 of this paper.  For $j=0\nek d$, define
\[
    C^j_+ \ := \ \Ker \big(\,\pa\circ\Gam\,\big)\cap C^j, \qquad C^j_- \ := \ \Ker\pa\cap C^j
\]
and  set $C^{-1}_+= C^{d+1}_-= 0$. Let $\B_j$ and $\B_j^\pm$ denote the restriction of $\B$ to $C^j$ and $C^j_\pm$ respectively. Then, for each
$j=0\nek d$, one has
\eq{Bkpmintrod}\notag
    \B_j^+ \ = \ \Gam\circ\pa:\, C^j_+ \ \longrightarrow \ C^{d-j-1}_+, \qquad
    \B_j^- \ = \ \pa\circ\Gam:\, C^j_- \ \longrightarrow \ C^{d-j+1}_-.
\end{equation}
Denote $C^\even:= \bigoplus_{j\ \even}\,C^j$,\ $C^\even_\pm:= \bigoplus_{j\ \even}\,C^j_\pm$. Set
\eq{Bevenoddintrod}\notag
    \B_\even\ :=\ \bigoplus_{j\ \even}\, \B_j:\,C^\even\ \to C^\even, \qquad \B_\even^\pm\ :=\ \bigoplus_{j\ \even} \B_j^\pm:\, C^\even_\pm\ \to \ C^\even_\pm,
\end{equation}
and define $\B_\odd, \ \B_\odd^\pm$ similarly. As $\B_\even = \Gam\circ{}\B_\odd\circ{}\Gam$, it turns out that it suffices to study $\B_\even$.

Suppose, first, that the signature operator $\B$ is bijective. Then, cf. \refl{Bsym-Cacyclic}, the complex $(C^\b,\pa)$ is acyclic. Hence,
$\Det(H^\b(\pa))$ is canonically identified with $\CC$ and the refined torsion $\rho_{{}_\Gam}$ can be viewed as a number in $\CC$. In
\refp{reftor-grdetfd} we compute this number and show that
\[
    \rho_{{}_\Gam} \ = \ \Detgrnp(\B_\even),
\]
where the {\em graded determinant}\/ $\Detgrnp(\B_\even)$ is defined by the formula
\[
    \Detgrnp(\B_\even) \ \overset{\text{Def}}{=}\ \Det(\B_\even^+)/\Det(-\B_\even^-).
\]
Note that in our definition of the graded determinant the quotient is the determinant of the {\em negative} of\/ $\B_\even^-$.

To compute the refined torsion in the case $\B$ is {\em not}\/ bijective, note that $\B^2=(\Gam\pa)^2+(\pa\Gam)^2$ maps $C^j$ ($j=0\nek d$) into
itself. For each $j=0\nek d$ and an arbitrary interval $\calI$, denote by $C^j_\calI\subset C^j$ the linear span of the generalized eigenvectors of
the restriction of $\B^2$ to $C^j$, corresponding to eigenvalues $\lam$ with $|\lam|\in \calI$. Since both operators, $\Gam$ and $\pa$, commute with
$\B$ (and, hence, with $\B^2$), $\Gam(C^j_\calI)\subset C^{d-j}_\calI$ and $\pa(C^j_\calI)\subset C^{j+1}_\calI$. Hence, we obtain a subcomplex
$C^\b_\calI$ of $C^\b$ and the restriction $\Gam_\calI$ of $\Gam$ to $C^\b_\calI$ is a chirality operator for $C^\b_\calI$. We denote by
$H^\b_\calI(\pa)$ the cohomology of the complex $(C^\b_\calI,\pa_\calI)$.

Denote by $\pa_\calI$ and $\B_\calI$ the restrictions of $\pa$ and $\B$ to $C^\b_\calI$. Then $B_\calI= \Gam_\calI\pa_\calI+\pa_\calI\Gam_\calI$ and
one easily shows (cf. \refl{BI>0}) that $(C^\b_\calI,\pa_\calI)$ is acyclic if $0\not\in\calI$.

For each $\lam\ge0$, $C^\b= C^\b_{[0,\lam]}\oplus C^\b_{(\lam,\infty)}$ and $H^\b_{(\lam,\infty)}(\pa)= 0$ whereas $H^\b_{[0,\lam]}(\pa)\simeq
H^\b(\pa)$. Hence, there are canonical isomorphisms
\[
        \Phi:\, \Det(H^\b_{(\lam,\infty)}(\pa))\ \longrightarrow\ \CC, \qquad \Psi:\,\Det(H^\b_{[0,\lam]}(\pa))\ \longrightarrow \
        \Det(H^\b(\pa)).
\]
In the sequel, we will write $t$ for $\Phi(t)\in \CC$ and denote by $h$ also the element $\Psi(h)\in \Det(H^\b(\pa))$. Then, cf.
\refp{antor-grdetfd2}, for any $\lam\ge0$, the refined torsion can be computed to be
\eq{antor-grdetfd2introd}
    \rho_{{}_\Gam} \ = \ \Detgrnp(\B_{(\lam,\infty)})\cdot\rho_{{}_{\Gam_{\hskip-1pt{}_{[0,\lam]}}}}.
\end{equation}
In particular, the element $\Detgrnp(\B_{(\lam,\infty)})\cdot\rho_{{}_{\Gam_{\hskip-1pt{}_{[0,\lam]}}}}$ is independent of $\lam$. It is this
property which allows us to define the refined analytic torsion.

\subsection{The canonical element in the determinant line of the cohomology of a flat vector bundle over a Riemannian manifold}\Label{SS:canelementintrod}
In the second part of this paper (Sections~\ref{S:determinants}--\ref{S:Ray-Singer}) we apply the notion of refined torsion to define and investigate
the refined analytic torsion of the (twisted) de Rham complex. Let $E\to M$ be a complex vector bundle over a closed manifold of {\em odd} dimension
$d=2r-1$ and let $\n$ be a flat connection on $E$. Let $\Ome^\b(M,E)$ denote the de Rham complex of $E$-valued differential forms on $M$. For a given
Riemannian metric $g^M$ on $M$ denote by
\[
    \Gam \ = \ \Gam(g^M):\, \Ome^\b(M,E) \ \longrightarrow \ \Ome^\b(M,E)
\]
the chirality operator (cf. \cite[\S3]{BeGeVe}), defined in terms of the Hodge $*$-operator by the formula
\eq{Gamintrod}\notag
    \Gam\, \ome \ := \ i^r\,(-1)^{\frac{k(k+1)}2}\,*\,\ome, \qquad \ome\in \Ome^k(M,E).
\end{equation}
The odd signature operator introduced by Atiyah, Patodi, and Singer \cite{APS1,APS2} (see also \cite{Gilkey84}) is the first order elliptic
differential operator $\B:\Ome^\b(M,E)\to \Ome^\b(M,E)$, given by
\[
    \B \ \overset{\text{Def}}{=} \ \Gam\,\n\ +\ \n\,\Gam.
\]
Notice, that $\B^2$ maps $\Ome^j(M,E)$ into itself for every $j=0\nek d$.

For an interval $\calI\subset [0,\infty)$ we denote by $\Ome^j_\calI(M,E)$ the image of the spectral projection of $\B^2$ corresponding to the
eigenvalues whose absolute value lies in $\calI$, cf. Subsections~\ref{SS:spectralsubspace} and \ref{SS:decompos} for details. The space
$\Ome^j_\calI(M,E)$ contains the span of the generalized eigenforms of $\B^2$ corresponding to eigenvalues whose absolute value lies in $\calI$ and
coincides with this span if the interval $\calI$ is bounded. In particular, if $\calI$ is bounded, then  the dimension of $\Ome^j_\calI(M,E)$ is
finite. Note that, since $\B^2$ and $\n$ commute, the space $\Ome^j_\calI(M,E)$ is a subcomplex of the de Rham complex $\Ome^\b(M,E)$.

For each $\lam\ge0$, we have
\eq{Ome=Ome0+Ome>0introd}\notag
    \Ome^\b(M,E) \ = \  \Ome^\b_{[0,\lam]}(M,E)\,\oplus\,\Ome^\b_{(\lam,\infty)}(M,E).
\end{equation}
The complex $\Ome^\b_{(\lam,\infty)}(M,E)$ is clearly acyclic. Hence, the cohomology $H^\b_{[0,\lam]}(M,E)$ of the complex $\Ome^\b_{[0,\lam]}(M,E)$
is naturally isomorphic to the cohomology $H^\b(M,E)$ of $\Ome^\b(M,E)$. Further, as $\Gam$ commutes with $\B^2$, it preserves the space
$\Ome_{[0,\lam]}(M,E)$ and the restriction $\Gam_{\hskip-1pt{}_{[0,\lam]}}$ of $\Gam$ to this space is a chirality operator on
$\Ome^\b_{[0,\lam]}(M,E)$.

To define the refined analytic torsion we need to introduce the notion of  a graded determinant of $\B$. For every interval $\calI\subset [0,\infty)$
and for each $k=0\nek d$, set
\eq{ome+-introd}\notag
 \begin{aligned}
  \Ome^k_{+,\calI}(M,E) \ &:= \ \Ker\,(\n\,\Gam)\,\cap\,\Ome^k_{\calI}(M,E) \ = \ \big(\,\Gam\,(\Ker\,\n)\,\big)\,\cap\,\Ome^k_{\calI}(M,E);\\
  \Ome^k_{-,\calI}(M,E) \ &:= \ \Ker\,(\Gam\,\n)\,\cap\,\Ome^k_{\calI}(M,E)
  \ = \ \Ker\,\n\,\cap\,\Ome^k_{\calI}(M,E).
 \end{aligned}
\end{equation}
If $0\not\in \calI$, then, clearly,
\eq{Ome>0=directsumintrod}\notag
    \Ome^k_{\calI}(M,E) \ = \ \Ome^k_{+,\calI}(M,E) \,\oplus \Ome^k_{-,\calI}(M,E).
\end{equation}
The latter decomposition is considered as a {\em grading} on $\Ome^\b_{\calI}(M,E)$. As both, $\Gam$ and $\n$, commute with $\B^2$, we conclude that
for all $j=0\nek d$,
\[
  \begin{aligned}
    \Gam&:\,\Ome^k_{\pm,\calI}(M,E) \ \longrightarrow \ \Ome^{d-k}_{\mp,\calI}(M,E),\\
    \n&:\,\Ome^k_{\pm,\calI}(M,E) \ \longrightarrow \ \Ome^{k+1}_{\mp,\calI}(M,E).
  \end{aligned}
\]

Denote by $\B^\calI_\even$ and $\B^{\pm,\calI}_{\even}$ the restrictions of\/ $\B$ to $\Ome^\even_\calI(M,E):=
\bigoplus_{p=0}^{r-1}\Ome^{2p}_\calI(M,E)$ and\linebreak $\Ome^{\even,\pm}_\calI(M,E):= \bigoplus_{p=0}^{r-1}\Ome^{2p}_{\pm,\calI}(M,E)$,
respectively.

Let $\tet\in (-\pi,0)$ be an Agmon angle for $\B^\calI$, cf. \refd{Agmon}. For each $\calI$ with $0\not\in\calI$ define the graded determinant of the
operator $\B_\even^\calI$ by the formula
\eq{Det=DetDetintrod}
    \Detgrtetnp(\B_\even^\calI) \ \overset{\text{Def}}{=}\ \Det_\tet(\B_\even^{+,\calI})/\Det_\tet(-\B_\even^{-,\calI}),
\end{equation}
where $\Det_\tet$ denotes the $\zet$-regularized determinant, cf. Sections~\ref{S:determinants} and \ref{S:grdet} for details. One verifies easily
that for any $0\le \lam\le \mu<\infty$,
\eq{DetlamDetmuintrod}
    \Detgrtetnp(\B^{(\lam,\infty)}_\even) \ = \ \Detgrtetnp(\B^{(\lam,\mu]}_\even)\cdot \Detgrtetnp(\B^{(\mu,\infty)}_\even).
\end{equation}
For any given $\lam\ge0$, denote by $\rho_{{}_{\Gam_{\hskip-1pt{}_{[0,\lam]}}}}$ the refined torsion of the finite dimensional complex
$\Ome^\b_{[0,\lam]}(M,E)$ and the chirality operator $\Gam_{\hskip-1pt{}_{[0,\lam]}}$. In view of \refe{antor-grdetfd2introd} and
\refe{Det=DetDetintrod}, the product
\eq{rhointrod}
    \rho \ = \ \rho(\n,g^M) \ := \ \Detgrtetnp(\B^{(\lam,\infty)}_{\even})\cdot \rho_{{}_{\Gam_{\hskip-1pt{}_{[0,\lam]}}}}
    \ \in  \ \Det(H^\b(M,E))
\end{equation}
is independent of the choice of $\lam\ge0$. It is also independent of the choice of the Agmon angle $\tet\in (-\pi,0)$ of $\B_\even$.

\subsection{The metric anomaly of $\rho(\n,g^M)$}\Label{SS:anomalyintrod}
The element $\rho(\n,g^M)$ is very close to our notion of the refined analytic torsion. However, in general, it is not a differential invariant of
the flat bundle $E$, since it does depend on the choice of the Riemannian metric $g^M$. To compute the metric anomaly of $\rho(\n,g^M)$ we show (cf.
\refp{xietad}) that
\eq{xietadintrod}\notag
    \Detgrtetnp(\B_\even^{(\lam,\infty)}) \ =\ \exp\Big(\, \xi_\lam \ - \ i\pi\,\eta_\lam \ - \ \frac{i\pi}2\,\sum_{j=0}^d\,(-1)^jj\,d_{j,\lam}\,\Big).
\end{equation}
Here
\eq{xiintrod}\notag
    \xi_\lam \ = \ \frac12\,\sum_{j=0}^{d-1}\,(-1)^{j+1}\frac{d}{ds}{\Big|_{s=0}}\,\zet_{2\tet} \Big(\,s,
    {(\Gam\n)^2}{\big|_{\Ome^j_{+,{(\lam,\infty)}}(M,E)}}\,\Big),
\end{equation}
where $\tet\in (-\pi/2,0)$ is an Agmon angle for $\B$ such that there are no eigenvalues of $\B$ in the solid angles $L_{(-\pi/2,\tet]}$ and
$L_{(\pi/2,\tet+\pi]}$,
\[
    \eta_\lam \ = \ \eta\big(\,\B^{(\lam,\infty)}_\even\,\big)
\]
is the $\eta$-invariant of $\B^{(\lam,\infty)}_\even$, and
\eq{djlamintrod}\notag
    d_{j,\lam} \ := \ \dim \Ome^j_{[0,\lam]}(M,E).
\end{equation}

We then analyze the dependence on the Riemannian metric separately for $\xi_\lam\cdot\rho_{{}_{\Gam_{\hskip-1pt{}_{[0,\lam]}}}}$ and $\eta_\lam$.
The metric anomaly for $\eta_\lam$ has been already computed by Gilkey \cite{Gilkey84}. Using his result, we compute the metric anomaly of
$\rho(\n,g^M)$. The refined analytic torsion is then defined by correcting $\rho(\n,g^M)$ by its anomaly. More precisely, we prove (cf.
\reft{metricindep}) the following
\th{metricindepintrod}
Let $E$ be a flat vector bundle over a closed oriented odd-dimensional manifold $M$ and let $\n$ denote the flat connection on $E$. Let $N$ be an
oriented manifold whose oriented boundary is the disjoint union of two copies of $M$. Then the product
\eq{metricindepintrod}
    \rho(\n,g^M)\cdot e^{i\pi\,\frac{\rank E}2\,\int_N\,L(p,g^M)}\ \in \ \Det(H^\b(M,E)),
\end{equation}
is independent of the metric $g^M$. Here $L(p,g^M)$ is the Hirzebruch $L$-polynomial in the Pontrjagin forms of any Riemannian metric on $N$ which
near $M$ is the product of $g^M$ and the standard metric on the half-line.

In particular, if\/ $\dim{}M\equiv 1 (\MOD 4)$, then $\int_NL(p,g^M)=0$ and, hence, $\rho(\n,g^M)$ is independent of $g^M$.
\eth

\subsection{Definition of the refined analytic torsion}\Label{SS:refantorintrod}
We now define the refined analytic torsion $\rat= \rat(\n)$ to be the element \refe{metricindepintrod} of $\Det(H^\b(M,E))$.

The refined analytic torsion is independent of the choice of the Agmon angle $\tet\in (-\pi,0)$ and of the metric $g^M$.  {\em It does depend,
however, on the choice of the manifold $N$}, but only up to a factor $i^{k\cdot\rank{E}}$ $(k\in \ZZ)$, cf.  the discussion at the end of
\refss{Mnotbounds}. If\/ $\dim{}M\equiv 1 (\MOD 4)$, then $\rat= \rho(\n,g^M)$ and $\rat$ is independent of any choices.

\subsection{The Ray-Singer metric of the refined analytic torsion}\Label{SS:Ray-Singerintrod}
Let $\RSn{\cdot}$ denote the Ray-Singer norm on the determinant line $\Det(H^\b(M,E))$, cf. \cite{BisZh92} and \refss{RSnorm}. In \refs{duality}, we
compute the refined analytic torsion associated to the connection $\n'$ on $E$ dual to a given connection $\n$. In \refs{Ray-Singer}, we use this
calculation to calculate the Ray-Singer norm of the refined analytic torsion. More precisely we prove (cf. \reft{RSnormRAT}) the following
\th{RSnormRATintrod}
Let $E$ be a complex vector bundle over a closed oriented odd-dimensional manifold $M$ and let $\n$ be a flat connection on $E$. Then
\eq{RSnormRATintrod}\notag
    \RSn{\rat} \ = \ e^{\pi\,\IM \eta(\B_\even(\n,g^M))}.
\end{equation}
In particular, if\/ $\n$ is a Hermitian connection, then the operator $\B_\even(\n,g^M)$ is self-adjoint, hence, its $\eta$-invariant
$\eta(\B_\even(\n,g^M))$ is real, and
\eq{RSnormRAThermintrod}\notag
    \RSn{\rat} \ = \ 1.
\end{equation}
\eth

\subsection*{Acknowledgements}
We would like to thank Guangxiang Su for pointing out a sign mistake in a preliminary version of the paper. The first author would like to thank the Institut des
Hautes \'Etudes Scientifiques, where part of this work was completed, for hospitality and providing excellent working conditions.

\section{The Determinant Line of a Finite Dimensional Complex}\Label{S:FinDim}

In this section we review some standard material about determinant lines of finite dimensional spaces and complexes and define a sign refined version
of the isomorphism between the determinant line of a complex and the determinant line of its cohomology similar to the one introduced by Turaev,
\cite{Turaev86} (see also \cite{Turaev90}, \cite{Turaev01}, and  \cite{FarberTuraev00}). We also discuss some properties of this isomorphism.

Let $\bfk$ be a field of characteristic zero.
\subsection{Determinant lines}\Label{SS:det}
Let $V$ be a $\bfk$-vector space of dimension $\dim{}V= n$. The {\em determinant line} of $V$ is the line $\Det(V):= \Lam^nV$, where $\Lam^nV$
denotes the $n$-th exterior power of $V$. By definition, we set $\Det(0):= \bfk$.

More generally, if $V^\b=V^0\oplus{}V^1\oplus\cdots\oplus{}V^d$ is a graded $\bfk$-vector space, we define the determinant line of $V^\b$ by the
formula
\eq{weighteddetb}
    \Det(V^\b) \ := \ \bigotimes_{j=0}^d\,\Det(V^j)^{(-1)^j},
\end{equation}
where for a $\bfk$-line $L$ we denote by $L^{-1}:= \Hom_\bfk(L,\bfk)$ the dual line.

If $L$ is a $\bfk$-line and $l\in L$ is a non-zero element, we denote by $l^{-1}\in L^{-1}$ the unique $\bfk$-linear map $L\to \bfk$ such that
$l^{-1}(l)= 1$.

\subsection{The determinant line of a finite dimensional complex}\Label{SS:FinDimCompl}
Let
\eq{Cp}
    \begin{CD}
       (C^\b,\pa):\quad  0 \ \to C^0 @>{\pa}>> C^1 @>{\pa}>>\cdots @>{\pa}>> C^d \ \to \ 0
    \end{CD}
\end{equation}
be a  complex of finite dimensional $\bfk$-vector spaces. We call the integer $d$ the {\em length} of the complex $(C^\b,\pa)$ and we denote by
$H^\b(\pa)=\bigoplus_{i=0}^d H^{i}(\pa)$ the cohomology of $(C^{\b},\pa)$. Set
\eq{DetCH}
    \Det(C^\b) \ := \ \bigotimes_{j=0}^d\,\Det(C^j)^{(-1)^j}, \qquad \Det(H^\b(\pa))\ := \ \bigotimes_{j=0}^d\,\Det(H^j(\pa))^{(-1)^j}.
\end{equation}

\subsection{The determinant line of a direct sum}\Label{SS:directsum}
For two finite dimensional $\bfk$-vector spaces $V$ and $W$ we define the canonical  {\em fusion} isomorphism
\eq{directsum}
    \mu_{V,W}:\, \Det(V)\,\otimes\,\Det(W)\, \longrightarrow \ \Det(V\oplus W)
\end{equation}
by the formula
\eq{directsum2}
    \mu_{V,W}:\, (v_1\wedge v_2\wedge\cdots\wedge v_k)\,\otimes\,(w_1\wedge w_2\wedge\cdots\wedge w_l)
    \ \mapsto \ v_1\wedge v_2\wedge\cdots\wedge v_k\wedge w_1\wedge w_2\wedge\cdots\wedge w_l,
\end{equation}
where $k=\dim{V},\ l=\dim W,\ v_j\in V,\ w_j\in W$. Clearly,
\eq{muVW-muWV}
    \mu_{V,W}(v\otimes w)  \ = \ (-1)^{\dim V\cdot\,\dim W}\, \mu_{W,V}(w\otimes v), \qquad v\in \Det(V), w\in \Det(W).
\end{equation}

By a slight abuse of notation, we denote by
\eq{directsum3}
    \mu_{V,W}^{-1}:\, \Det(V)^{-1}\,\otimes\,\Det(W)^{-1}\, \longrightarrow \ \Det(V\oplus W)^{-1}
\end{equation}
the transpose of the inverse of $\mu_{V,W}$. It than follows that, for any $v\in \Det(V)$ and $w\in \Det(W)$,
\eq{directsum4}
    \mu_{V,W}^{-1}(v^{-1}\otimes w^{-1}) \ = \ \big(\, \mu_{V,W}(v\otimes w)\,\big)^{-1}.
\end{equation}

Similarly, if $V_1\nek V_r$ are finite dimensional $\bfk$-vector spaces, we define an isomorphism
\eq{directsumr}
    \mu_{V_1\nek V_r}:\, \Det(V_1)\,\otimes\cdots\,\otimes\,\Det(V_r) \ \longrightarrow\ \Det(V_1\oplus\cdots\oplus\,V_r).
\end{equation}
One easily checks that, for every $j\in 1\nek r-1$,
\eq{V1Vr}
    \mu_{V_1\nek V_r} \ = \ \mu_{V_1\nek V_{j-1},V_j\oplus V_{j+1}, V_{j+2}\nek V_r}\circ
    \big(\, 1\otimes\cdots\otimes1\otimes \mu_{V_j,V_{j+1}}\otimes1\otimes\cdots\otimes1\,\big)
\end{equation}

\subsection{The isomorphism between the determinant line of a complex and the determinant line of its cohomology}\Label{SS:isomorphismn}
Fix a direct sum decomposition
\eq{decomposition}
    C^j \ = \ B^j\,\oplus\, H^j\,\oplus\, A^j, \qquad j=0\nek d,
\end{equation}
such that $B^j\oplus{}H^j= (\Ker\pa)\cap C^j$ and $B^j= \pa(C^{j-1})= \pa(A^{j-1})$, for all $j$. Note that $A^d=\{0\}$. Set $A^{-1}=\{0\}$. Then
$H^j$ is naturally isomorphic to the cohomology $H^j(\pa)$ and $\pa$ defines an isomorphism $\pa:A^j\to B^{j+1}$.

For each $j=0\nek d$, fix $c_j\in \Det(C^j)$ and $a_j\in \Det(A^j)$. Let $\pa(a_j)\in \Det(B^{j+1})$ denote the image of $a_j$ under the map
$\Det(A^j)\to \Det(B^{j+1})$ induced by the isomorphism $\pa:A^j\to B^{j+1}$. Then, for each $j=0\nek d$, there is a unique element $h_j\in
\Det(H^j)$ such that
\eq{c=ahb}
    c_j \ = \ \mu_{B^j,H^j,A^j}\big(\,\pa(a_{j-1})\otimes h_j\otimes a_j\,\big).
\end{equation}


Define the isomorphism
\eq{isomorphism}
     \phi_{C^\b} \ = \ \phi_{(C^\b,\pa)}:\, \Det(C^\b) \ \longrightarrow \ \Det(H^\b(\pa)) \ \simeq \ \Det(H^\b),
\end{equation}
by the formula
\eq{isomorphism2}
     \phi_{C^\b}:\, c_0\,\otimes\, c_1^{-1}\,\otimes\,\cdots\, \otimes\, c_d^{(-1)^d} \ \mapsto \
    (-1)^{\calN(C^\b)}\, h_0\,\otimes\, h_1^{-1}\,\otimes\,\cdots\, \otimes\, h_d^{(-1)^d},
\end{equation}
where
\eq{N(C)}
    \calN(C^\b) \ := \ \frac12\ \sum_{j=0}^d\, \dim A^j\cdot\big(\,\dim A^j+(-1)^{j+1}\,\big).
\end{equation}
One easily checks that $ \phi_{C^\b}$ is independent of the choices of $c_j$ and $a_j$.

\rem{N(C)}
a. \ We have
\eq{N(C)2}
   \sum_{k=0}^{j}\,(-1)^{k}\,\dim C^k \ + \ (-1)^{j+1}\,\dim A^j\ = \ \sum_{k=0}^j\,(-1)^{k}\,\dim H^k,
\end{equation}
since both sides of this equality are equal to the Euler characteristic of the complex
\[
 \begin{CD}
    0@>>> C^0 @>\pa>> C^{1}@>\pa>> \cdots @>\pa>> C^j @>\pa>> \pa(A^j) @>>> 0.
 \end{CD}
\]
Hence, $\calN(C^\b)$ can be expressed exclusively in terms of the dimensions of the spaces  $C^j$ and $H^j(\pa)$.

b.\ The isomorphism $ \phi_{C^\b}$ is a sign refined version of the standard construction, cf. \cite{Milnor66}. The idea to introduce a sign factor
in the definition of $ \phi_{C^\b}$ is due to Turaev \cite{Turaev86}. It allows to obtain various compatibility properties, cf., for example,
\refl{fusionFT} and \refp{reftor-grdetfd} below. Our sign is slightly different from \cite{Turaev86} but is consistent with \cite{Nicolaescu03RT}. We
refer the reader to \cite{Deligne85} and \cite{Nicolaescu03RT} for the motivation of this choice of sign, based on the theory of weighted determinant
lines.
\erem
\subsection{The fusion isomorphism for graded vector spaces}\Label{SS:fusiongr}
Let  $V^\b=V^0\oplus{}V^1\oplus\cdots\oplus{}V^d$ and  $W^\b= W^0\oplus{}W^1\oplus\cdots\oplus{}W^d$ be finite-dimensional graded $\bfk$-vector
spaces. The {\em fusion isomorphism}
\eq{fusiongr}
    \mu_{V^\b,W^\b}:\, \Det(V^\b)\,\otimes\,\Det(W^\b) \ \longrightarrow\ \Det(V^\b\oplus W^\b),
\end{equation}
is defined by the formula
\eq{fusiongr2}
    \mu_{V^\b,W^\b} \ := \ (-1)^{\calM(V^\b,W^\b)}\,\bigotimes_{q=0}^d\,\mu_{V^q,W^q}^{(-1)^q},
\end{equation}
where $\mu_{V^j,W^j}^{+1}= \mu_{V^j,W^j}$ and $\mu_{V^j,W^j}^{-1}$ are defined in \refss{directsum}, and
\eq{M(V,W)}
    \calM(V^\b,W^\b) \ := \
    \sum_{0\le k<j\le d}\, \dim V^j\cdot\dim  W^k.
\end{equation}

The following lemma is a precise analogue of  Lemma~2.4 of \cite{FarberTuraev00}. 
\lem{fusionFT}
Let $(C^\b,\pa)$ and $(\tilC^\b,\tilpa)$ be length $d$ complexes of finite dimensional $\bfk$-vector spaces. Then the following diagram commutes:
\eq{fusinonFT}
    \begin{CD}
        &\Det(C^\b)\,\otimes\,\Det(\tilC^\b) & @>{ \phi_{C^\b}\otimes \phi_{\tilC^\b}}>>& \Det(H^\b(\pa))\,\otimes
        \Det(H^\b(\tilpa))&\\
        &@V{\mu_{C^\b,\tilC^\b}}VV&& @VV{\mu_{H^\b(\pa),H^\b(\tilpa)}}V&\\
        &\Det(C^\b\oplus \tilC^\b)& @>{\phi_{C^\b\oplus \tilC^\b}}>>& \Det\big(\,H^\b(\pa\oplus \tilpa)\,\big) \ \simeq \
        \Det\big(\,H^\b(\pa))\oplus H^\b(\tilpa)\,\big)&
    \end{CD}
\end{equation}
\elem
\prf
As in \refe{decomposition}, write
\eq{decomp2}
  C^j\ = \ B^j\oplus H^j\oplus A^j,\qquad \tilC^j \ = \ \tilB^j\oplus\tilH^j\oplus\tilA^j.
\end{equation}
For each $j= 0\nek d$, choose
\[
 \begin{aligned}
    &c_j\ \in\ \Det(C^j), \quad &a_j\ \in\ \Det(A^j),\quad &h_j\ \in\ \Det(H^j),\\
    &\tilc_j\ \in\ \Det(\tilC^j), \quad &\tila_j\ \in\ \Det(\tilA^j), \quad &\tilh_j\ \in\ \Det(\tilH^j),
 \end{aligned}
\]
such that
\eq{c=ahb2}
    c_j \ = \ \mu_{B^j,H^j,A^j}\big(\,\pa(a_{j-1})\otimes h_j\otimes a_j\,\big), \qquad
    \tilc_j \ = \ \mu_{\tilB^j,\tilH^j,\tilA^j}\big(\,\tilpa(\tila_{j-1})\otimes \tilh_j\otimes\tila_j\,\big).
\end{equation}
Set
\[
    \hatC^j\ =\ C^j\oplus \tilC^j, \quad \hatH^j\ =\ H^j\oplus \tilH^j, \quad \hatA^j\ =\ A^j\oplus \tilA^j,\quad  \hatB^j\ =\ B^j\oplus\tilB^j.
\]
Also denote $\hatpa= \pa\oplus\tilpa$. Further, set $\hatc_j= \mu_{C^j,\tilC^j}(c_j\otimes\tilc_j)$ and $\hata_j=
\mu_{A^j,\tilA^j}(a_j\otimes\tila_j)$. Then, for all $j= 0\nek d$, the unique element $\hath_j\in \Det(\hatH^j)$, satisfying
\eq{hatc=hath}
    \hatc_j \ = \ \mu_{\hatB^j,\hatH^j,\hatA^j}\big(\,\hatpa(\hata_{j-1})\otimes \hath_j\otimes \hata_j\,\big),
\end{equation}
is given by
\eq{hathj}
    \hath_j \ = \ (-1)^{\dim A^j\cdot \dim\tilA^{j-1}+\dim H^j\cdot\dim\tilA^{j-1} + \dim A^j\cdot\dim \tilH^j}\,
    \mu_{H^j,\tilH^j}(h_j\otimes\tilh_j).
\end{equation}

Set $c := \otimes_{j=0}^d\,c_j^{(-1)^j}$ and define $\tilc$ in a similar way. Then, by definitions \refe{fusiongr2} and \refe{hatc=hath},
\eq{mu(cc)}
    \mu_{C^\b,\tilC^\b}(c\otimes\tilc) \ = \ (-1)^{\calM(C^\b,\tilC^\b)}\, \otimes_{j=0}^d\,\mu_{C^j,\tilC^j}(c_j\otimes\tilc_j)^{(-1)^j}
    \ = \ (-1)^{\calM(C^\b,\tilC^\b)}\, \otimes_{j=0}^d\, \hatc_j^{(-1)^j}.
\end{equation}

From \refe{hatc=hath} and \refe{mu(cc)} we conclude that
\eq{phimu}
    \phi_{C^\b\oplus \tilC^\b}\circ \mu_{C^\b,\tilC^\b}(c\otimes\tilc) \ = \ (-1)^{\calK(C^\b,\tilC^\b)}\,
    \otimes_{j=0}^d\,\mu_{H^j,\tilH^j}(h_j\otimes \tilh_j)^{(-1)^j},
\end{equation}
where
\meq{K(C)}
    \calK(C^\b,\tilC^\b) \ = \ \calN(C^\b\oplus\tilC^\b) \ + \ \calM(C^\b,\tilC^\b) \\ + \
    \sum_{j=0}^d\,\Big[\,\dim A^j\cdot \dim\tilA^{j-1}+\dim H^j\cdot\dim\tilA^{j-1} + \dim A^j\cdot\dim \tilH^j\,\Big]
\end{multline}

Since, clearly,
\meq{muphi}
    \mu_{H^\b(\pa),H^\b(\tilpa)}\circ\big(\, \phi_{C^\b}\otimes \phi_{\tilC^\b}\,\big)(c\otimes\tilc) \\ = \
    (-1)^{\calN(C^\b)+\calN(\tilC^\b)+\calM(H^\b,\tilH^\b)}\, \otimes_{j=0}^d\,\mu_{H^j,\tilH^j}(h_j\otimes \tilh_j)^{(-1)^j},
\end{multline}
to prove the commutativity of the diagram \refe{fusinonFT} it remains to show that, $\MOD 2$,
\meq{signfusino}
    \calN(C^\b\oplus\tilC^\b) \ + \ \calN(C^\b) \ + \ \calN(\tilC^\b) \ +\ \calM(H^\b,\tilH^\b) \ + \ \calM(C^\b,\tilC^\b)
    \\ \equiv  \
        \sum_{j=0}^d\,\Big[\,\dim A^j\cdot \dim\tilA^{j-1}+\dim H^j\cdot\dim\tilA^{j-1} + \dim A^j\cdot\dim \tilH^j\,\Big].
\end{multline}

Using the identity
\eq{x(x+1)}
    \frac{(x+y)\,(x+y+(-1)^j)}2 \ - \ \frac{x\,(x+(-1)^j)}2 \ - \ \frac{y\,(y+(-1)^j)}2 \ = \ xy, \qquad\qquad x,y\in \bfk, \ j\in \ZZ_{\ge0},
\end{equation}
we obtain
\eq{N-N-N}
    \calN(C^\b\oplus\tilC^\b) \ - \ \calN(C^\b) \ - \ \calN(\tilC^\b) \ =\ \sum_{j=0}^{d}\,\dim A^j\cdot \dim \tilA^j.
\end{equation}

On the other hand, using \refe{decomp2} and the equalities $\dim\tilA^{j-1}= \dim\tilB^{j}$ we see that the following equality holds modulo 2
\meq{2terms}
    \sum_{j=0}^d\,\Big[\,\dim A^j\cdot \dim\tilA^{j-1} + \dim A^j\cdot\dim \tilH^j\,\Big]
    \ = \
    \sum_{j=0}^d\,\dim A^j\cdot\big(\, \dim \tilA^{j-1}+\dim \tilH^j\,\big)
    \\ \equiv \
    \sum_{j=0}^d\,\dim A^j\cdot\big(\, \dim \tilA^{j}+\dim \tilC^j\,\big) \ = \ \sum_{j=0}^d\,\dim A^j\cdot\dim\tilA^j
    \ + \ \sum_{j=0}^d\,\dim A^j\cdot\dim\tilC^j.
\end{multline}

By \refe{N(C)2},
\[
    \dim A^j \ \equiv \ \sum_{k=0}^j\,\big(\,\dim H^k+ \dim C^k\,\big).
\]
A similar equality holds for $\dim\tilA^j$. Hence, we get from \refe{2terms}
\meq{2terms2}
    \sum_{j=0}^d\,\Big[\,\dim A^j\cdot \dim\tilA^{j-1} + \dim A^j\cdot\dim \tilH^j\,\Big]
    \\ \equiv \  \sum_{j=0}^d\,\dim A^j\cdot\dim\tilA^j \ + \ \sum_{0\le k\le j\le d}\, \dim C^k\cdot\dim\tilC^j
    \ + \ \sum_{0\le k\le j\le d}\, \dim H^k\cdot\dim\tilC^j.
\end{multline}
Similarly,
\eq{1term}
    \sum_{j=0}^d\,\dim H^j\cdot\dim\tilA^{j-1} \ = \ \sum_{0\le k<j\le d}\,\dim H^j\cdot\dim\tilC^k \ + \  \sum_{0\le k<j\le d}\,\dim
    H^j\cdot\dim\tilH^k.
\end{equation}
Combining \refe{M(V,W)}, \refe{2terms2}, and \refe{1term} we obtain that modulo 2
\meq{3term+M}
    \sum_{j=0}^d\,\Big[\,\dim A^j\cdot \dim\tilA^{j-1}+\dim H^j\cdot\dim\tilA^{j-1} + \dim A^j\cdot\dim \tilH^j\,\Big]
    \ + \ \calM(C^\b,\tilC^\b) \ + \ \calM(H^\b,\tilH^\b)
    \\ \equiv \   \sum_{j=0}^d\,\dim A^j\cdot\dim\tilA^j \ + \  \sum_{j,k=0}^d\,\dim C^k\dim\tilC^j
     \ + \  \sum_{j,k=0}^d\,\dim H^k\dim\tilC^j
    \\ \equiv \ \sum_{j=0}^d\,\dim A^j\cdot\dim\tilA^j \ + \ \Big(\,\sum_{k=0}^d\,\dim C^k\,\Big)\cdot\Big(\,\sum_{j=0}^d\,\dim \tilC^j\,\Big)
     \ + \   \Big(\,\sum_{k=0}^d\,\dim H^k\,\Big)\cdot\Big(\,\sum_{j=0}^d\,\dim \tilC^j\,\Big).
\end{multline}
Both $\sum_{k=0}^d\,\dim C^k$ and $\sum_{k=0}^d\,\dim H^k$ are equivalent modulo 2 to the Euler characteristic of the complex $(C^\b,\pa)$. Hence, we
conclude that the left hand side of \refe{3term+M} is equivalent modulo 2 to $$\sum_{j=0}^d\,\dim A^j\cdot\dim\tilA^j.$$ Combining this with
\refe{N-N-N}, we obtain \refe{signfusino}.
\eprf

\section{The Determinant Line of the Dual Complex}\Label{S:dualcomplex}

In this section we introduce the dual of a complex and, for the case when the length of the complex is odd, construct a natural isomorphism between
the determinant lines of a complex and that of its dual. We also show that this isomorphism is compatible with the canonical isomorphism
\refe{isomorphism}.

Throughout the section, $\bfk$ is a field of characteristic zero endowed with an involutive automorphism
\[
    \tau:\bfk\to \bfk.
\]
The main examples are $\bfk= \CC$ with $\tau$ being the complex conjugation and $\bfk= \RR$ with $\tau$ being the identity map.

\subsection{The determinant line of the $\tau$-dual space}\Label{SS:dualspace}
If $V,W$ are $\bfk$-vector spaces, a map $f:V\to W$ is said to be {\em $\tau$-linear} if
\[
    f(x_1v_1+x_2v_2)\ = \ \tau(x_1)\,v_1 \ + \ \tau(x_2)\,v_2, \qquad \text{for any}\quad  v_1,v_2\in V,\ x_1,x_2\in \bfk.
\]

Let $V$ be an $n$-dimensional $\bfk$-vector space. The linear space $V^*= \tV$ of all $\tau$-linear maps $V\to \bfk$  is called the {\em $\tau$-dual
space to $V$}. There is a natural $\tau$-linear isomorphism
\eq{u(V*)}
    \alp_{V}:\, \Det(V^*) \ \longrightarrow \ \Det(V)^{-1},
\end{equation}
defined by the formula
\eq{u(V*)2}
    \big(\,\alp_{V}(v^1\wedge\cdots \wedge v^n)\,\big)(v_1\wedge\cdots\wedge v_n)
    \ = \ \sum_{\sig}\,(-1)^{|\sig|}\, \tau\big(v^1(v_{\sig(1)})\big)\cdot \tau\big(v^2(v_{\sig(2)})\big)\cdots \tau\big(v^n(v_{\sig(n)})\big),
\end{equation}
where $v_1\nek v_n\in V, \ v^1\nek v^n\in V^*$, and the sum is taken over all permutations $\sig$ of $\{1\nek n\}$. Similarly, we define the
$\tau$-linear map
\eq{u(V)}
    \bet_V:\, \Det(V) \ \longrightarrow \ \Det(V^*)^{-1},
\end{equation}
defined by the formula
\meq{u(V)2}
    \big(\,\bet_V(v_1\wedge\cdots \wedge v_n)\,\big)(v^1\wedge\cdots\wedge v^n)
    \\ = \ (-1)^n\cdot\sum_{\sig}\,(-1)^{|\sig|}\, \tau\big(v^1(v_{\sig(1)})\big)\cdot \tau\big(v^2(v_{\sig(2)})\big)\cdots \tau\big(v^n(v_{\sig(n)})\big).
\end{multline}

\rem{sign--=+}
The sign factor $(-1)^{n}$ in \refe{u(V)2} simplifies the statements of various compatibility relations with the fusion isomorphism \refe{directsum},
cf. below. It is motivated by the fact that, in \refe{u(V)2}, we interchange $v_1\wedge\cdots \wedge v_n$ and $v^1\wedge\cdots\wedge v^n$, which both
are forms of degree\/ $n$.
\erem

{}Formulae \refe{u(V*)2} and \refe{u(V)2} can be simplified by choosing an appropriate basis. Let $e_1\nek e_n$ be a basis of\/ $V$. Denote by
$e^1\nek e^n$ the dual basis of\/ $V^*$, i.e., the unique set of elements of $V^*$ such that $e^j(e_i)= \del^j_i$ for all $i,j=1\nek n$. Then
\begin{gather}
    \big(\,\alp_V(e^1\wedge\cdots\wedge e^n)\,\big)(e_1\wedge\cdots \wedge e_n) \ = \
    \tau\big(e^1(e_1)\big)\cdot \tau\big(e^2(e_2)\big)\cdots \tau\big(e^n(e_n)\big); \Label{E:u(V*)3}\\
    \big(\,\bet_V(e_1\wedge\cdots \wedge e_n)\,\big)(e^1\wedge\cdots\wedge e^n)
    \ = \ (-1)^n\cdot \tau\big(e^1(e_1)\big)\cdot \tau\big(e^2(e_2)\big)\cdots \tau\big(e^n(e_n)\big).\Label{E:u(V)3}
\end{gather}
Recall from \refss{det} that for a non-zero element $v\in \Det(V)$ we denote by $v^{-1}$ the unique element of $\big(\Det(V)\big)^{-1}$ such that
$v^{-1}(v)= 1$. It follows from \refe{u(V*)3} that
\eq{alpV-1}
    \alp_V^{-1}\big(\,(e_1\wedge\cdots \wedge e_n)^{-1}\,\big) \ = \ e^1\wedge \cdots \wedge e^n.
\end{equation}

Using \refe{u(V)3} and \refe{alpV-1} we conclude that for any $v\in \Det(V)$
\eq{u-1(v)-1}
    \big(\,\alp_V^{-1}(v^{-1})\,\big)^{-1} \ = \ (-1)^{\dim V}\, \bet_V(v).
\end{equation}

Let $V$ and $W$ be $\bfk$-vector spaces. From \refe{directsum2}, \refe{u(V*)3}, and \refe{u(V)3}, we obtain
\eq{mu-u}
    \big(\,\mu_{V,W}(v\otimes w)\,\big)^{-1} \ = \ \alp_{V\oplus W}\circ\mu_{V^*,W^*}\big(\,\alp_V^{-1}(v^{-1})\otimes \alp_W^{-1}(w^{-1})\,\big),
\end{equation}
for any $v\in \Det(V),\ w\in \Det(W)$.

\subsection{The $\tau$-adjoint map}\Label{SS:tauadj}
Let $T:V\to W$ be a $\bfk$-linear map. The {\em $\tau$-adjoint of\/ $T$} is the linear map
\[
    T^*:\, W^* \ \longrightarrow \ V^*
\]
such that
\eq{tauadj2}
    (T^*w^*)(v) \ = \ w^*(Tv), \qquad \text{for all}\quad v\in V,\ w^*\in W^*.
\end{equation}

If $\dim V= \dim W$ then $T$ and $T^*$ induce $\bfk$-linear maps $\Det(V)\to \Det(W)$ and $\Det(W^*)\to \Det(V^*)$, which, by a slight abuse of
notation, we also denote by $T$ and $T^*$ respectively. If $T$ is bijective then, for any non-zero $v\in \Det(V)$, we have
\eq{tauadj-1}
    T^*\alp_W^{-1}\big(\,(Tv)^{-1}\,\big)\ =\ \alp_V^{-1}(v^{-1}).
\end{equation}

\subsection{The $\tau$-dual graded space}\Label{SS:dualgraded}
Let now $V^\b= V^0\oplus V^1\oplus \cdots\oplus V^d$  be a finite dimensional graded $\bfk$-vector space. We define the {\em ($\tau$-)dual graded
space}\/ $\hatV= \hatV^0\oplus\hatV^1\oplus\cdots\oplus\hatV^d$ by
\eq{hatVb}\notag
    \hatV^j \ := \ (V^{d-j})^*, \qquad j=0\nek d.
\end{equation}
{Assume now that the number $d=2r-1$ is odd}. Then \refe{u(V*)} and \refe{u(V)} induce a $\tau$-linear isomorphism
\eq{u(Vb)}
    \alp_{V^\b}:\, \Det(V^\b) \ \longrightarrow \ \Det(\hatV^\b),
\end{equation}
defined by the formula
\eq{u(Vb)2}
    \alp_{V^\b}\big(\, v_0\otimes (v_{1})^{-1}\otimes \cdots \otimes(v_d)^{-1}\,\big) \ = \
    (-1)^{\calM(V^\b)}\cdot \alp_{V^d}^{-1}(v_d^{-1})\otimes \bet_{V^{d-1}}(v_{d-1})\otimes\cdots\otimes \bet_{V^0}(v_0),
\end{equation}
where $v_j\in \Det(V^j)$ $(j=0\nek d)$ and
\eq{calF}
    \calM(V^\b) \ = \ \calM(V^\b,V^\b) \ = \  \sum_{0\le j<k\le d}\, \dim V^j\cdot \dim V^k,
\end{equation}
cf. \refe{M(V,W)}. We refer to \cite{Nicolaescu03RT} for the motivation of the choice of the sign in \refe{u(Vb)2}.

\subsection{The dual complex}\Label{SS:dualcomplex}
Consider the complex \refe{Cp} of finite dimensional $\bfk$-vector spaces. The {\em dual complex} is the complex
\eq{hatCp}
    \begin{CD}
       (\hatC^\b,\pa):\quad  0 \ \to \hatC^0 @>{\pa^*}>> \hatC^1 @>{\pa^*}>>\cdots @>{\pa^*}>> \hatC^d \ \to \ 0,
    \end{CD}
\end{equation}
where $\hatC^j= (C^{d-j})^*$ and $\pa^*$ is the $\tau$-adjoint of $\pa$. Then the cohomology $H^j(\pa^*)$ of $\hatC^\b$ is naturally isomorphic to
the $\tau$-dual space to $H^{d-j}(\pa)$ $(j=0\nek d)$. Hence, if the length $d$ of the complex $C^\b$ is odd, then, by \refe{u(Vb)}, we obtain
$\tau$-linear isomorphisms
\eq{uC}
    \alp_{C^\b}:\, \Det(C^\b) \ \longrightarrow \ \Det(\hatC^\b), \qquad
    \alp_{H^\b(\pa)}:\, \Det(H^\b(\pa))\ \longrightarrow \ \Det(H^\b(\pa^*)).
\end{equation}

\lem{hatC-C}
Let $(C^\b,\pa)$ be a complex of finite dimensional $\bfk$-vector spaces and assume that its length $d=2r-1$ is odd. Then the following diagram
commutes
\eq{hatC-C}
    \begin{CD}
        &\Det(C^\b) & @>{ \phi_{C^\b}}>>& \Det\big(H^\b(\pa)\big)&\\
        &@V{\alp_{C^\b}}VV&& @VV{\alp_{H^\b(\pa)}}V&\\
        &\Det(\hatC^\b)& @>{\phi_{\hatC^\b}}>>& \Det\big(H^\b(\pa^*)\big)&
    \end{CD},
\end{equation}
where the isomorphisms $\phi_{C^\b}$ and $\phi_{\hatC^\b}$ are as in \refe{isomorphism}.
\elem
\prf
We shall use the notation of \refss{isomorphismn}. For $j=0\nek d$, set
\eq{hatAhatBhatH}
    \hatA^j\ := \ (B^{d-j})^*, \quad \hatB^j\ := \ (A^{d-j})^*, \quad \hatH^j\ := \ (H^{d-j})^*
\end{equation}
and identify these spaces with subspaces of $\hatC^j$ in a natural way. Then $\pa^*(\hatA^j)= \hatB^{j+1}$.

Let $c_j, a_j, h_j$ ($j=0\nek d$) be as in \refe{c=ahb}. For each $j=0\nek d$, set
\begin{align}
    \hatc_j \ &= \ \alp_{C^{d-j}}^{-1}(c_{d-j}^{-1}),\Label{E:hatc}\\
    \hata_j\  &= \ \alp_{B^{d-j}}^{-1}\big((\pa a_{d-j-1})^{-1}\big),\Label{E:hata}\\
    \hath_j\ &= \ \alp_{H^{d-j}}^{-1}(h_{d-j}^{-1}).\Label{E:hath}
\end{align}
Then, from the equality \refe{u-1(v)-1}, we obtain,
\eq{bet-alp-1}\notag
 \begin{aligned}
    \bet_{C^{d-j}}(c_{d-j}) \ &= \ (-1)^{\dim C^j}\cdot \hatc_j^{-1},\\
    \bet_{H^{d-j}(\pa)}(h_{d-j}) \ &= \ (-1)^{\dim H^j(\pa)}\cdot \hath_j^{-1}.
 \end{aligned}
\end{equation}
Hence, from the definition \refe{u(Vb)2}, we get
\begin{align}
    \alp_{C^\b}\big(\, c_0\otimes c_1^{-1}\otimes\cdots\otimes c_d^{-1}\,\big) \ &= \ (-1)^{\calM(C^\b)+\sum_{p=0}^{r-1}\dim C^{2p}}\cdot
    \hatc_0\otimes \hatc_1^{-1}\otimes\cdots\otimes \hatc_d^{-1},\Label{E:alpc=hatc}\\
    \alp_{H^\b(\pa)}\big(\, h_0\otimes h_1^{-1}\otimes\cdots\otimes h_d^{-1}\,\big) \ &= \ (-1)^{\calM(H^\b(\pa))+\sum_{p=0}^{r-1}\dim H^{2p}(\pa)}\cdot
    \hath_0\otimes \hath_1^{-1}\otimes\cdots\otimes \hath_d^{-1}.\Label{E:alph=hath}
\end{align}

From the identity \refe{tauadj-1} and the definition \refe{hata} of $\hata_j$, we get
\eq{pa*hata}\notag
    \pa^*(\hata_{j-1}) \ = \ \alp_{A^{d-j}}^{-1}(a_{d-j}^{-1}), \qquad j=1\nek d.
\end{equation}
Hence, from \refe{c=ahb} and \refe{mu-u}, we obtain
\eq{hatcj=ahb}\notag
 \begin{aligned}
    \hatc_j \ &= \ \alp_{C^{d-j}}^{-1}\Big(\,\mu_{B^{d-j},H^{d-j},A^{d-j}}\big(\,\pa(a_{d-j-1})\otimes h_{d-j}\otimes a_{d-j}\,\big)\,\Big)^{-1}
    \\ &= \  \mu_{\hatA^j,\hatH^j,\hatB^j}\big(\,\hata_j\otimes \hath_j\otimes \pa^*(\hata_{j-1})\,\big).
 \end{aligned}
\end{equation}
Using \refe{muVW-muWV}, we now conclude that
\eq{hatcj=ahb2}
    \hatc_j \ = \ (-1)^{\calG_j}\cdot
      \mu_{\hatB^j,\hatH^j,\hatA^j}\big(\,\pa^*(\hata_{j-1})\otimes \hath_j\otimes \hata_{j-1}\,\big),
\end{equation}
where
\eq{calGj}
    \calG_j \ = \ \dim\hatA^j\cdot\dim\hatH^j+\dim\hatA^j\cdot\dim\hatA^{j-1}+\dim\hatA^{j-1}\cdot\dim\hatH^j.
\end{equation}

Thus, from \refe{isomorphism2}, we obtain
\eq{phihatCc}
    \phi_{\hatC^\b}\big(\, \hatc_0\otimes \hatc_1^{-1}\otimes\cdots\otimes\hatc_d^{-1}\,\big) \ = \
    (-1)^{\calN(\hatC^\b)+\sum_{j=0}^d\calG_j}\cdot
    \hath_0\otimes\hath_1^{-1}\otimes\cdots\otimes\hath_d^{-1}.
\end{equation}
Hence, by \refe{u(Vb)2} and \refe{alpc=hatc},
\meq{uphihatCc}
    \phi_{\hatC^\b}\circ\alp_{C^\b}\big(\, c_0\otimes c_1^{-1}\otimes\cdots\otimes c_d^{-1}\,\big) \\ = \
    (-1)^{\calM(C^\b)+\calN(\hatC^\b)+\sum_{j=0}^d\calG_j+\sum_{p=0}^{r-1}\dim C^{2p}}\cdot
    \hath_0\otimes \hath_1^{-1}\otimes\cdots\otimes \hath_d^{-1}.
\end{multline}
From \refe{u(Vb)2} and \refe{alph=hath}, we get
\meq{phiuhatc}
    \alp_{H^\b(\pa)}\circ \phi_{C^\b}\big(\, c_0\otimes c_1^{-1}\otimes\cdots\otimes c_d^{-1}\,\big)
     \\ = \ (-1)^{\calM(H^\b(\pa))+\calN(C^\b)+\sum_{p=0}^{r-1}\dim H^{2p}(\pa)}\cdot
    \hath_0\otimes \hath_1^{-1}\otimes\cdots\otimes \hath_d^{-1}.
\end{multline}
From \refe{uphihatCc} and \refe{phiuhatc}, we conclude that to prove \refe{hatC-C} it remains to show that, modulo 2,
\eq{FGN}
    \calM(C^\b)+\calN(\hatC^\b)+\sum_{j=0}^d\calG_j +\sum_{p=0}^{r-1}\dim C^{2p}\
     \equiv \ \calM(H^\b(\pa))+\calN(C^\b) + \ \sum_{p=0}^{r-1}\dim H^{2p}(\pa).
\end{equation}
Using the equality
\[
    \dim \hatA^j \ = \ \dim B^{d-j} \ = \ \dim A^{d-j-1},
\]
we easily see that $\calN(\hatC^\b)= \calN(C^\b)$. In addition, note that if we set $\tilC^\b= C^\b$ in \refe{signfusino}, then the right hand side
of \refe{signfusino} is equal to  $\sum_{j=0}^d\calG_j$. Hence, \refe{signfusino} and \refe{calF} imply that \refe{FGN} is equivalent to
\eq{N(C+C)+C+H}
    \calN(C^\b\oplus C^\b) \ + \ \sum_{p=0}^{r-1}\,\big(\,\dim C^{2p}- \dim H^{2p}(\pa)\,\big)\equiv \ 0\qquad \MOD\ 2.
\end{equation}
By \refe{N(C)},
\eq{N(C+C)}
    \calN(C^\b\oplus C^\b) \ \equiv \ \sum_{j=0}^d\, \dim A^j \qquad \MOD\ 2.
\end{equation}
From \refe{decomposition} and the equality $\dim B^j= \dim A^{j-1}$ we conclude that
\[
    \dim C^{2p}\ -\ \dim H^{2p}(\pa)\ = \ \dim A^{2p} \ + \ \dim A^{2p-1},
\]
and, hence,
\eq{sumC2p-H2p}
    \sum_{p=0}^{r-1}\,\big(\,\dim C^{2p}- \dim H^{2p}(\pa)\,\big) \ = \ \sum_{j=0}^d\, \dim A^j.
\end{equation}
Combining \refe{N(C+C)} and \refe{sumC2p-H2p} we obtain \refe{N(C+C)+C+H}.
\eprf

\section{The Refined Torsion of a Finite Dimensional Complex with a Chirality Operator}\Label{S:FinDimchirality}

On the de Rham complex of a Riemannian manifold acts a canonical involution $\Gam$, called the {\em chirality operator}, cf. \cite[Ch.~3]{BeGeVe}. In
this section we consider a finite dimensional complex with such an operator. We show that the chirality operator defines a canonical element, the
{\em refined torsion}, of the determinant line of the cohomology of this complex and we study some properties of this element.

In this section $\bfk$ is a field of characteristic zero.

\subsection{The refined torsion associated to a chirality operator}\Label{SS:chiral}
Let $d= 2r-1$ be an odd integer and let $(C^\b,\pa)$ be a length $d$ complex of finite dimensional $\bfk$-vector spaces. A {\em chirality operator}
is an involution $\Gam:C^\b\to C^\b$ such that $\Gam(C^j)= C^{d-j}$, $j=0\nek d$. For $c_j\in \Det(C^j)$ $(j=0\nek d)$ we denote by $\Gam{}c_j\in
\Det(C^{d-j})$ the image of $c_j$ under the isomorphism $\Det(C^j)\to \Det(C^{d-j})$ induced by $\Gam$.

Fix non-zero elements $c_j\in \Det(C^j)$, $j=0\nek r-1$, and consider the element
\eq{Gamc}
    c_{{}_\Gam} \ := \ (-1)^{\calR(C^\b)}\cdot
    c_0\otimes c_1^{-1}\otimes \cdots \otimes c_{r-1}^{(-1)^{r-1}}\otimes (\Gam c_{r-1})^{(-1)^r}\otimes (\Gam c_{r-2})^{(-1)^{r-1}}
    \otimes \cdots\otimes (\Gam c_0)^{-1},
\end{equation}
of $\Det(C^\b)$, where

\eq{R(C)}
  \calR(C^\b) \ = \ \frac12\ \sum_{j=0}^{r-1}\, \dim C^{j}\cdot\big(\, \dim C^j+(-1)^{r+j}\,\big).
\end{equation}
It follows from the definition of $c_j^{-1}$ that $c_{{}_\Gam}$ is independent of the choice of $c_j$ ($j=0\nek r-1$).

\rem{Gamc}
Using the isomorphisms $\Gam:C^j\to C^{d-j}$ one can define a natural trace functional $\Tr:\Det(C^\b)\to \bfk$, cf. \cite{Nicolaescu03RT}. The sign
factor $(-1)^{\calR(C^\b)}$ is defined so that the equality $\Tr(c_{{}_\Gam})=1$ holds.
\erem
\defe{refinedtorsion}
The {\em refined torsion} of the pair $(C^\b,\Gam)$ is the element
\eq{refinedtor}
    \rho_{{}_\Gam} \ = \ \rho_{{}_{C^\b,\Gam}} \ := \  \phi_{C^\b}(c_{{}_\Gam}),
\end{equation}
where $ \phi_{C^\b}$ is the canonical map defined in \refss{isomorphismn}.
\edefe

\subsection{The norm of the refined torsion}\Label{SS:norm}
In this subsection we assume that $\bfk=\RR$ or $\CC$. Suppose that the spaces $C^j,\ j=0\nek d$, are endowed with a Euclidean (if $\bfk= \RR$) or a
Hermitian (if $\bfk=\CC$) scalar products $\<\cdot,\cdot\>_j$. These scalar products induce a metric $\|\cdot\|_{\Det(C^\b)}$ on the determinant line
$\Det(C^\b)$. Let $\fdm$ be the metric on the determinant line $\Det(H^\b(\pa))$ such that the canonical isomorphism $\phi_{C^\b}$, defined in
\refe{isomorphism}, is an isometry.

\lem{normrho}
Let $\<\cdot,\cdot\>_j$ be scalar products on $C^j,\ j=0\nek d$,  such that the chirality operator $\Gam$ is self-adjoint. Then
\eq{normrho}
    \|\rho_{{}_\Gam}\|_{\Det(H^\b(\pa))} \ = \ 1.
\end{equation}
\elem
\prf
By definition,
\eq{nrho=nc}
    \|\rho_{{}_\Gam}\|_{\Det(H^\b(\pa))} \ = \ \|c_{{}_\Gam}\|_{\Det(C^\b)}.
\end{equation}
Let $\|\cdot\|_j$ denote the norm on $\Det(C^j)$ induced by $\<\cdot,\cdot\>_j$. Since $\Gam$ is a self-adjoint involution it is also a unitary
operator, i.e., for every $x\in \Det(C^j)$ we have $\|\Gam{x}\|_{d-j}= \|x\|_j$. Hence, from \refe{Gamc} we get $\|c_{{}_\Gam}\|_{\Det(C^\b)}= 1$.
The lemma follows now from \refe{nrho=nc}.
\eprf

The above lemma explains why we call $\rho_{{}_\Gam}$ the refined torsion: the classical combinatorial torsion \cite{Milnor66} is an element $\rho$
of $\Det(H^\b(\pa))$, defined up to a multiplication by $t\in \bfk$ with $|t|=1$, such that $\|\rho\|_{\Det(H^\b(\pa))} \ = \ 1$. The refined torsion
$\rho_{{}_\Gam}$ is a choice of a particular element of $\Det(H^\b(\pa))$ with norm 1.

\subsection{The refined torsion of a direct sum}\Label{SS:rhodirectsum}
\lem{rhodirectsum}
Let $(C^\b,\pa)$ and $(\tilC^\b,\tilpa)$ be length $d=2r-1$ complexes of finite dimensional $\bfk$-vector spaces and let\/ $\Gam:C^\b \to C^\b$,
$\tilGam:\tilC^\b \to \tilC^\b$ be chirality operators. Then $\hatGam:= \Gam\oplus\tilGam:C^\b\oplus\tilC^\b\to C^\b\oplus \tilC^\b$ is a chirality
operator on the direct sum complex $(C^\b\oplus{}\tilC^\b,\pa\oplus{}\tilpa)$ and
\eq{rhodirectsum}
    \rho_{{}_\hatGam} \ = \ \mu_{H^\b(\pa),H^\b(\tilpa)}(\rho_{{}_\Gam}\otimes\rho_{{}_\tilGam}).
\end{equation}
\elem
\prf
Clearly, $\hatGam^2=1$ and $\hatGam(C^j\oplus{}\tilC^j)= C^{d-j}\oplus{}\tilC^{d-j}$. Hence, $\hatGam$ is a chirality operator. By \refl{fusionFT},
to prove \refe{rhodirectsum} it is enough to show that
\eq{cGamdirectsum}
    c_{\hatGam} \ = \ \mu_{C^\b,\tilC^\b}(c_{\Gam}\otimes c_{\tilGam}).
\end{equation}

For each $j=0\nek r-1$, fix non-zero elements $c_j\in \Det(C^j),\ \tilc_j\in \tilC^j$ and set $\hatc_j= \mu_{C^j,\tilC^j}(c_j\otimes \tilc_j)$.
Recall that we denoted the operators induced by $\Gam$ and $\tilGam$ on $\Det(C^\b)$ and $\Det(\tilC^\b)$ by the same letters. Thus,
\[
   \hatGam\hatc_j \ = \ (\Gam\oplus \tilGam)\circ \mu_{C^j,\tilC^j}(c_j\otimes \tilc_j) \ = \
   \mu_{C^{d-j},\tilC^{d-j}}(\Gam c_j\otimes \tilGam\tilc_j).
\]

By \refe{directsum4},
\[
   \mu_{C^j,\tilC^j}^{-1}(c_j^{-1}\otimes \tilc_j^{-1}) \ = \ \big(\,\mu_{C^j,\tilC^j}(c_j\otimes \tilc_j)\,\big)^{-1}.
\]
Hence, it follows from \refe{fusiongr2} and \refe{Gamc} that
\eq{mucc}
  \begin{aligned}
    \mu_{C^\b,\tilC^\b}(c_{\Gam}\otimes c_{\tilGam}) \ &= \ (-1)^{\calM(C^\b,\tilC^\b)+\calR(C^\b)+\calR(\tilC^\b)}\, \times\\
       &\times\, \hatc_0\otimes \hatc_1^{-1}\otimes \cdots \otimes \hatc_{r-1}^{(-1)^{r-1}}\otimes
            \hatGam(\hatc_{r-1})^{(-1)^r}\otimes \hatGam(\hatc_{r-2})^{(-1)^{r+1}}
    \otimes \cdots\otimes \hatGam(\hatc_0)^{-1}\\
    &= \ (-1)^{\calM(C^\b,\tilC^\b)+\calR(C^\b)+\calR(\tilC^\b)-\calR(C^\b\oplus\tilC^\b)}\cdot c_\hatGam.
  \end{aligned}
\end{equation}


Using the isomorphisms $\Gam:C^\b\to C^{d-\b}$ and $\tilGam:\tilC^\b\to \tilC^{d-\b}$ one sees that $\dim C^j= \dim C^{d-j}$ and $\dim \tilC^j= \dim
\tilC^{d-j}$. Therefore,
\meq{M=CtilC}\notag
    \calM(C^\b,\tilC^\b) \ = \ \sum_{0\le k<j\le d}\, \dim C^j\cdot \dim \tilC^k \\ = \ \sum_{0\le k<j\le d}\, \dim C^{d-j}\cdot \dim \tilC^{d-k}
    \ = \ \sum_{0\le k<j\le d}\, \dim C^k\cdot \dim \tilC^j.
\end{multline}
Hence,
\eq{M=CtilC2}
\begin{aligned}
    \calM(C^\b,\tilC^\b) \ &= \ \frac12\,\sum_{0\le k\not= j\le d}\, \dim C^j\cdot \dim \tilC^k
    \\& = \
    \frac12\,\Big[\, \big(\,\sum_{j=0}^d\,\dim C^j\,\big)\cdot\big(\,\sum_{j=0}^d\,\dim \tilC^j\,\big) \ - \ \sum_{j=0}^d\,\dim C^j\cdot\dim
    \tilC^j\,\Big]
    \\&= \ \frac12\,\Big[\, \big(\,2\,\sum_{j=0}^{r-1}\,\dim C^j\,\big)\cdot\big(\,2\,\sum_{j=0}^{r-1}\,\dim \tilC^j\,\big) \ - \
            2\,\sum_{j=0}^{r-1}\,\dim C^j\cdot\dim\tilC^j\,\Big]
    \\&\equiv \ \sum_{j=0}^{r-1}\,\dim C^j\cdot\dim\tilC^j \qquad\qquad\qquad \MOD 2.
\end{aligned}
\end{equation}

Using the identity \refe{x(x+1)}, we obtain from \refe{R(C)}
\eq{MRR=R}
  \calR(C^\b\oplus\tilC^\b) \ - \ \calR(C^\b) \ - \ \calR(\tilC^\b) \ = \ \sum_{j=0}^{r-1}\,\dim C^j\cdot \dim \tilC^j.
\end{equation}
Combining \refe{M=CtilC2} and \refe{MRR=R}, we get
\eq{M+R+R+R}
   \calM(C^\b,\tilC^\b) \ + \ \calR(C^\b) \ + \ \calR(\tilC^\b) \ - \ \calR(C^\b\oplus\tilC^\b) \ \equiv \ 0 \qquad\qquad \MOD 2.
\end{equation}
The identity \refe{cGamdirectsum} follows now from \refe{mucc}.
\eprf

\subsection{Dependence of the refined torsion on the chirality operator}\Label{SS:deponGam}
In this subsection $\bfk= \CC$ or $\RR$. Suppose that $\Gam_t$, $t\in \RR$, is a smooth family of chirality operators on the complex $(C^\b,\pa)$.
Let $\dot\Gam_t:C^\b\to C^{d-\b}$ denote the derivative of $\Gam_t$ with respect to $t$. Then, for each $k=0\nek d$, the composition
$\dot\Gam_t\circ\Gam_t$ maps $C^k$ into itself. In particular $\dot\Gam_t\circ\Gam_t:C^\even\to C^\even$ and $\dot\Gam_t\circ\Gam_t:C^\odd\to
C^\odd$. Define the {\em supertrace} $\Tr_s(\dot\Gam_t\circ\Gam_t)$ of $\dot\Gam_t\circ\Gam_t$ by the formula
\eq{supertrace}
    \Tr_s\,(\dot\Gam_t\circ\Gam_t) \ := \ \Tr\,(\dot\Gam_t\circ\Gam_t|_{C^\even})\ - \ \Tr\,(\dot\Gam_t\circ\Gam_t|_{C^\odd})
    \ = \ \sum_{j=0}^d\, (-1)^j\,\Tr\,(\dot\Gam_t\circ\Gam_t|_{C^j}).
\end{equation}

\prop{deponGam}
Let $(C^\b,\pa)$ be a length $d=2r-1$ complex of finite dimensional $\bfk$-vector spaces and let\/ $\Gam_t$, $t\in\RR$, be a smooth family of
chirality operators on $C^\b$. Then the following equality holds
\eq{deponGam}
    \frac{d}{dt}\,\rho_{{}_{\Gam_t}} \ = \ \frac12\,\Tr_s\,(\dot\Gam_t\circ\Gam_t).
\end{equation}
\eprop
\prf
Let $\Gam_{t,j}$ denote the restriction of $\Gam_t$ to $C^j$. Above we denoted the map $\Det(C^j)\to \Det(C^{d-j})$ induced by $\Gam_t$ by the same
symbol $\Gam_t$. To avoid confusion we will not use that convention in this proof and denote this map by $\Gam_{t,j}^{\Det}$.

For each $j=0\nek r-1$, $t_0\in \RR$ we have $\Gam_{t,j}= \Gam_{t,j}\circ \Gam_{t_0,d-j}\circ\Gam_{t_0,j}$ and, hence,
\eq{dGamGamj}
 \begin{aligned}
    \frac{d}{dt}{\big|_{t=t_0}}\,\Gam_{t,j}^{\Det} \ &= \ \frac{d}{dt}{\big|_{t=t_0}}\,\Big[\,
      \Det(\Gam_{t,j}\circ \Gam_{t_0,d-j})\cdot\Gam_{t_0,j}^{\Det}\,\Big]
    \\ &= \ \Tr(\dot\Gam_{t_0,j}\circ \Gam_{t_0,d-j})\cdot\Gam_{t_0,j}^{\Det},
 \end{aligned}
\end{equation}
where for the latter equality we used that for any smooth family of operators $A_t:C^{d-j}\to C^{d-j}$ one has $\frac{d}{dt}\Det(A_t)=
\Tr(\dot{A}_tA_t^{-1})$ and that $\Gam_{t,j}^{-1}= \Gam_{t,d-j}$. Hence, for any non-zero element $c_j\in \Det(C^j)$, we have
\eq{ddtGamcj}
    \frac{d}{dt}\big(\,\Gam_{t,j}^{\Det}(c_j)\,\big)^{\pm1} \ = \ \pm\,\Tr(\dot{\Gam}_{t,j}\circ\Gam_{t,d-j})\cdot
            \big(\,\Gam_{t,j}^{\Det}(c_j)\,\big)^{\pm1}.
\end{equation}
Using \refe{ddtGamcj} and the equality $(-1)^{j+1}= (-1)^{d-j}$ we conclude from the definition \refe{Gamc} of the element $c_{{}_{\Gam_t}}$ that
\eq{ddtcGam}
    \frac{d}{dt}\,c_{{}_{\Gam_t}} \ = \ \sum_{j=0}^{r-1}\,(-1)^{d-j}\,\Tr(\dot\Gam_{t,j}\circ \Gam_{t,d-j})\cdot c_{{}_{\Gam_t}}.
\end{equation}

Since $\Gam_{t,j}\circ \Gam_{t,d-j}=1$ we obtain
\[
    0 \ = \ \frac{d}{dt}\, \Tr(\Gam_{t,j}\circ \Gam_{t,d-j}) \ = \ \Tr(\dot\Gam_{t,j}\circ \Gam_{t,d-j})\ + \ \Tr(\Gam_{t,j}\circ
    \dot\Gam_{t,d-j}).
\]
Hence,
\eq{dotGamGam=-GamdotGam}
    \Tr(\dot\Gam_{t,j}\circ \Gam_{t,d-j}) \ = \ - \Tr(\dot\Gam_{t,d-j}\circ \Gam_{t,j}).
\end{equation}
Combining \refe{ddtcGam} with \refe{dotGamGam=-GamdotGam}, we obtain \refe{deponGam}.
\eprf

\subsection{The refined torsion of the dual complex}\Label{SS:reftordual}
Suppose now that $\bfk$ is endowed with an involutive endomorphism $\tau$, cf. \refs{dualcomplex}. Let $\hatC^\b$ be the $\tau$-dual complex of $C$
and let $\alp_{C^\b}: \Det(C^\b)\to \Det(\hatC^\b)$ denote the $\tau$-isomorphism defined in \refe{uC}. Let $\hGam$ be the $\tau$-adjoint of $\Gam$,
cf. \refss{tauadj}. Then $\hatGam$ is a chirality operator for the complex $\hatC^\b$.

\lem{reftordual}
In the situation described above,
\eq{reftordual}
    \displaystyle \rho_{{}_\hGam} \ = \ \alp_{H^\b(\pa)}(\rho_{{}_\Gam}).
\end{equation}
\elem
\prf
Fix $c_j\in \Det(C^j)$, $j=0\nek r-1$, and set
\eq{rthatc}
    \hatc_j \ = \ \alp_{C^{d-j}}^{-1}\big(\,(\Gam c_j)^{-1}\,\big)\ \in \ \Det(\hatC^j), \qquad\qquad  j=0\nek r-1.
\end{equation}
Then, by \refe{tauadj-1},
\eq{rtGamhatc}
    \hGam\hatc_j \ = \ \alp_{C^j}^{-1}(c_j^{-1})\ \in \ \Det(\hatC^{d-j}), \qquad \qquad j=0\nek r-1.
\end{equation}

Using \refe{u-1(v)-1}, we obtain from \refe{rthatc} and \refe{rtGamhatc}, that, for $j=0\nek r-1$,
\eq{rtbetc}
    \bet_{C^j}(c_j) \ = \ (-1)^{\dim C^j}\cdot (\hGam\hatc_j)^{-1},\qquad
    \bet_{C^{d-j}}(\Gam c_j) \ = \ (-1)^{\dim C^j}\cdot \hatc_j^{-1}.
\end{equation}
Combining \refe{u(Vb)2}, \refe{Gamc}, \refe{rthatc}, \refe{rtGamhatc}, and \refe{rtbetc}, we get
\eq{alpCcGam}
    \alp_{C^\b}(c_{{}_\Gam}) \ = \ (-1)^{\calM(C^\b)+\sum_{p=0}^{r-1}\dim C^{2p}}\cdot c_{{}_{\Gam^*}}.
\end{equation}
By definition \refe{refinedtor}, $\rho_{{}_\Gam} =  \phi_{C^\b}(c_\Gam)$. Therefore, from \refl{hatC-C}, we obtain
\eq{alpHrhoGam}
    \alp_{H^\b(\pa)}(\rho_{{}_\Gam}) \ = \  \phi_{\hatC^\b}\circ \alp_{C^\b}(c_{{}_\Gam})
    \ = \ (-1)^{\calM(C^\b)+\sum_{p=0}^{r-1}\dim C^{2p}}\cdot \rho_{{}_{\Gam^*}}.
\end{equation}
Since, by assumption $\Gam(C^j)= C^{d-j}$, we have $\dim C^j= \dim C^{d-j}$. Hence,
\[
    \sum_{p=0}^{r-1}\dim C^{2p} \ = \ \frac12\,\sum_{j=0}^{d}\,\dim C^j,
\]
and, by \refe{calF},
\eq{F(C)+sum=}
 \begin{aligned}
    \calM(C^\b)\ +\ \sum_{p=0}^{r-1}\, \dim C^{2p} \ &= \ \sum_{0\le j<k\le d}\,\dim C^j\cdot \dim C^k \ + \
    \frac12\,\sum_{j=0}^{d}\,\dim C^j \\ &= \ \frac12\,\sum_{j,k=0}^d\,\dim C^j\cdot \dim C^k \ = \
    \frac12\,\big(\,\sum_{j=0}^d\,\dim C^j\,\big)^2.
 \end{aligned}
\end{equation}
Using again the equality $\dim C^j= \dim C^{d-j}$, we obtain
\[
    \sum_{j=0}^d\,\dim C^j\ = \ 2\,\sum_{j=0}^{r-1}\,\dim C^j.
\]
Hence, from \refe{F(C)+sum=} we get
\eq{F(C)+sum=0}
    \calM(C^\b)\ +\ \sum_{p=0}^{r-1}\, \dim C^{2p} \ = \ 2\,\big(\,\sum_{j=0}^{r-1}\,\dim C^j\,\big)^2
    \ \equiv \ 0,\qquad \MOD\ 2.
\end{equation}
Combining \refe{alpHrhoGam} with \refe{F(C)+sum=0}, we obtain \refe{reftordual}.
\eprf

\section{Calculation of the Refined Torsion of a Finite Dimensional Complex}\Label{S:calcreftor}

In this section we introduce a finite dimensional analogue of the Atiyah-Patodi-Singer odd signature operator and express the refined torsion in
terms of the determinant of this operator.

Throughout the section we work under the assumptions of \refss{chiral}.

\subsection{The signature operator}\Label{SS:signfd}

The {\em signature operator} $\B$ is defined by the formula
\eq{signaturefd}
    \B \ := \ \Gam\,\pa\ +\ \pa\, \Gam.
\end{equation}
This is a direct analogue of the signature operator of an odd-dimensional manifold, cf. \cite[p.~44]{APS1}, \cite[p.~405]{APS2},
\cite[p.~64--65]{Gilkey84}. See also \refs{grdet} below.

Define
\eq{Cpm}
    C^j_+ \ := \ \Ker \big(\,\pa\circ\Gam\,\big)\cap C^j \ = \ \Gam\big(\,\Ker\, \pa\cap C^{d-j}\,\big), \qquad C^j_- \ := \ \Ker\pa\cap C^j,
    \qquad\qquad j=0\nek d,
\end{equation}
and set $C^{-1}_+= C^{d+1}_-= 0$.  Let $\B_j$ and $\B_j^\pm$ denote the restriction of $\B$ to $C^j$ and $C^j_\pm$ respectively. Then, for each
$j=0\nek d$, one has
\begin{align}
    \IM \B_j^+ \ &\subseteq \ \IM\big(\,\Gam\circ \pa|_{C^j}\,\big)
        \ \subseteq \ \Gam\big(\,\Ker \pa|_{C^{j+1}}\,\big) \ \subseteq \ C^{d-j-1}_+;\Label{E:IMB+}\\
    \IM \B_j^- \ &\subseteq \ \IM\big(\,\pa\circ \Gam|_{C^j}\,\big) \ \subseteq \ \IM\big(\,\pa|_{C^{d-j}}\,\big) \ \subseteq \ C^{d-j+1}_-.\Label{E:IMB-}
\end{align}
Hence,
\eq{Bkpm}
    \B_j^+ \ = \ \Gam\circ\pa:\, C^j_+ \ \longrightarrow \ C^{d-j-1}_+, \qquad
    \B_j^- \ = \ \pa\circ\Gam:\, C^j_- \ \longrightarrow \ C^{d-j+1}_-.
\end{equation}
Denote $C^\even:= \bigoplus_{j\ \even}\,C^j$,\ $C^\even_\pm:= \bigoplus_{j\ \even}\,C^j_\pm$. Set
\eq{Bevenodd}
    \B_\even\ :=\ \bigoplus_{j\ \even}\, \B_j:\,C^\even\ \to C^\even, \qquad \B_\even^\pm\ :=\ \bigoplus_{j\ \even} \B_j^\pm:\, C^\even_\pm\ \to \ C^\even_\pm,
\end{equation}
and define $\B_\odd, \ \B_\odd^\pm$ similarly. Note that $\B_\even = \Gam\circ{}\B_\odd\circ{}\Gam$. Hence, the whole information about $\B$ is
encoded in its {\em even part}\/ $\B_{\even}$.

\lem{Bsym-Cacyclic}
Suppose that the signature operator\/ $\B:C^\b\to C^\b$ is bijective. Then the complex $(C^\b,\pa)$ is acyclic and, for all $j=0\nek d$,
\eq{decompositionfd}
    C^j \ = \ C^j_+ \,\oplus \, C^j_-.
\end{equation}
\elem
\prf
If $c\in C^j_+\cap C^j_-$, then it follows from \refe{Cpm} and the definition \refe{signaturefd} of $\B$ that $\B{}c=0$. Hence, since $\B$ is
injective,  we obtain
\eq{C+capC-}
     C_+^j\cap C_-^j \ = \ \{0\}.
\end{equation}
Similarly, from \refe{Cpm} and \refe{signaturefd} we obtain $\IM\B\subset C^\b_++C^\b_-$. Hence, since $\B$ is surjective,
\eq{C=C+=C-}
    C^\b\ =\  C^\b_+\, +\, C^\b_-.
\end{equation}
Combining \refe{C+capC-} and \refe{C=C+=C-} we obtain \refe{decompositionfd}.

Let us now show that the complex $(C^\b,\pa)$ is acyclic. By \refe{decompositionfd},
\eq{B=B++B-}
    \B\ =\ \B^+\ +\ \B^-.
\end{equation}
By \refe{IMB+} and \refe{IMB-}, $\IM\B^\pm\subset C^\b_\pm$. Thus, since $\B=\B^++\B^-$ is surjective, it follows from \refe{decompositionfd}, that
$\IM\B^\pm= C^\b_\pm$. The equality $\B^-= \pa\circ\Gam$ implies now that $\IM\pa\supset C^\b_-$. Hence, by \refe{Cpm}, $\IM\pa= C^\b_-$. Combining
the latter equality with the definition \refe{Cpm} of $C^\b_-$, we obtain $\IM\pa= \Ker\pa$, proving the acyclicity of $(C^\b,\pa)$.
\eprf
\rem{Bdegenerate}
It is easy to construct an acyclic complex $(C^\b,\pa)$ and a chirality operator $\Gam$, so that the corresponding signature operator $\B$ is not
bijective.
\erem

\subsection{Calculation of the refined torsion in case $\B$ is bijective}\Label{SS:acylic}
Assume that the signature operator $\B:C^\b\to C^\b$ is bijective.  Then, by \refl{Bsym-Cacyclic}, the complex $(C^\b,\pa)$ is acyclic. Hence,
$\Det(H^\b(\pa))$ is canonically isomorphic to $\bfk$ and the refined torsion $\rho_{{}_\Gam}$ can be viewed as a number in $\bfk$. In this
subsection we calculate this number.

\defe{detgrfd}
The {\em graded determinant} of the even part of the signature operator is defined by the formula
\eq{detgrfd}
    \Detgrnp(\B_\even) \ = \ \Det(\B_\even^+)/\Det(-\B_\even^-).
\end{equation}
\edefe
Since \/ $\Gam\circ \B_\even^-\circ \Gam = \B_\odd^+$ \/  and $\Gam^2=\Id$, we have \/ $\Det(-\B_\even^-)= \Det(-\B_\odd^+)$ \/ and
\eq{detgrfd2}
 \begin{aligned}
    \Detgrnp(\B_\even) \ &= \ \Det(\B_\even^+)/\Det(-\B_\odd^+) \\ &= \ \Det\big(\,(-1)^{r-1}\,\Gam\pa|_{C^{r-1}_+}\,\big)^{(-1)^{r-1}}\cdot
     \prod_{j=1}^{r-1}\, \Det\big(\, (-1)^{j-1}\,\Gam\pa|_{C^{j-1}_+\oplus C^{d-j}_+}\,\big)^{(-1)^{j-1}}
    \\ &=\ (-1)^{(r-1)\dim C^{r-1}_+}\cdot \Det\big(\,\Gam\pa|_{C^{r-1}_+}\,\big)^{(-1)^{r-1}}\cdot
    \prod_{j=1}^{r-1}\, \Det\big(\, \Gam\pa|_{C^{j-1}_+\oplus C^{d-j}_+}\,\big)^{(-1)^{j-1}},
 \end{aligned}
\end{equation}
where in the last equality we used that
\[
    \Det\big(\, (-1)^{j-1}\,\Gam\pa|_{C^{j-1}_+\oplus C^{d-j}_+}\,\big) \ = \ \Det\big(\,\Gam\pa|_{C^{j-1}_+\oplus C^{d-j}_+}\,\big)
\]
since $\dim C^{j-1}_+= \dim C^{d-j}_+$.

\prop{reftor-grdetfd}
Suppose that the signature operator $\B$ is invertible and, hence, the complex $(C^\b,\pa)$ is acyclic. Then
\eq{reftor-grdetfd}
    \rho_{{}_\Gam} \ = \  \Detgrnp(\B_\even),
\end{equation}
\eprop
\prf
Recall that $\rho_{{}_\Gam}= \phi_{C^\b}(c_{{}_\Gam})$, where $c_{{}_\Gam}$ is the element in $\Det(C^\b)$ given by the formula \refe{Gamc} and
\[
    \phi_{C^\b}:\, \Det(C^\b) \ \longrightarrow \ \Det\big(H^\b(\pa)\big) \ \simeq \ \bfk
\]
is defined by \refe{isomorphism}.

To compute $\phi_{C^\b}(c_{{}_\Gam})$, we choose the decomposition \refe{decomposition} to be $C^j= C^j_-\oplus C^j_+$ and define elements $c_0\nek
c_d$ as follows: For each $j=0\nek d-1$, fix a non-zero element $a_j\in \Det(C^j_+)$ and set
\begin{gather}
    c_0 \ = \ a_0, \qquad c_d=\Gam a_{0},\notag\\
    c_j \ =\  \mu_{C^j_-,C^j_+}(\Gam a_{d-j}\otimes a_j), \qquad j=1\nek d-1.\notag
\end{gather}
Note that, for each $j=1\nek d$,
\meq{Gamc2}
    \Gam c_j \ = \ \mu_{C^{d-j}_+,C^{d-j}_-}( a_{d-j}\otimes \Gam a_j) \\ = \ (-1)^{\dim C^j_+\cdot \dim C^j_-}\,
    \mu_{C^{d-j}_-,C^{d-j}_+}(\Gam a_j\otimes  a_{d-j}) \ = \  (-1)^{\dim C^j_+\cdot \dim C^j_-}\, c_{d-j}.
\end{multline}
Thus, from \refe{Gamc}, we obtain
\eq{cj-cGam}
  \begin{aligned}
    c_{{}_\Gam} \ &= \
    (-1)^{\calR(C^\b)}\cdot
    c_0\otimes c_1^{-1}\otimes \cdots \otimes c_{r-1}^{(-1)^{r-1}}\otimes (\Gam c_{r-1})^{(-1)^r}\otimes (\Gam c_{r-2})^{(-1)^{r-1}}
    \otimes \cdots\otimes (\Gam c_0)^{-1}
    \\ &= \ (-1)^{\calR(C^\b)+\sum_{j=1}^{r-1}\dim C^j_+\cdot \dim C^j_-}\,
    c_0\wedge c_1^{-1}\wedge\cdots\wedge c_{d-1}\wedge c_d^{-1}.
  \end{aligned}
\end{equation}

To compute $\rho_{{}_\Gam}$ we now need to calculate $\phi_{C^\b}(c_0\wedge c_1^{-1}\wedge\cdots\wedge c_{d-1}\wedge c_d^{-1})$, i.e., in view of
\refe{isomorphism2}, we need to determine the elements $h_j\in \Det(H^j)\simeq \bfk$ which satisfy \refe{c=ahb}.

If $L$ is a $\bfk$-line and $x,y\in L$ with $y\not=0$, we denote by $[x:y]\in \bfk$ the unique number such that $x= [x:y]\,y$. Then, by \refe{c=ahb},
the elements $h_j\in \bfk$, which appear in \refe{isomorphism2}, are given by
\begin{align}
    h_0 \ &= \ 1, \Label{E:h0=1}\\
    h_d \ &= \ [c_d:\pa a_{d-1}] \ = \ [a_0:\Gam\pa a_{d-1}], \Label{E:hd=c:mu}
\end{align}
and, for $j=1\nek d-1$,
\eq{hj=c:mu}
 \begin{aligned}
    h_j\ &=  \ \big[\,c_j:\mu_{C^j_-,C^j_+}(\pa a_{j-1}\otimes a_j)\,\big] \\ &= \
    \big[\,\mu_{C^j_-,C^j_+}(\Gam a_{d-j}\otimes a_j):\mu_{C^j_-,C^j_+}(\pa a_{j-1}\otimes a_j)\,\big] \\ &= \
    \big[\,\Gam a_{d-j}:\pa a_{j-1}\,\big] \ = \ \big[\, a_{d-j}:\Gam\pa a_{j-1}\,\big].
 \end{aligned}
\end{equation}
By \refe{hd=c:mu} and \refe{hj=c:mu}, for $j=1\nek r-1$, we obtain
\eq{c:mu2}
 \begin{aligned}
    h_j\cdot h_{d-j+1}\ &= \ \big[\, a_{d-j}:\Gam\pa a_{j-1}\,\big]\cdot \big[\, a_{j-1}:\Gam\pa a_{d-j}\,\big]
    \\ &= \
    \big[\, \mu_{C^{d-j}_+,C^{j-1}_+}(a_{d-j}\otimes a_{j-1}):\mu_{C^{d-j}_+,C^{j-1}_+}(\Gam\pa a_{j-1}\otimes \Gam\pa a_{d-j}\,\big]
    \\ &= \ (-1)^{\dim C^{d-j}_+\cdot\dim C^{j-1}_+}\, \Det(\Gam\pa|_{C^{j-1}_+\oplus C^{d-j}_+})^{-1}.
 \end{aligned}
\end{equation}
By \refe{hj=c:mu} we have
\eq{c:mu3}
    h_r \ = \ \Det(\Gam\pa|_{C^{r-1}_+})^{-1}.
\end{equation}

Combining \refe{isomorphism2}, \refe{h0=1}, \refe{c:mu2}, and \refe{c:mu3}, we obtain
\eq{phi(c0...cd)}
    \phi_{C^\b}(c_0\wedge c_1^{-1}\wedge\cdots\wedge c_{d-1}\wedge c_d^{-1}) \ = \ (-1)^{\calN(C^\b)}\,\Det(\Gam\pa|_{C^{r-1}})^{(-1)^{r-1}}\cdot
    \prod_{j=1}^{r-1}\, \Det(\Gam\pa|_{C^{j-1}_+\oplus C^{d-j}_+})^{(-1)^{j-1}},
\end{equation}
From the definition \refe{refinedtor} of $\rho_{{}_\Gam}$ and the identities \refe{cj-cGam}, \refe{phi(c0...cd)}, we get
\eq{rho=}
    \rho_{{}_\Gam} \ = \ (-1)^{\calF(C^\b)}\, \Det(\Gam\pa|_{C^{r-1}})^{(-1)^{r-1}}\cdot
    \prod_{j=1}^{r-1}\, \Det(\Gam\pa|_{C^{j-1}_+\oplus C^{d-j}_+})^{(-1)^{j-1}},
\end{equation}
where
\eq{F(C)}
    \calF(C^\b) \ = \
    \calN(C^\b)\ + \ \calR(C^\b) \ + \ \sum_{j=0}^{r-1}\, \dim C^j_+\cdot \dim C^j_-\ + \ \sum_{j=0}^{r-1}\, \dim C^{d-j}_+\cdot \dim C^{j-1}_+.
\end{equation}

As the maps $\pa:C^{j-1}_+\to C^j_-$ and $\Gam:C^{d-j}_+\to C^j_-$ are isomorphisms, the last two terms in \refe{F(C)} can be computed to be
\meq{F(C)1}
    \sum_{j=0}^{r-1}\, \dim C^j_+\cdot \dim C^j_-\ + \ \sum_{j=0}^{r-1}\, \dim C^{d-j}_+\cdot \dim C^{j-1}_+
    \\ = \ \sum_{j=0}^{r-1}\, \dim C^j_+\cdot \dim C^{j-1}_+\ + \ \sum_{j=0}^{r-1}\, \dim C^{j}_-\cdot \dim C^{j-1}_+
     \ = \   \sum_{j=0}^{r-1}\, \dim C^j\cdot \dim C^{j-1}_+.
\end{multline}

Since the map $\Gam\pa:C^j_+\to C^{d-j-1}_+$ is an isomorphism, we have
\[
    \dim C^j_+\cdot\big(\,\dim C^j_++(-1)^{j+1}\,\big) \ = \ \dim C^{d-j-1}_+\cdot\big(\,\dim C^{d-j-1}_++(-1)^{d-j}\,\big),
    \qquad j=0\nek r-2.
\]
Hence, by formula \refe{N(C)} for $\calN(C)$ and the fact that for any $x\in \ZZ$, $x(x\pm 1)\equiv 0 \ (\MOD 2)$, we obtain
\eq{N(C)Gam}
 \begin{aligned}
    \calN(C^\b) \ &= \ \sum_{j=0}^{r-2}\,\dim C^j_+\cdot\big(\,\dim C^j_++(-1)^{j+1}\,\big) \ + \
    \frac12\ \dim C^{r-1}_+\cdot\big(\,\dim C^{r-1}_++(-1)^{r}\,\big)
    \\ &\equiv \ \frac12\ \dim C^{r-1}_+\cdot\big(\,\dim C^{r-1}_++(-1)^{r}\,\big)
    \\ &= \ \frac12\ \dim C^{r-1}_+\cdot\big(\,\dim C^{r-1}_+-1\,\big) \  + \ \frac12\,(1+(-1)^r)\cdot \dim C^{r-1}_+\qquad\qquad \MOD 2.
 \end{aligned}
\end{equation}

Next, using the isomorphism $\pa:C^{j-1}_+\to C^j_-$ we obtain from \refe{decompositionfd},
\eq{Cj=Cj++Cj-}
    \dim C^j  \ = \ \dim C^{j-1}_+ \ + \ \dim C^j_+.
\end{equation}
Hence, from definition \refe{R(C)} of $\calR(C^\b)$  and from identity \refe{x(x+1)}, we get
\meq{R(C)=CC+}
    \calR(C^\b) \ = \ \sum_{j=0}^{r-1}\,\Big[\, \frac12\,\dim C^j_+\cdot\big(\, \dim C^j_++(-1)^{r+j}\,\big) \\ + \
    \frac12\,\dim C^{j-1}_+\cdot\big(\, \dim C^{j-1}_++(-1)^{r+j-1}\,\big) \ + \ \dim C^j_+\cdot \dim C^{j-1}_+\,\Big]
    \\ = \
    \frac12\,\dim C^{r-1}_+\cdot\big(\, \dim C^{r-1}_+-1\,\big) \ + \ \sum_{j=0}^{r-2}\, (\dim C^j_+)^2
    \ + \ \sum_{j=0}^{r-1}\,\dim C^j_+\cdot \dim C^{j-1}_+
\end{multline}
By \refe{Cj=Cj++Cj-},
\meq{C2+C+C}
    \sum_{j=0}^{r-2}\, (\dim C^j_+)^2 \ + \ \sum_{j=0}^{r-1}\,\dim C^j_+\cdot \dim C^{j-1}_+ \\ = \
    \sum_{j=1}^{r-1}\,\,\Big[\, (\dim C^{j-1}_+)^2+\dim C^j_+\cdot\dim C^{j-1}_+\,\Big]
    \ \equiv \
    \sum_{j=0}^{r-1}\,\dim C^j\cdot \dim C^{j-1}_+.
\end{multline}
Hence, from \refe{R(C)=CC+}, we get
\eq{R(C)=CC+2}
    \calR(C^\b) \ \equiv \ \frac12\ \dim C^{r-1}_+\cdot\big(\, \dim C^{r-1}_+-1\,\big) \ + \
    \sum_{j=0}^{r-1}\,\dim C^j\cdot \dim C^{j-1}_+, \qquad\qquad \MOD 2
\end{equation}
Combining \refe{N(C)Gam} and \refe{R(C)=CC+2} and using again that $x(x\pm1)\equiv0$ for $x\in \ZZ$, we conclude that
\eq{N+R}
 \begin{aligned}
    \calN(C^\b)\ + \ \calR(C^\b) \  &\equiv \ \frac12\,(1+(-1)^r)\cdot \dim C^{r-1}_+ \ + \ \sum_{j=0}^{r-1}\, \dim C^j\cdot\dim C^{j-1}_+, \qquad\qquad\qquad \MOD 2.
 \end{aligned}
\end{equation}

Since $\frac{1+(-1)^r}2\equiv r-1$ modulo 2, we conclude from \refe{F(C)}, \refe{F(C)1}, and \refe{N+R}, that
\[
    \calF(C^\b)\ \equiv\ (r-1)\dim C^{r-1}_+, \qquad\qquad \MOD 2.
\]
The equality \refe{reftor-grdetfd} follows now from \refe{detgrfd2} and \refe{rho=}.
\eprf

\subsection{Calculation of the refined torsion in case $\B$ is not bijective}\Label{SS:generalcase}
In this section we don't assume that $\B$ is bijective. In particular, the complex $(C^\b,\pa)$ is not necessarily acyclic. For simplicity, we
restrict to the case $\bfk= \CC$.

Consider the operator $\B^2$. Note that $\B^2= (\Gam\pa)^2+(\pa\Gam)^2$ and $\B^2(C^j)\subset C^j$ for all $j=0\nek d$. For an arbitrary interval
$\calI\subset [0,\infty)$ and $j=0\nek d$,  we denote by $C^j_{\calI}\subset C^j$ the span of the generalized eigenvectors of the restriction of
$\B^2$ to $C^j$ corresponding to eigenvalues $\lam$ with $|\lam|\in \calI$.  Since both operators $\Gam$ and $\pa$ commute with $\calB$ and, hence,
with  $\B^2$ we have
\[
    \Gam:\, C^j_\calI\ \longrightarrow \ C^{d-j}_\calI, \quad \pa:\, C^j_\calI\ \longrightarrow \ C^{j+1}_\calI.
\]
Hence,  we obtain a subcomplex $C^\b_\calI$ of $C^\b$ and the restriction $\Gam_\calI$ of $\Gam$ to $C^\b_\calI$ is a chirality operator on this
complex. Let
$\pa_\calI$ and $\B_\calI$  denote the restriction of $\pa$ and  $\B$ to $C^\b_\calI$. Then $\B_\calI= \Gam_\calI\pa_\calI+ \pa_\calI\Gam_\calI$. 
\lem{BI>0}
If\/ $0\not\in \calI$ then the complex $(C^\b_\calI,\pa_\calI)$ is acyclic.
\elem
\prf
If $x\in \Ker\pa_\calI$ then $\B^2x= (\pa\Gam)^2x\in \IM\pa_\calI$. Hence,
\[
    \B^2_\calI:\, \Ker\pa_\calI \ \longrightarrow \ \IM\pa_\calI \ \subset \Ker\pa_\calI.
\]
Since the operator $\B^2_\calI:C^\b_\calI\to C^\b_\calI$ is invertible, we conclude that $\Ker\pa_\calI= \IM\pa_\calI$.
\eprf

For each $\lam\ge 0$, the complex $C^\b$ is a direct sum of the complex $C^\b_{[0,\lam]}$ and the acyclic complex $C^\b_{(\lam,\infty)}$. In
particular, $H^\b_{(\lam,\infty)}(\pa)=0$ and $H^\b_{[0,\lam]}(\pa)\simeq H^\b(\pa)$.  Hence, there are canonical isomorphisms
\[
        \Phi:\, \Det(H^\b_{(\lam,\infty)}(\pa))\ \longrightarrow\ \CC, \qquad \Psi:\,\Det(H^\b_{[0,\lam]}(\pa))\ \longrightarrow \
        \Det(H^\b(\pa)).
\]

\lem{tcdoth}
For every $t\in \Det(H^\b_{(\lam,\infty)}(\pa)),\ h\in \Det(H^\b_{[0,\lam]}(\pa))$
\eq{actiononH}
    \Phi(t)\cdot \Psi(h) \ = \ \mu_{H^\b_{(\lam,\infty)}(\pa),H^\b_{[0,\lam]}(\pa)}(t\otimes h).
\end{equation}
\elem
\prf
Since $H^\b_{(\lam,\infty)}(\pa)=0$, it follows from \refe{M(V,W)}, that $\calM(H^\b_{(\lam,\infty)}(\pa),H^\b_{[0,\lam]})= 0 $. The lemma follows
now from the definition \refe{fusiongr2} of the fusion isomorphism.
\eprf

In the sequel we will not distinguish between $t\in \Det(H^\b_{(\lam,\infty)}(\pa))$ and $\Phi(t)\in \CC$ and write simply $t$ for $\Phi(t)$.
Similarly, for $h\in \Det(H^\b_{[0,\lam]}(\pa))$ we will denote by $h$ also the element $\Psi(h)\in \Det(H^\b(\pa))$.

From \refl{rhodirectsum}, \refp{reftor-grdetfd}, and \refl{tcdoth}, we immediately obtain the following 
\prop{antor-grdetfd2}
Let $(C^\b,\pa)$ be a complex of finite dimensional complex vector spaces and let $\Gam$ be a chirality operator on $C^\b$. Then, for each
$\lam\ge0$,
\eq{antor-grdetfd2}
    \rho_{{}_\Gam} \ = \ \Detgrnp(\B_{(\lam,\infty)})\cdot\rho_{{}_{\Gam_{\hskip-1pt{}_{[0,\lam]}}}},
\end{equation}
where, as above, we view $\rho_{{}_{\Gam_{\hskip-1pt{}_{[0,\lam]}}}}$ as an element of\/ $\Det(H^\b(\pa))$ via the canonical isomorphism\linebreak
$\Psi:\Det(H^\b_{[0,\lam]}(\pa))\to \Det(H^\b(\pa))$.
\eprop

\section{Preliminaries on Determinants and the $\eta$-invariant of Elliptic Operators}\Label{S:determinants}

In this section we briefly review the main facts about the $\zet$-regularized determinants and $\eta$-invariants of non self-adjoint elliptic
operators. In particular, we define a sign-refined version of the graded determinant -- a notion, which plays a central role in this paper. We refer
the reader to Sections~3 and 4 of \cite{BrKappelerRAT} for a more detailed discussion of the subject.

Let $E$ be a complex vector bundle over a smooth compact manifold $M$ and let $D:C^\infty(M,E)\to C^\infty(M,E)$ be an elliptic differential operator
of order $m\ge 1$. Denote by $\symb(D)$ the leading symbol of $D$.

\subsection{Choice of an angle}\Label{SS:agmonangle}
Our aim is to define the $\zeta$-function and the determinant of $D$. For this we will need to define the complex powers of $D$. As usual, to define
complex powers we need to choose a {\em spectral cut} in the complex plane. We restrict ourselves to the spectral cuts given by a ray
\eq{Rtet}
    R_\tet \ = \ \big\{\, \rho e^{i\tet}:\, 0\le \rho<\infty\, \big\},
    \qquad 0\le \tet < 2\pi.
\end{equation}
Consequently, we have to choose an angle $\tet\in [0,2\pi)$.

\defe{prangle}
The angle $\tet$ is a {\em principal angle} for an elliptic operator $D$ if
\eq{prangle}\notag
    \spec\big(\, \symb(D)(x,\xi)\, \big)\ \cap\ R_\tet \ = \ \emptyset,
    \qquad \text{for all}\quad x\in M,\ \xi\in T^*_xM\backslash\{0\}.
\end{equation}
\edefe

If $\calI\subset \RR$ we denote by $L_{\calI}$ the solid angle
\eq{Lab}\notag
    L_{\calI} \ = \  \big\{\, \rho e^{i\tet}:\, 0 < \rho<\infty,\, \tet\in \calI\,  \big\}.
\end{equation}
The existence of a principal angle is an additional assumption on $D$. Since $M$ is compact every operator which possesses a principal angle has a
discrete spectrum.

\defe{Agmon}
The angle $\tet$ is an {\em Agmon angle}
\footnote{Note that in the literature the notion of Agmon angle is often defined differently, namely it is required, in addition,  that zero is nor in the spectrum
of the operator.} for an elliptic operator $D$ if it is a principal angle for $D$ and there exists $\eps>0$ such that
\[
    \spec(D) \, \cap L_{[\tet-\eps,\tet+\eps]} \ = \ \emptyset.
\]
\edefe

If  $\tet$ is a principal angle for $D$, then, cf. \cite{Seeley67, ShubinPDObook}, there exists $\eps>0$ such that $\spec(D)\cap
L_{[\tet-\eps,\tet+\eps]}$ is finite and $\spec(\symb(D))\cap L_{[\tet-\eps,\tet+\eps]}=\emptyset$. Hence there exists an Agmon angle $\tet'\in
(\tet-\eps,\tet+\eps)$ for $D$.

\subsection{$\zet$-function and determinant}\Label{SS:zet-det}
Assume that $\tet$ is an Agmon angle for $D$. Let $\Pi:L^2(M,E)\to L^2(M,E)$ denote the spectral projection of $D$ corresponding to all non-zero
eigenvalues of $D$. The $\zeta$-function $\zeta_\tet(s,D)$ of $D$ is defined as follows.

Since, by assumption, $D$ possesses a principal angle, its spectrum is discrete. Hence, there exists a small number $\rho_0>0$ such that
\eq{rho0}\notag
    \spec(D) \, \cap \, \big\{\, z\in \CC;\, |z|<2\rho_0\, \big\} \ \subseteq \ \{0\}.
\end{equation}
Define the contour $\Gam= \Gam_{\tet,\rho_0}\subset \CC$ consisting of three curves $\Gam= \Gam_1\cup \Gam_2\cup \Gam_3$, where
\begin{gather}\Label{E:Gamtetrho}\notag
    \Gam_1 \ = \ \big\{\, \rho e^{i\tet}:\, \infty >\rho\ge \rho_0\, \big\},
    \quad
    \Gam_2 \ = \ \big\{\, \rho_0 e^{i\alp}:\, \tet< \alp<\tet+2\pi\, \big\},\\
    \quad
    \Gam_3 \ = \ \big\{\, \rho e^{i(\tet+2\pi)}:\, \rho_0\le \rho<\infty\, \big\}.
\end{gather}

For $\RE s> \frac{\dim M}m$, the operator
\eq{Ds}
    \Pi\,D_\tet^{-s} \ = \ \frac{i}{2\pi}\, \int_{\Gam_{\tet,\rho_0}}\, \lam^{-s}(D-\lam)^{-1}\, d\lam
\end{equation}
is a pseudo-differential operator with continuous kernel $K_\tet(s;x,y)$, cf. \cite{Seeley67, ShubinPDObook}. In particular, the operator
$\Pi\,D_\tet^{-s}$ is of trace class. 

We define
\eq{zeta}
    \zeta_{\tet}(s,D) \ = \ \Tr \Pi\,D_\tet^{-s} \ = \ \int_M\, \tr K_\tet(s;x,x)\, dx,
    \qquad \RE s> \frac{\dim M}m.
\end{equation}
It was shown by Seeley \cite{Seeley67} (see also \cite{ShubinPDObook}) that $\zeta_{\tet}(s,D)$ has a meromorphic extension to the whole complex
plane and that 0 is a regular value of $\zeta_{\tet}(s,D)$.

More generally, let $Q$ be a pseudo-differential operator of order $q$. We set
\eq{zetaQ}
    \zeta_{\tet}(s,Q,D) \ = \  \Tr\, Q\,\Pi\,D_\tet^{-s},
    \qquad \RE s> (q+\dim M)/m.
\end{equation}
This function also has a meromorphic extension to the whole complex plane, see \cite[\S3.22]{Wodzicki87} and \cite[Th.~2.7]{Grubbseeley95} (see also
\cite{Guillemin85}). Moreover, if $Q$ is a $0$-th order pseudo-differential projection, i.e. a 0-th order pseudo-differential operator satisfying
$Q^2=Q$, then by \cite[\S7]{Wodzicki84}, \cite{Wodzicki87} (see also \cite{BruningLesch99,Ponge-asymetry} for a shorter proof), $\zeta_{\tet}(s,Q,D)$
is regular at 0.

If the dimension of $M$ is odd and $D$ is a {\em bijective differential}\/ operator of even order, then $\zet_\tet(0,D)=0$, cf. \cite{Seeley67}. More
generally, we
have the following 
\prop{zet=0}
Suppose $\dim M$ is odd, $D:C^\infty(M,E)\to C^\infty(M,E)$ is an elliptic differential operator of even order $m\ge2$, $\tet$ is an Agmon angle for
$D$, and $P$ is a finite rank pseudo-differential projection which commutes with $D$. Set $Q=\Id-P$ and assume that the restriction $D|_{\IM Q}$ of
$D$ to the image of $Q$ defines an invertible operator $D|_{\IM Q}:\IM{}Q\to \IM{Q}$. Then,
\eq{zet=0}
    \zet_\tet(0,Q,D) \ = \ -\,\rank\, (\Id-Q).
\end{equation}
In particular, if $m_0$ denotes the dimension of the span of the generalized eigenvectors of $D$ corresponding to the eigenvalue  $\lam=0$ (i.e.
$m_0$ is the algebraic multiplicity of the eigenvalue $\lam=0$ of $D$), then
\eq{zet=02}
    \zet_\tet(0,D) \ = \ -\,m_0.
\end{equation}
\eprop
\prf
If $\eps\not=0$ is a small enough real number, then $D+\eps$ is an invertible differential operator of even order and $\tet$ is an Agmon angle for
$D+\eps$. Hence, $\zet_\tet(0,D+\eps) = 0$, cf. \cite{Seeley67}. Clearly,
\eq{zet=zeteps+P}
    \zet_\tet(0,Q,D+\eps) \ = \ \zet_\tet(0,D+\eps) \ - \ \rank\,(\Id- Q) \ = \ -\,\rank\,(\Id- Q).
\end{equation}
Since 0 is not an eigenvalue of the restriction of $D$ to the image of $Q$, we have
\eq{limzeteps}
    \lim_{\eps\to 0}\,\zet_\tet(0,Q,D+\eps) \ = \ \zet_\tet(0,Q,D).
\end{equation}
Combining \refe{zet=zeteps+P} and \refe{limzeteps}, we obtain \refe{zet=0}.
\eprf

\defe{determinant}
The  $\zeta$-regularized determinant of $D$ is defined by the formula
\eq{logdet}
    \Det'_\tet (D) \ := \ \exp\,\left(\,-\frac{d}{ds}\big|_{s=0}\zeta_{\tet}(s,D)\,\right).
\end{equation}
\edefe
Roughly speaking, \refe{logdet} says that the logarithm $\log\Det'_\tet(D)$ of the determinant of $D$ is equal to $-\zeta_{\tet}'(0,D)$. However, the
logarithm is a multivalued function. Hence,  $\log\Det'_\tet(D)$ is defined only up to a multiple of $2\pi{i}$, while $-\zeta_{\tet}'(0,D)$ is a well
defined complex number. We denote by $\LD_\tet(D)$ the particular value of the logarithm of the determinant such that
\eq{logdet2}
    \LD_\tet(D) \ = \ -\,\frac{d}{ds}\big|_{s=0}\,\zeta_{\tet}(s,D).
\end{equation}
Let us emphasize that the equality \refe{logdet2} is the definition of the number $\LD_\tet(D)$.

\rem{Det'-Det}
The prime in $\Det'_\tet(D)$ and $\LD_\tet(D)$ indicates that we ignore the zero eigenvalues of $D$ in the definition of the regularized determinant.
If the operator $D$ is invertible we usually omit the prime and write $\Det_\tet(D)$ and $\operatorname{LDet}_\tet(D)$ instead.
\erem

We will need the following generalization of \refd{determinant}.
\defe{logdetV}
Suppose $Q$ is a $0$-th order pseudo-differential projection commuting with $D$. Then $V:= \IM{}Q$ is a $D$ invariant subspace of $C^\infty(M,E)$.
The $\zet$-regularized determinant of the restriction $D{|_V}$ of $D$ to $V$ is defined by the formula
\eq{logdetV}
    \Det'_\tet (D{|_V}) \ := \ e^{\LD_\tet(D|_V)},
\end{equation}
where
\eq{logdetV2}
    \LD_\tet(D|_V) \ = \ -\frac{d}{ds}\big|_{s=0}\zeta_{\tet}(s,Q,D).
\end{equation}
\edefe
As in \refr{Det'-Det}, if the restriction of $D$ to $V$ defines an invertible operator $D|_V:V\to V$, we usually omit the prime in the notation for
the numbers \refe{logdetV} and \refe{logdetV2} and write $\Det_\tet(D{|_V})$ and $\LDnp_\tet(D|_V)$ instead.
\rem{DetDV}
The  right hand side of \refe{logdetV2} is independent of $Q$ except through $\IM(Q)$. This justifies the notation $\LD_\tet(D|_V)$. However, we need
to know that $V$ is the image of a 0-th order pseudo-differential projection $Q$ to ensure that $\zeta_{\tet}(s,D)$ has a meromorphic extension to
the whole $s$-plane with  $s=0$ being a regular point.
\erem

\subsection{Spectral subspaces}\Label{SS:spectralsubspace}
We will often use \refd{logdetV} in the following special situation. For $\lam\ge0$ let $\Pi_{D,[0,\lam]}$ denote the spectral projection of $D$
corresponding to the set $\{z\in \CC:\,|z|\le\lam\}$. It is given by the Cauchy integral
\eq{PiD0lam}\notag
    \Pi_{D,[0,\lam]} \ = \ \frac{i}{2\pi}\,\int_0^{2\pi}\,\big(\,D-(\lam+\eps)\,e^{i\phi}\,\big)^{-1}\,d\phi,
\end{equation}
where $\eps>0$ is small enough so that there are no eigenvalues of $D$ with absolute value in the interval $(\lam,\lam+\eps]$. Since the operator $D$
is an elliptic differential operator of order $>0$, the image of $\Pi_{D,[0,\lam]}$ is finite dimensional and consists of smooth sections. We denote
by $C^\infty_{[0,\lam]}(M,E)\subset C^\infty(M,E)$ the image of $\Pi_{D,[0,\lam]}$. Note that $C^\infty_{[0,\lam]}(M,E)$ is equal to the span of the
generalized eigenvectors of $D$ corresponding to eigenvalues with absolute value $\le\lam$.

Define the projections
\eq{Pimu}
 \begin{aligned}
    \Pi_{D,(\lam,\infty)}\ &= \ \Id\ -\ \Pi_{D,[0,\lam]},\\
    \Pi_{D,(\lam,\mu]}\ &=\ \Pi_{D,[0,\mu]}\ -\ \Pi_{D,[0,\lam]}, \qquad\text{for}\quad \mu\ge\lam.
 \end{aligned}
\end{equation}
The range of  $\Pi_{D,(\lam,\mu]}$ is finite dimensional and contained in $C^\infty(M,E)$. It is equal to the span of the generalized eigenvectors of
$D$ with eigenvalues $\chi$ such that $\lam<|\chi|\le \mu$. The range of $\Pi_{D,(\lam,\infty)}$ is infinite dimensional and contains the span of the
generalized eigenvectors of $D$ with eigenvalues whose absolute value is greater than $\lam$, cf. \cite[Appendix~B]{Ponge-asymetry}.

Let now $\calI$ be an interval of the form $[0,\lam],\ (\lam,\mu]$, or $(\lam,\infty)$. Then $\Pi_{D,\calI}$ maps smooth sections to smooth sections
and the space
\eq{CcalI}\notag
    C^\infty_\calI(M,E) \ := \ \Pi_{D,\calI}\big(\,C^\infty(M,E)\,\big) \ \subset \ C^\infty(M,E)
\end{equation}
is $D$ invariant. Let $D_\calI$ denote the restriction of $D$ to the space $C^\infty_\calI(M,E)$. Note also that $D_\calI$ is invertible whenever
$0\not\in \calI$.

\refd{logdetV} gives us the determinant $\Det_\tet'(D_\calI)$. Clearly, for any $0\le\lam\le \mu$,
\eq{Det=DetlamDetmu}
    \Det_\tet\big(\,D_{(\lam,\infty)}\,\big) \ = \ \Det_\tet\big(\,D_{(\lam,\mu]}\,\big)\cdot \Det_\tet\big(\,D_{(\mu,\infty)}\,\big).
\end{equation}

\subsection{Dependence of the determinant on the angle}\Label{SS:det-tet}
Assume that  $\tet$ is a principal angle for $D$. Then, cf. \cite{Seeley67, ShubinPDObook}, for any $\eps>0$, we can choose an Agmon angle $\tet'\in
(\tet-\eps,\tet+\eps)$ for $D$. Let $\tet''> \tet'$ be another Agmon angle for $D$ such that all the angles in the interval $[\tet',\tet'']$ are
principal for $D$. Then, cf., for example, \cite[\S3.10]{BrKappelerRAT},
\eq{zet-zet}
    \frac{d}{ds}\big|_{s=0}\,\zet_{\tet'}(s,D) \ \equiv \ \frac{d}{ds}\big|_{s=0}\,\zet_{\tet''}(s,D) \qquad \MOD\ 2\pi i.
\end{equation}
Hence, by \refd{determinant},
\eq{det-det}
        \Det'_{\tet''} (D) \ = \ \Det'_{\tet'}(D).
\end{equation}

Note that the equality \refe{det-det} holds because both angles, $\tet'$ and $\tet''$, are close to a given principal angle $\tet$ so that the
intersection \/ $\spec(D)\cap L_{[\tet',\tet'']}$\, is finite. If there are infinitely many eigenvalues of $D$ in the solid angle
$L_{[\tet',\tet'']}$ then $\Det'_{\tet'}(D)$ and $\Det'_{\tet''}(D)$ might be different.

\subsection{Graded determinant}\Label{SS:grdeterminant}
Let $D:C^\infty(M,E)\to C^\infty(M,E)$ be a differential operator. Suppose that $Q_j:C^\infty(M,E)\to C^\infty(M,E)$ ($j=0\nek d$) are $0$-th order
pseudo-differential projections commuting with $D$. Set $V_j:= \IM{}Q_j$ and assume that
\[
    C^\infty(M,E)= \bigoplus_{j=0}^dV_j.
\]
\defe{grdeterminantV}
Assume that $\tet\in [0,2\pi)$ is an Agmon angle for the operator $(-1)^jD|_{V_j}$, for every $j=0\nek d$. The {\em graded determinant}
$\Detgrtet(D)$ of $D$ (with respect to the grading defined by the pseudo-differential projections $Q_j$) is defined by the formula
\eq{grdeterminantV}
    \Detgrtet(D) \ := \ e^{\LDgrtet(D)}.
\end{equation}
where
\eq{grldeterminantV}
    \LDgrtet(D) \ := \ \sum_{j=0}^d\, (-1)^j\,\LD_\tet\,\big(\, (-1)^jD|_{V_j}\,\big).
\end{equation}
\edefe

The following is an important example of the above situation: Let $E=\bigoplus_{j=0}^d E_j$ be a {\em graded} vector bundle over $M$. Suppose that
for each $j=0\nek d$, there is a bijective elliptic differential operator
\[
    D_j:\,C^\infty(M,E_j) \ \longrightarrow \  C^\infty(M,E_j),
\]
such that $\tet\in [0,2\pi)$ is an Agmon angle for $(-1)^jD_j$ for all $j=0\nek d$. We denote by
\eq{D=bigoplus}
    D \ = \ \bigoplus_{j=0}^d\,D_j:\, C^\infty(M,E)\ \longrightarrow \  C^\infty(M,E)
\end{equation}
the direct sum of the operators $D_j$. Then \refe{grldeterminantV} reduces to
\eq{grldeterminant}
    \LDgrtet(D) \ = \ \sum_{j=0}^d\, (-1)^j\,\LD_\tet\,\big(\, (-1)^jD_j\,\big).
\end{equation}

\subsection{Case of a self-adjoint leading symbol}\Label{SS:det-sa}
Let $h^E$ be a Hermitian metric on the bundle $E\to M$ and assume that the principal symbol $\symb(D)(x,\xi)$ of the elliptic operator $D$ is
self-adjoint, i.e.,
\eq{close2sa}
    \symb(D)^*(x,\xi)  \ = \ \symb(D)(x,\xi), \qquad (x,\xi)\in T^*M,
\end{equation}
where $\symb(D)^*(x,\xi)$ denotes the adjoint of the operator $\symb(D)(x,\xi)$ with respect to the scalar product $h^{E}$. This assumption implies
that $D$ can be written as a sum $D=D'+A$ where $D'$ is a self-adjoint differential operator of order $m$ and $A$ is a differential operator of order
smaller than $m$. If the leading symbol of $D$ is self-adjoint then any angle $\tet\not= 0,\pi$ is principal for $D$.

Though the operator $D$ is not self-adjoint in general, the assumption \refe{close2sa} guarantees that it has nice spectral properties, cf.
\cite[\S{}I.6]{Markus88} and \cite[\S3.9]{BrKappelerRAT}. Though many of the results of this paper remain valid for arbitrary elliptic differential
operators which possess an Agmon angle, for simplicity of notation we will often assume that our operators have a self-adjoint leading symbol.

\subsection{$\eta$-invariant}\Label{SS:etainv}
It is well known, cf. \cite{Singer85Dirac,Wojciechowski99}, that the phase of the determinant of a self-adjoint elliptic differential operator $D$
can be expressed in terms of the $\eta$-invariant of $D$ and the $\zet$-function of $D^2$. We now extend this result to non self-adjoint operators.

First, we recall the definition of the $\eta$-function of $D$ for a non-self-adjoint operator, cf. \cite{Gilkey84}. 
\defe{eta}
Let $D:C^\infty(M,E)\to C^\infty(M,E)$ be an elliptic differential operator of order $m\ge 1$ with self-adjoint leading symbol. Assume that $\tet$ is
an Agmon angle for $D$ (cf. \refd{Agmon}). Let $\Pi_>$  (resp. $\Pi_<$) be a pseudo-differential projection whose image contains the span of all
generalized eigenvectors of $D$ corresponding to eigenvalues $\lam$ with $\RE\lam>0$ (resp. with $\RE\lam<0$) and whose kernel contains the span of
all generalized eigenvectors of $D$ corresponding to eigenvalues $\lam$ with $\RE\lam\le0$ (resp. with $\RE\lam\ge0$), cf.
\cite[Appendix~B]{Ponge-asymetry}. We define the $\eta$-function of $D$ by the formula
\eq{eta}
    \eta_{\tet}(s,D) \ = \ \zet_\tet(s,\Pi_>,D) \ - \ \zet_\tet(s,\Pi_<,-D).
\end{equation}
\edefe
Note that, by definition, the purely imaginary eigenvalues of $D$ do not contribute to $\eta_\tet(s,D)$.

It was shown by Gilkey, \cite{Gilkey84}, that $\eta_\tet(s,D)$ has a meromorphic extension to the whole complex plane $\CC$ with isolated simple
poles, and that it is regular at $0$. Moreover, the number $\eta_\tet(0,D)$ is independent of the Agmon angle $\tet$.

Since the leading symbol of $D$ is self-adjoint, the angles $\pm\pi/2$ are principal angles for $D$, cf. \refd{prangle}. In particular, there are at
most finitely many eigenvalues of $D$ on the imaginary axis.

Let $m_+(D)$ (resp., $m_-(D)$) denote the number of eigenvalues of $D$, counted with their algebraic multiplicities, on the positive (resp.,
negative) part of the imaginary axis. Let $m_0(D)$ denote algebraic multiplicity of 0 as an eigenvalue of $D$.

\defe{etainv}
The $\eta$-invariant $\eta(D)$ of $D$ is defined by the formula
\eq{etainv}
    \eta(D) \ = \ \frac{\eta_\tet(0,D)+m_+(D)-m_-(D)+m_0(D)}2.
\end{equation}
\edefe
As $\eta_\tet(0,D)$ is independent of the choice of the Agmon angle $\tet$ for $D$, cf. \cite{Gilkey84}, so is $\eta(D)$.

Let $D(t)$ be a smooth 1-parameter family of elliptic operators with self-adjoint leading symbol. Then $\eta(D(t))$ is, in general, not smooth but
may have integer jumps when eigenvalues cross the imaginary axis or cross 0 along the imaginary axis. Because of this, the $\eta$-invariant is
usually considered modulo integers. However, in this paper we will be interested in the number $e^{i\pi\eta(D)}$, which changes its sign when
$\eta(D)$ is changed by an odd integer. Hence, we will consider the $\eta$-invariant as a complex number.

\rem{etainv}
Note that our definition of $\eta(D)$ is slightly different from the one proposed by Gilkey in \cite{Gilkey84}. In fact, in our notation, Gilkey's
$\eta$-invariant is given by $\eta(D)+m_-(D)$. Hence, reduced modulo integers, the two definitions coincide. However, the number $e^{i\pi\eta(D)}$
will be multiplied by $(-1)^{m_-(D)}$ if we replace one definition by the other. In this sense, \refd{etainv} can be viewed as a {\em sign
refinement} of the definition given in \cite{Gilkey84}.
\erem

\subsection{Relationship between the $\eta$-invariant and the determinant}\Label{SS:det-eta}
Since the leading symbol of $D$ is self-adjoint, the angles $\pm\pi/2$ are principal for $D$. Hence, there exists an Agmon angle $\tet\in (-\pi/2,0)$
such that there are no eigenvalues of $D$ in the solid angles $L_{(-\pi/2,\tet]}$ and $L_{(\pi/2,\tet+\pi]}$. Then $2\tet$ is an Agmon angle for the
operator $D^2$.

\th{det-eta}
Let $D:C^\infty(M,E)\to C^\infty(M,E)$ be an elliptic differential operator of order $m\ge 1$ with self-adjoint leading symbol. Assume $\tet\in
(-\pi/2,0)$ is an Agmon angle for $D$ such that there are no eigenvalues of $D$ in the solid angles $L_{(-\pi/2,\tet]}$ and $L_{(\pi/2,\tet+\pi]}$
(Hence, there are no eigenvalues of $D^2$ in the solid angle $L_{(-\pi,2\tet]}$). Then \footnote{Recall that we denote by $\LD_{\tet}(D)$ the
particular branch of the logarithm of the determinant of $D$ defined by formula \refe{logdet2}.}
\eq{det-eta}
   \LD_{\tet}(D) \ = \ \frac12\,\LD_{2\tet}(D^2)  \ - \
   i\pi\,\Big(\,\eta(D)- \frac{\zet_{2\tet}(0,D^2)+m_0(D)}2\,\Big).
\end{equation}
In particular,
\eq{det-eta2}
 \Det_{\tet}'(D) \ = \ e^{-\frac12\,\zet_{2\tet}'(0,D^2)}\cdot e^{-i\pi\,(\,\eta(D)- \frac{\zet_{2\tet}(0,D^2)+m_0(D)}2\,)}.
\end{equation}
\eth
In the case when $D$ is invertible the theorem is proven in Section~4 of \cite{BrKappelerRAT}. The same arguments without any changes prove
\reft{det-eta} in the general case.

\rem{det-eta}
a. \ Let $\tet$ be as in \reft{det-eta} and suppose that $\tet'\in (-\pi,0)$ is another angle such that both $\tet'$ and $\tet'+\pi$ are Agmon angles
for $D$. Then, by \refe{det-det},
\eq{tet-tet'}
  \begin{aligned}
    \Det_{\tet'}'(D) \ &= \ \Det_{\tet}'(D), \\ \zet_{2\tet}'(0,D^2)\ &\equiv \ \zet_{2\tet'}'(0,D^2) \qquad \MOD\, 2\pi i.
  \end{aligned}
\end{equation}
In particular,
\eq{tet-tet'2}
    e^{-\frac12\,\zet_{2\tet'}'(0,D^2)} \ = \ \pm\, e^{-\frac12\,\zet_{2\tet}'(0,D^2)}.
\end{equation}

Clearly, $\zet_{\tet_1}(0,D^2)= \zet_{\tet_2}(0,D^2)$ if there are finitely many eigenvalues of $D^2$ in the solid angle $L_{[\tet_1,\tet_2]}$.
Hence, $\zet_{2\tet}(0,D^2)= \zet_{2\tet'}(0,D^2)$. We then conclude from \refe{det-eta2}, \refe{tet-tet'}, and \refe{tet-tet'2} that
\eq{set-etatet'}
 \Det_{\tet'}'(D) \ = \ \pm\,e^{-\frac12\,\zet_{2\tet'}'(0,D^2)}\cdot e^{-i\pi\,(\,\eta(D)- \frac{\zet_{2\tet'}(0,D^2)+m_0(D)}2\,)}.
\end{equation}
In other words, for \refe{det-eta2} to hold we need the precise assumption on $\tet$ which are specified in \reft{det-eta}. But ``up to a sign" it
holds for every spectral cut in the lower half plane.

b. \ If instead of the spectral cut $R_{\tet}$ in the lower half-plane we use the spectral cut $R_{\tet+\pi}$ in the upper half-plane we will get a
similar formula
\eq{det-eta1}
   \LD_{\tet+\pi}(D) \ = \ \frac12\,\LD_{2\tet}(D^2)  \ + \
   i\pi\,\Big(\,\eta(D)- \frac{\zet_{2\tet}(0,D^2)+m_0(D)}2\,\Big),
\end{equation}
whose proof is a verbatim repetition of the proof of \refe{det-eta}, cf. Section~4 of \cite{BrKappelerRAT}.

c. \ {\em If the dimension of $M$ is odd}, then, in view of \refp{zet=0},\/ $\zet_{2\tet}(0,D^2)= -m_0(D)$. Hence, \refe{det-eta} simplifies to
\eq{det-eta-odd}
   \LD_{\tet}(D) \ = \ \frac12\,\LD_{2\tet}(D^2)  \ - \ i\pi\,\eta(D).
\end{equation}
\erem

\subsection{$\eta$-invariant and graded determinant}\Label{SS:eta-grdet}
Suppose now that $D= \bigoplus_{j=0}^dD_j$ is as in \refe{D=bigoplus}. Choose  $\tet\in (-\pi/2,0)$ such that there are no eigenvalues of $D_j$ in
the solid angles $L_{(-\pi/2,\tet]}$ and $L_{(\pi/2,\tet+\pi]}$ for every $0\le j\le  d$. From \refd{etainv} of the $\eta$-invariant it follows that
\[
    \eta( -D_j) \ = \ -\, \eta(D_j) \ + \ m_0(D_j).
\]
Combining this latter equality with \refe{grldeterminant} and \refe{det-eta} we obtain
\eq{grdet-eta}
    \LDgrtet(D) \ = \ \frac12\,\sum_{j=0}^d\, (-1)^j\,\LD_{2\tet}(D^2_j)  \ - \
   i\pi\,\Big(\,\eta(D)-\frac{m_0(D)}2-\frac12\sum_{j=0}^d\,(-1)^j\,\zet_{2\tet}(0,D^2_j)\,\Big),
\end{equation}
where $\eta(D)= \sum_{j=0}^d\eta(D_j)$ is the $\eta$-invariant of the operator $D= \bigoplus_{j=0}^dD_j$ and $m_0(D)= \sum_{j=0}^dm_0(D_j)$ is the
algebraic multiplicity of 0 as an eigenvalue of $D$.

Finally, note that, by \refr{det-eta}.c, if the dimension of $M$ is odd, and all the operators $D_j$ are invertible (so that $m_0(D_j)=0$), then
\refe{grdet-eta} takes the form
\eq{grdet-eta-odd}
    \LDgrtetnp(D) \ = \ \frac12\,\sum_{j=0}^d\, (-1)^j\,\LDnp_{2\tet}(D^2_j)  \ - \
   i\pi\,\eta(D).
\end{equation}

\subsection{Generalization}\Label{SS:general}
The definition \refe{etainv} of the $\eta$-invariant easily generalizes to operators acting on a subspace of the space $C^\infty(M,E)$ of smooth
sections of the vector bundle $E$, cf. \cite[\S4.10]{BrKappelerRAT}.

Let $D:C^\infty(M,E)\to C^\infty(M,E)$ be an elliptic differential operator with a self-adjoint leading symbol. Let $Q:C^\infty(M,E)\to
C^\infty(M,E)$ be a 0-th order pseudo-differential projection commuting with $D$. Then $V:= \IM{Q}\subset C^\infty(M,E)$ is a $D$-invariant subspace.
Let $\Pi_>$  and $\Pi_<$ be as in \refd{eta}. Let $\tet$ be as in \refss{eta-grdet} and set
\begin{gather}
    \eta_\tet(s,D|_V) \ = \ \zet_\tet(s,Q\Pi_>,D) \ - \ \zet_\tet(s,Q\Pi_<,-D), \Label{E:etasV}\\
    \eta(D|_V) \ = \ \frac{\eta_\tet(0,D|_V)+m_+(D|_V)-m_-(D|_V)+m_0(D|_V)}2.  \Label{E:etaV}
\end{gather}
Then, cf. \cite[\S4.10]{BrKappelerRAT},
\eq{det-etaV}
   \LD_{\tet}(D|_V) \ = \ \frac12\,\LD_{2\tet}(D^2|_V)  \ - \
   i\pi\,\Big(\,\eta(D|_V)- \frac{\zet_{2\tet}(0,D^2|_V)+m_0(D|_V)}2\,\Big),
\end{equation}
where we used the notation
\eq{zetDV}
    \zet_{2\tet}(s,D^2|_V) \ = \ \zet_{2\tet}(s,Q,D^2),
\end{equation}
cf. \refe{zetaQ}.

Note, however, that an analogue of \refe{det-eta-odd} does not necessarily hold in this case even if $\dim{M}$ is odd, because
$\zet_{2\tet}(s,D^2|_{V_j})$ defined by \refe{zetDV}, is not a $\zet$-function of a {\em differential} operator and \refp{zet=0} does not necessarily
hold.

Finally, suppose that $V=\bigoplus_{j=0}^dV_j$ is given as in \refd{grdeterminantV}. Then
\meq{grdet-etaV}
    \LDgrtet(D) \ = \ \frac12\,\sum_{j=0}^d\, (-1)^j\,\LD_{2\tet}(D^2|_{V_j})  \\ - \
   i\pi\,\Big(\,\eta(D|_V)-\frac{m_0(D|_V)}2- \frac12\sum_{j=0}^d\,(-1)^j\zet_{2\tet}(0,D^2|_{V_j})\,\Big),
\end{multline}
where $\eta(D|_V)= \sum_{j=0}^d\eta(D|_{V_j})$ and $m_0(D|_V)= \sum_{j=0}^dm_0(D|_{V_j})$.

\section{The Graded Determinant of the Odd Signature Operator}\Label{S:grdet}

In this section we define the graded determinant of the Atiyah-Patodi-Singer odd signature operator, \cite{APS2,Gilkey84}, of a flat vector bundle
$E$ over a closed oriented Riemannian manifold $M$. We also use this determinant to define an element $\rho$ of the determinant line of the
cohomology of the bundle $E$. Our definition is based on the formula which relates the graded determinant of the signature operator and the refined
torsion in the finite dimensional setting, cf. \refp{antor-grdetfd2}. In \refs{Ray-Singer} we will show that, if $E$ admits an invariant Hermitian
metric, then the Ray-Singer norm of $\rho$ is equal to 1. Thus the element $\rho$ can be viewed as  a {\em refinement of the Ray-Singer metric}. In
general, however, it might depend on the Riemannian metric on $M$. In subsequent sections we study the {\em metric anomaly} of $\rho$, use it to
``correct" $\rho$, and then define a differential invariant of the flat bundle $E$ -- a metric independent element of the determinant line of the
cohomology, called the {\em refined analytic torsion}.

\subsection{Setting}\Label{SS:setgrdet}
Let $M$  be a smooth closed oriented manifold of odd dimension $d=2r-1$ and let $E\to M$ be a complex vector bundle over $M$ endowed with a flat
connection $\n$. We denote by $\n$ also the induced differential
\[
    \n:\, \Ome^\b(M,E) \ \longrightarrow \Ome^{\b+1}(M,E),
\]
where $\Ome^k(M,E)$ denotes the space of smooth differential forms on $M$  of degree $k$ with values in $E$.

\subsection{Odd signature operator}\Label{SS:oddsign}
Fix a Riemannian metric $g^M$ on $M$ and let $*:\Ome^\b(M,E)\to \Ome^{d-\b}(M,E)$ denote the Hodge $*$-operator. Define the {\em chirality operator}
$\Gam= \Gam(g^M):\Ome^\b(M,E)\to \Ome^\b(M,E)$ by the formula
\eq{Gam}
    \Gam\, \ome \ := \ i^r\,(-1)^{\frac{k(k+1)}2}\,*\,\ome, \qquad \ome\in \Ome^k(M,E),
\end{equation}
with $r$ given as above by \/ $r=\frac{d+1}2$. This operator is equal to the operator defined in \S3.2 of \cite{BeGeVe} as one can see by applying
Proposition~3.58 of \cite{BeGeVe} in the case $\dim{M}$ is odd. In particular, $\Gam^2=1$.

\defe{oddsign}
The {\em odd signature operator} is the operator
\begin{equation} \Label{E:oddsignGam}
    \B=\B(\n,g^M) \ := \ \Gam\,\n \ + \ \n\,\Gam:\,\Ome^\b(M,E)\ \longrightarrow \  \Ome^\b(M,E).
\end{equation}
We denote by $\B_k$ the restriction of $\B$ to the space $\Ome^{k}(M,E)$.
\edefe
More explicitly, the value of the odd signature operator  on a form $\ome\in \Ome^{k}(M,E)$ is given by the formula
\begin{equation} \Label{E:oddsign}
    \B_k\,\ome \ := \ i^r(-1)^{\frac{k(k+1)}2+1}\,\big(\,(-1)^k*\n-\n*\,\big)\,\ome
    \ \in \ \Ome^{d-k-1}(M,E)\oplus  \Ome^{d-k+1}(M,E),
\end{equation}

The odd signature operator was introduced by Atiyah, Patodi, and Singer, \cite[p.~44]{APS1}, \cite[p.~405]{APS2}, in the case when $E$ is endowed
with a Hermitian metric {\em invariant with respect to $\n$} (i.e. invariant under parallel transport by $\n$). The general case was studied by
Gilkey, \cite[p.~64--65]{Gilkey84}.

\lem{symbolofB}
Suppose that $E$ is endowed with a Hermitian metric $h^{E}$. Denote by $\<\cdot,\cdot\>$ the scalar product on $\Ome^\b(M,E)$ induced by $h^{E}$ and
the Riemannian metric $g^M$ on $M$. Then

1. \  $\B$ is elliptic and its leading symbol is symmetric with respect to the Hermitian metric $h^{E}$.

2. \ If, in addition, the metric $h^{E}$ is invariant with respect to the connection $\n$, then $\B$ is symmetric with respect to the scalar product
$\<\cdot,\cdot\>$, $$\B^*= \B.$$ If\, the metric $h^{E}$ is not invariant, then, in general, $\B$ is not symmetric.
\elem
The proof of the lemma is a simple calculation. The first part is already stated in \cite[p.~405]{APS2}. The second part is proven in the Remark on
page~65 of \cite{Gilkey84}.

\subsection{Decomposition of the odd signature operator}\Label{SS:decompos}
Set
\begin{equation}\Label{E:decompos}\notag
 \begin{aligned}
    {}&\Ome^{\even}(M,E) \ := \ \bigoplus_{p=0}^{r-1}\, \Ome^{2p}(M,E), \quad
    \Ome^{\odd}(M,E) \ := \ \bigoplus_{p=1}^r\, \Ome^{2p-1}(M,E),\\
    {}&\B_{\even} \ := \ \bigoplus_{p=0}^{r-1}\, \B_{2p}:\, \Ome^{\even}(M,E)\ \longrightarrow \ \Ome^{\even}(M,E), \\
    {}&\B_{\odd} \ := \ \bigoplus_{p=1}^{r}\, \B_{2p-1}:\, \Ome^{\odd}(M,E)\ \longrightarrow \ \Ome^{\odd}(M,E).
 \end{aligned}
\end{equation}

Since $\Gam^2=1$ we obtain
\eq{B=*B*}
    \B_{\odd}\ =\ \Gam\circ \B_{\even}\circ \Gam{\big|_{\Ome^\odd(M,E)}}.
\end{equation}
Hence, the whole information about the odd signature operator is encoded in its {\em even} part $\B_{\even}$. The operator $\B_{\even}$ can be
expressed by the following formula, which is slightly simpler than \refe{oddsign}:
\eq{oddsigneven}
    \B_{\even }\,\ome \ := \ i^r(-1)^{p+1}\,\big(\,*\n-\n*\,\big)\,\ome,\qquad\text{for}\quad
                 \ome\in \Ome^{2p}(M,E).
\end{equation}

Note that for each $k=0\nek d$, the operator $\B^2$ maps $\Ome^k(M,E)$ into itself. Suppose $\calI$ is an interval of the form $[0,\lam],\
(\lam,\mu]$, or $(\lam,\infty]$ ($\mu\ge\lam\ge0$). Then $\Pi_{\B^2,\calI}$ is the spectral projection of $\B^2$ corresponding to $\calI$, cf.
\refss{spectralsubspace}.  Set
\eq{OmecalI}\notag
    \Ome^\b_{\calI}(M,E)\ := \ \Pi_{\B^2,\calI}\big(\, \Ome^\b(M,E)\,\big)\ \subset\  \Ome^\b(M,E).
\end{equation}
Recall from \refss{spectralsubspace}, that if the interval $\calI$ is bounded, then the space $\Ome^\b_\calI(M,E)$ is finite dimensional and is equal
to the span of the generalized eigenforms of $\B^2$ corresponding to eigenvalues with absolute value $\le\lam$. In general, $\Ome^\b_\calI(M,E)$
contains the span of the eigenforms of $\B^2$ corresponding to eigenvalues whose absolute value lies in $\calI$.

For each $k=0\nek d$, set
\eq{ome+-}
 \begin{aligned}
  \Ome^k_{+,\calI}(M,E) \ &:= \ \Ker\,(\n\,\Gam)\,\cap\,\Ome^k_{\calI}(M,E) \ = \ \big(\,\Gam\,(\Ker\,\n)\,\big)\,\cap\,\Ome^k_{\calI}(M,E);\\
  \Ome^k_{-,\calI}(M,E) \ &:= \ \Ker\,(\Gam\,\n)\,\cap\,\Ome^k_{\calI}(M,E)
  \ = \ \Ker\,\n\,\cap\,\Ome^k_{\calI}(M,E).
 \end{aligned}
\end{equation}
Clearly,
\eq{Ome>0=directsum}
    \Ome^k_{\calI}(M,E) \ = \ \Ome^k_{+,\calI}(M,E) \,\oplus \Ome^k_{-,\calI}(M,E), \qquad\text{if}\quad 0\not\in \calI.
\end{equation}
We consider the decomposition \refe{Ome>0=directsum} as a {\em grading} \footnote{Note, that our grading is opposite to the one considered in
\cite[\S2]{BFK3}.} of the space $\Ome^\b_{\calI}(M,E)$, and refer to $\Ome^k_{+,\calI}(M,E)$ and $\Ome^k_{-,\calI}(M,E)$ as the positive and negative
subspaces of $\Ome^k_{\calI}(M,E)$.

As both, $\Gam$ and $\n$, commute with $\B^2$, we conclude that, for $k=0\nek d$,
\eq{GamOme}
    \Gam:\,\Ome^k_{\pm,\calI}(M,E) \ \longrightarrow \ \Ome^{d-k}_{\mp,\calI}(M,E),
\end{equation}
and
\eq{naOme}
    \n:\,\Ome^k_{\pm,\calI}(M,E) \ \longrightarrow \ \Ome^{k+1}_{\mp,\calI}(M,E).
\end{equation}
Denote by $\B_k^\calI$ the restriction of $\B$ to $\Ome^k_\calI(M,E)$ and by $\B_{k}^{\pm,\calI}$ the restriction of $\B$ to
$\Ome^k_{\pm,\calI}(M,E)$. Then
\[
 \begin{aligned}
    \B_k^{+,\calI}:\, \Ome^k_{+,\calI}(M,E) \ &\longrightarrow \ \Ome^{d-k-1}_{+,\calI}(M,E), \qquad \ome\mapsto \Gam\n\ome;\\
    \B_k^{-,\calI}:\, \Ome^k_{-,\calI}(M,E) \ &\longrightarrow \ \Ome^{d-k+1}_{-,\calI}(M,E), \qquad \ome\mapsto \n\Gam\ome.
 \end{aligned}
\]

\subsection{Graded determinant of the odd signature operator}\Label{SS:drdetoddsign}
Let $\calI$ be an interval of the form $[0,\lam]$, $(\lam,\mu]$, or $(\lam,\infty]$ ($\mu\ge\lam\ge0$) and define
\[
    \Ome^{\even}_{\pm,\calI}(M,E)\ =\ \bigoplus_{p=0}^{r-1}\, \Ome^{2p}_{\pm,\calI}(M,E).
\]
Let $\B^\calI$, $\B^\calI_\even$, and $\B^\calI_\odd$ denote the restrictions of $\B$ to the subspaces $\Ome^\b_{\calI}(M,E)$,
$\Ome^\even_{\calI}(M,E)$, and $\Ome^\odd_{\calI}(M,E)$ respectively. Then
\[
    \B^\calI_{\even}:\, \Ome^{\even}_{\pm,\calI}(M,E)\ \longrightarrow \  \Ome^{\even}_{\pm,\calI}(M,E).
\]
Let $\B_{{\even}}^{\pm,\calI}$ denote the restriction of $\B^\calI_{{\even}}$ to the space $\Ome^{\even}_{\pm,\calI}(M,E)$. Clearly, the operators
$\B_{{\even}}^{\pm,\calI}$ are bijective whenever $0\not\in \calI$.


By Definitions~\ref{D:detgrfd} and \ref{D:grdeterminantV}, for every $\calI$, the graded determinant of $\B^{\calI}_{\even}$ is given by the formula
\eq{grdetB}
    \Detgrtet(\B^{\calI}_{\even}) \ := \ e^{\LDgrtet(\B^{\calI}_{\even})},
\end{equation}
where $\tet\in (-\pi,0)$ is an Agmon angle for the operator $\B^{\calI}_{\even}$, and
\eq{grldetB}
    \LDgrtet(\B^{\calI}_{\even}) \ := \ \LD_\tet\big(\B_{\even}^{+,{\calI}}\big)
      \ - \ \LD_\tet\big(-\B_{\even}^{-,{\calI}}\big)
    \ \in \ \CC.
\end{equation}

Clearly, for $0\le \lam\le \mu$, we have
\eq{DetlamDetmu}
    \Detgrtetnp(\B^{(\lam,\infty)}_\even) \ = \ \Detgrtetnp(\B^{(\lam,\mu]}_\even)\cdot \Detgrtetnp(\B^{(\mu,\infty)}_\even).
\end{equation}
Note also that since the rank of $\B^{(\lam,\mu]}_\even$ is finite, $\Detgrtetnp(\B^{(\lam,\mu]}_\even)$ is independent of $\tet$ and is equal to the
product of the eigenvalues of $\B^{(\lam,\mu]}_\even$.

\subsection{The canonical element of the determinant line}\Label{SS:rho}
Since the covariant differentiation $\n$ commutes with $\B$, the subspace $\Ome^\b_\calI(M,E)$ is a subcomplex of the twisted de Rham complex
$(\Ome^\b(M,E),\n)$. Clearly, for each $\lam\ge0$, the complex $\Ome^\b_{(\lam,\infty)}(M,E)$ is acyclic. Since
\eq{Ome=Ome0+Ome>0}
    \Ome^\b(M,E) \ = \  \Ome^\b_{[0,\lam]}(M,E)\,\oplus\,\Ome^\b_{(\lam,\infty)}(M,E),
\end{equation}
the cohomology $H^\b_{[0,\lam]}(M,E)$ of the complex $\Ome^\b_{[0,\lam]}(M,E)$ is naturally isomorphic to the cohomology $H^\b(M,E)$ of
$\Ome^\b(M,E)$.

Let $\Gam_{\hskip-1pt{}\calI}$ denote the restriction of $\Gam$ to $\Ome^\b_\calI(M,E)$. For each $\lam\ge0$, let
\eq{thoGam0lam}
    \rho_{{}_{\Gam_{\hskip-1pt{}_{[0,\lam]}}}}\ = \ \rho_{{}_{\Gam_{\hskip-1pt{}_{[0,\lam]}}}}(\n,g^M) \in \ \Det(H^\b_{[0,\lam]}(M,E))
\end{equation}
denote the refined torsion of the finite dimensional complex $(\Ome^\b_{[0,\lam]}(M,E),\n)$ corresponding to the chirality operator
$\Gam_{\hskip-1pt{}_{[0,\lam]}}$, cf. \refd{refinedtorsion}. We view $\rho_{{}_{\Gam_{\hskip-1pt{}_{[0,\lam]}}}}$ as an element of $\Det(H^\b(M,E))$
via the canonical isomorphism between $H^\b_{[0,\lam]}(M,E)$ and $H^\b(M,E)$.

From \refp{antor-grdetfd2}, \refe{DetlamDetmu}, and \refe{det-det}, we immediately obtain the following
\prop{antor-grdet}
Assume that $\tet\in (-\pi,0)$ is an Agmon angle for the operator $\B_\even$. Then the element
\eq{rho}
    \rho \ = \ \rho(\n,g^M) \ := \ \Detgrtetnp(\B^{(\lam,\infty)}_{\even})\cdot \rho_{{}_{\Gam_{\hskip-1pt{}_{[0,\lam]}}}}
    \ \in  \ \Det(H^\b(M,E))
\end{equation}
is independent of the choice of $\lam\ge0$. Further, $\rho$ is independent of the choice of the Agmon angle $\tet\in (-\pi,0)$ of\/ $\B_\even$.
\eprop

If the odd signature operator is invertible then $H^\b(M,E)=0$. In this case, $\Det(H^\b(M,E))$ is canonically isomorphic to $\CC$ and
$\rho_{{}_{\Gam_{\hskip-1pt{}_{\{0\}}}}}= 1$. Hence, $\rho$ is a complex number which coincides with the graded determinant
$\Det_{\operatorname{gr},\tet}(\B_{\even})= \Det_{\operatorname{gr},\tet}(\B^{(0,\infty)}_{\even})$. This case was studied in \cite{BrKappelerRAT}.

\section{Relationship with the $\eta$-invariant}\Label{S:witheta}

In this section, we study the relationship between the graded determinant \refe{grdetB} and the $\eta$-invariant of $\B^{(\lam,\infty)}_\even$. For
the special case when $\B$ is bijective and $\lam=0$ this relationship was established in Section~7 of \cite{BrKappelerRAT}.

To simplify the notation set
\eq{eta=eta}
    \eta_\lam \ = \ \eta_\lam(\n,g^M) \ := \  \eta(\B_{\even}^{(\lam,\infty)}),
\end{equation}
and
\eq{xi}
\begin{aligned}
    \xi_\lam \ = \ \xi_\lam(\n,g^M,\tet) \ &:= \ \frac12\,\sum_{j=0}^{d-1}\,(-1)^j\LDnp_{2\tet}\big(\,\B_{d-j-1}^{+,{(\lam,\infty)}}\circ \B_{j}^{+,{(\lam,\infty)}}\,\big)
    \\ &= \ \frac12\,\sum_{j=0}^{d-1}\,(-1)^j\LDnp_{2\tet} \Big(\,
    {(\Gam\n)^2}_{\big|_{\Ome^j_{+,{(\lam,\infty)}}(M,E)}}\,\Big),
\end{aligned}
\end{equation}
where $\tet\in (-\pi/2,0)$ and both, $\tet$ and $\tet+\pi$, are Agmon angles for $\B_\even$ (hence, $2\tet$ is an Agmon angle for $\B_\even^2$).

It is shown in Section~8.4 of \cite{BrKappelerRAT} that \footnote{In \cite{BrKappelerRAT} we only considered the case when $\B$ is bijective and
$\lam=0$. But the arguments leading to formula (8.7) of \cite{BrKappelerRAT} work without any changes in our more general situation.}
\eq{xi2}
\begin{aligned}
    \xi_\lam  \ &:= \ \frac12\,\sum_{j=0}^{d}\,(-1)^{j+1}\,j\,\LDnp_{2\tet}\big[\,{(\B^{(\lam,\infty)})^2}_{\big|_{\Ome^j_{{(\lam,\infty)}}(M,E)}}\,\big]
    \\ &= \ \frac12\,\sum_{j=0}^{d}\,(-1)^{j+1}\,j\,\LDnp_{2\tet} \Big[\,\Big(\,
    (\Gam\n)^2+(\n\Gam)^2\,\Big)_{\big|_{\Ome^j_{{(\lam,\infty)}}(M,E)}}\,\Big].
\end{aligned}
\end{equation}

Set
\eq{djlam}
    d_{j,\lam} \ := \ \dim \Ome^j_{[0,\lam]}(M,E),
    \qquad j=0\nek d.
\end{equation}

\prop{xietad}
Let $\n$ be a flat connection on a vector bundle $E$ over a closed Riemannian manifold $(M,g^M)$ of odd dimension $d=2r-1$. Assume $\tet\in
(-\pi/2,0)$ is an Agmon angle for the odd signature operator $\B= \B(\n,g^M)$ such that there are no eigenvalues of\/ $\B$ in the solid angles
$L_{(-\pi/2,\tet]}$ and $L_{(\pi/2,\tet+\pi]}$. Then, for every $\lam\ge0$,
\eq{xietad}
    \LDgrtetnp(\B_\even^{(\lam,\infty)}) \ =\ \xi_\lam \ - \ i\pi\,\eta_\lam \ - \ \frac{i\pi}2\,\sum_{j=0}^d\,(-1)^jj\,d_{j,\lam}.
\end{equation}
\eprop
\prf
Since the operator $B_\even^{(\lam,\infty)}$ has no zero eigenvalues, we conclude from \refe{grdet-etaV}, that to prove \refe{xietad} it is enough to
show the following two identities
\begin{gather}
    2\,\xi_\lam \ = \ \LDnp_{2\tet}\,(\B_{\even}^{+,(\lam,\infty)})^2 \ - \ \LDnp_{2\tet}\,(\B_{\even}^{-,(\lam,\infty)})^2; \Label{E:calBalp=}\\
    \zet_{2\tet}\big(0,(\B_{\even}^{+,(\lam,\infty)})^2\big) \ - \ \zet_{2\tet}\big(0,(\B_{\even}^{-,(\lam,\infty)})^2\big)
    \ = \ \sum_{j=0}^d\,(-1)^{j}j\,d_{j,\lam}. \Label{E:zetB2=}
\end{gather}
A verbatim repetition of the arguments which led to formula (7.17) of \cite{BrKappelerRAT} implies that
\eq{zetB2=3}
    \zet_{2\tet}\big(s,(\B_{\even}^{+,(\lam,\infty)})^2\big) \ - \ \zet_{2\tet}\big(s,(\B_{\even}^{-,(\lam,\infty)})^2\big) \ = \
    \sum_{j=0}^d\, (-1)^{j+1}\,j\,
    \zet_{2\tet}\Big(\,s, {(\B^{(\lam,\infty)})^2}_{\big|_{\Ome^j(M,E)}}\, \Big).
\end{equation}

From \refe{xi2} and \refe{zetB2=3} we obtain
\eq{calBalp=2}
    2\,\xi_\lam \ = \ \frac{d}{ds}{\big|_{s=0}}\Big[\,
      \zet_{2\tet}\big(s,(\B_{\even}^+)^2\big) \ - \ \zet_{2\tet}\big(s,(\B_{\even}^-)^2\big) \,\Big].
\end{equation}
Hence \refe{calBalp=} is established.

Combining \refe{zetB2=3} and \refp{zet=0} we obtain \refe{zetB2=}.
\eprf
\section{The Metric Anomaly of $\rho$ and the Definition of the Refined Analytic Torsion}\Label{S:metric}

In this section we study the dependence of the element $\rho=\rho(\n,g^M)$ defined in \refe{rho} on the Riemannian metric $g^M$. In particular, we
show that, if $\dim M= 2r-1\equiv 1 \ (\MOD\ 4)$, then $\rho$ is independent of $g^M$. We then use these results to construct the {\em refined
analytic torsion} -- a canonical element of $\Det\big(H^\b(M,F)\big)$ which is independent of the metric, i.e., is a differential invariant of the
flat vector bundle $(E,\n)$.

\subsection{Relationship between $\rho(t)$ and the $\eta$-invariant}\Label{SS:rho-eta}
Suppose that $g^M_t$, $t\in\RR$, is a smooth family of Riemannian metrics on $M$. Let
\eq{rho(t)}
    \rho(t) \ = \ \rho(\n,g^M_t) \ \in \ \Det\big(\,H^\b(M,E)\,\big)
\end{equation}
be the canonical element defined in \refe{rho}.

Let $\Gam_t$ denote the chirality operator corresponding to the metric $g_t^M$, cf. \refe{Gam}, and let $\B(t)= \B(\n,g^M_t)$ denote the odd
signature operator corresponding to the Riemannian metric $g^M_t$.

Fix $t_0\in \RR$ and choose $\lam\ge 0$ so that there are no eigenvalues of $\B(t_0)^2$ whose absolute values are equal to $\lam$. Then there exists
$\del>0$ such that the same is true for all $t\in (t_0-\del,t_0+\del)$. In particular, if we denote by $\Ome^\b_{[0,\lam],t}(M,E)$ the span of the
generalized eigenvectors of $\B(t)^2$ corresponding to eigenvalues with absolute value $\le\lam$, then $\dim\Ome^\b_{[0,\lam],t}(M,E)$ is independent
of $t\in (t_0-\del,t_0+\del)$. We set
\eq{d(t)}
    d_{j,\lam} \ := \ \dim \Ome^j_{[0,\lam],t}(M,E), \qquad j=0\nek d, \quad t\in (t_0-\del,t_0+\del).
\end{equation}

By definition \refe{rho},
\eq{rho(t)=}
    \rho(t) \ = \ \Detgrtetnp\big(\B_\even^{(\lam,\infty)}(t)\big)\cdot\rho_{{}_{\Gam_{\hskip-1pt{}t,[0,\lam]}}}.
\end{equation}

For each $t\in (t_0-\del,t_0+\del)$ and $\tet\in (-\pi/2,0)$, such that $\tet$ and $\tet+\pi$ are Agmon angles for $\B^{(\lam,\infty)}(t)$, let us
introduce the following short notation for the quantities introduced in \refe{eta=eta} and \refe{xi}
\eq{xi(t)eta(t)}\notag
    \xi_\lam(t,\tet) \ := \ \xi_\lam(\n,g^M_t,\tet), \qquad \eta_\lam(t) \ := \ \eta_\lam(\n,g^M_t).
\end{equation}

Assume that $\tet_0\in (-\pi/2,0)$ is an Agmon angle for $\B(t_0)= \B(\n,g_{t_0}^M)$ such that there are no eigenvalues of $\B(t_0)$ in the solid
angles $L_{(-\pi/2,\tet_0]}$ and $L_{(\pi/2,\tet_0+\pi)}$. Choose $\del$, if necessary, smaller, so that for every $t\in (t_0-\del,t_0+\del)$ and
every $j=0\nek d$ both, $\tet_0$ and $\tet_0+\pi$, are Agmon angles of $\B^{(\lam,\infty)}_j(t)$. For $t\not= t_0$ it might happen that there are
eigenvalues of $\B^{(\lam,\infty)}(t)$ in $L_{(-\pi/2,\tet_0)}$ and/or $L_{(\pi/2,\tet_0+\pi)}$. Hence, \refe{xietad} is not necessarily true, in
general, for $t\not=t_0$. However, from \refe{zet-zet} and \refe{xi}, we conclude that for every angle $\tet\in (-\pi/2,0)$, so that $\tet$ and
$\tet+\pi$ are Agmon angles for $\B^{(\lam,\infty)}(t)$ (and, hence, $2\tet$ is an Agmon angle for $\B^{(\lam,\infty)}(t)^2$),
\eq{xitettet'}
    \xi_\lam(t,\tet) \ \equiv \ \xi_\lam(t,\tet_0) \qquad \MOD\ \pi i,
\end{equation}
Thus, from \refe{xietad}, we obtain
\eq{rho(t)=2}
    \rho(t) \ = \ \pm \, e^{\xi_\lam(t,\tet_0)}\cdot e^{-i\pi\eta_\lam(t)}\cdot
    e^{-\frac{i\pi}2\,\sum_{j=0}^d\,(-1)^jj\,d_{j,\lam}}\cdot\rho_{{}_{\Gam_{\hskip-1pt{}t,[0,\lam]}}}.
\end{equation}

\lem{xirho}
Under the above assumptions, the product\/ $e^{\xi_\lam(t,\tet_0)}\cdot\rho_{{}_{\Gam_{\hskip-1pt{}t,[0,\lam]}}}\in \Det\big(H^\b(M,F)\big)$ is
independent of\/ $t\in (t_0-\del,t_0+\del)$.
\elem
\prf
Recall that we have chosen $\lam\ge0$ and $\del>0$ so that there are no eigenvalues of $\calB(t)^2$ with absolute value $\lam$ for any $t\in
(t_0-\del,t_0+\del)$.

We shall use the following notation (cf., for example, Section~2 of \cite{BFK3}): Suppose $f(s)$ is a function of a complex parameter $s$ which is
meromorphic near $s=0$. We call the zero order term in the Laurent expansion of $f$ near $s=0$ the {\em finite part of $f$ at 0} and denote it by
$\Fp f(s)$. A verbatim repetition of the proof of formula (9.13) of \cite{BrKappelerRAT} shows that
\eq{ddtxi}
    \frac{d}{dt}\, \xi_\lam(t,\tet_0) \ = \ \frac12\,
    \sum_{j=0}^d\, (-1)^j\,\Fp\,\Tr\,\Big[\,\dot\Gam_t\,\Gam_t\,
        {\Big(\,\big(\B^{(\lam,\infty)}(t)\big)^2\,\Big)^{-s}_{2\tet_0}}_{\Big|_{\Ome^j_{(\lam,\infty)}(M,E)}}\,\Big].
\end{equation}
Since $\B^2(t)$ is an elliptic differential operator, the dimension of $\Ome^\b_{[0,\lam]}(M,E)$ is finite. Let $\eps\not=0$ be a small enough real
number so that $\B(t)^2+\eps$ is bijective and $2\tet_0$ is an Agmon angle for $(\B^{(\lam,\infty)}(t))^2+\eps$. Then, for each $j=0\nek d$, we have
\meq{Fplam=Fp+d}
    \Fp\,\Tr\,\Big[\,\dot\Gam_t\,\Gam_t\,
        {\Big(\,\big(\B^{(\lam,\infty)}(t)\big)^2\,\Big)^{-s}_{2\tet_0}}_{\Big|_{\Ome^j_{(\lam,\infty)}(M,E)}}\,\Big]
    \\ = \
    \Fp\,\Tr\,\Big[\,\dot\Gam_t\,\Gam_t\,{\Big(\,\big(\B^{(\lam,\infty)}(t)\big)^2+\eps\,\Big)^{-s}_{2\tet_0}}_{\Big|_{\Ome^j(M,E)}}\,\Big]
    \ - \ \Tr\,\Big[\,{\dot\Gam_t\,\Gam_t}_{\big|_{\Ome^j_{[0,\lam]}(M,E)}}\,\Big].
\end{multline}
By a slight generalization of a result of Seeley \cite{Seeley67}, which is discussed in \cite{RaSi1} and, in greater generality, in
\cite{Grubbseeley95}, the first summand on the right hand side of \refe{Fplam=Fp+d} is given by a {\em local formula}, i.e., by the integral
$\int_M\phi_t$ of a differential form $\phi_t$, whose value at any point $x\in M$ depends only on the values of the components of the metric tensor
$g_t^M$ and a finite number of their derivatives at $x$. Moreover, since the dimension of the manifold $M$ is odd, the differential form $\phi_t$
vanishes identically. Thus we obtain from \refe{ddtxi} and \refe{Fplam=Fp+d}
\eq{ddtxi2}
    \frac{d}{dt}\, \xi_\lam(t,\tet_0) \ = \ -\,\frac12\,\sum_{j=0}^d\,(-1)^j\,\Tr\,\big[\,{\dot\Gam_t\,\Gam_t}_{\big|_{\Ome^j_{[0,\lam]}(M,E)}}\,\big].
\end{equation}
Combining this equation with \refp{deponGam}, we get
\[
    \frac{d}{dt}\,e^{\xi_\lam(t,\tet_0)}\cdot\rho_{{}_{\Gam_{\hskip-1pt{}t,[0,\lam]}}}\ = \ 0.
\]
\eprf

\subsection{Dependence of the $\eta$-invariant on the metric}\Label{SS:etamet}
From \refe{rho(t)=2} and \refl{xirho} we see that the dependence of $\rho(t)$ on $t\in (t_0-\del,t_0+\del)$ is determined up to a sign  by the
dependence of $\eta_\lam(t)$ on $t$.

\lem{etalam-eta}
For any $t_1,t_2\in (t_0-\del,t_0+\del)$ we have
\eq{etalam-eta}
    \eta_\lam(t_1) \ - \ \eta_\lam(t_2) \ \equiv \ \eta\big(\B_\even(t_1)\big) \ - \ \eta\big(\B_\even(t_2)\big), \qquad \MOD\ \ZZ.
\end{equation}
\elem
\prf
Recall that $\Ome^\b_{[0,\lam],t}(M,E)$ denotes the span of the generalized eigenvectors of $\B(t)^2$ corresponding to eigenvalues with absolute
value $\le\lam$ and that $\lam$ and $\del$ were chosen so that
\eq{dimOme(t)}
    \dim \Ome^j_{[0,\lam],t}(M,E) \ = \ \const, \qquad t\in (t_0-\del,t_0+\del), \ j=0\nek d.
\end{equation}
Since the dimension of the space $\Ome^\even_{[0,\lam],t}(M,E)$ is finite, formula \refe{etasV} says that $\eta\big(0,\B_\even^{[0,\lam]}(t)\big)$ is
equal to the sum of the algebraic multiplicities of the eigenvalues of $\B_\even^{[0,\lam]}(t)$ with positive real parts minus the sum of the
algebraic multiplicities of the eigenvalues of $\B_\even^{[0,\lam]}(t)$ with  negative real parts. It then follows from \refe{etaV} that
\eq{eta0lam}
    \eta\big(\B_\even^{[0,\lam]}(t)\big) \ \equiv \ \frac12\,\dim \Ome^\even_{[0,\lam],t}(M,E), \qquad \MOD\ \ZZ.
\end{equation}
By the definition of the $\eta$-invariant,
\[
    \eta\big(\B_\even(t)\big) \ - \ \eta_\lam(t) \ = \ \eta\big(\B_\even^{[0,\lam]}(t)\big).
\]
Hence, from \refe{dimOme(t)} and \refe{eta0lam}, we conclude that, modulo $\ZZ$,
\[
    \eta\big(\B_\even(t)\big) \ - \ \eta_\lam(t) \ \equiv \ \const
\]
for $t\in (t_0-\del,t_0+\del)$.
\eprf

We now need to study the dependence of $\eta(\B_\even(\n,g^M))$ on the Riemannian metric $g^M$. Fortunately, this was essentially done in \cite{APS2}
and \cite{Gilkey84}. Below we present a brief review of the relevant results.

Let $\B_{\text{trivial}}= \B_{\text{trivial}}(g^M):\Ome^{\even}(M)\to \Ome^{\even}(M)$ denote the even part of the odd signature operator
corresponding to the metric $g^M$ and the trivial line bundle over $M$ endowed with the trivial connection.  It is shown on page~52 of
\cite{Gilkey84} (see also Theorem~2.4 of \cite{APS2} where the case of unitary connection is established) that modulo $\ZZ$ the difference
\eq{tileta}
    \eta\big(\B_{\even}(\n,g^M)\big) \ - \ (\rank E)\ \eta\big(\B_{\text{trivial}}(g^M)\big)
\end{equation}
is independent of the Riemannian metric.

Let us describe the dependence of $\eta\big(\B_{\text{trivial}}(g^M)\big)$ on the metric.

\subsubsection{Case when $M$ bounds an oriented manifold $N'$}\Label{SS:Mbounds}
Suppose, first, that $M$ is the oriented boundary of a smooth compact oriented manifold $N'$. Let $\sign(N')$ denote the signature of $N'$, cf.
\cite{APS1}. This is an integer defined in purely cohomological terms. In particular, it is independent of the metric. The signature theorem for
manifolds with boundary (cf. Theorem~4.14 of \cite{APS1} and Theorem~2.2 of \cite{APS2}) states that
\eq{signth}
    \sign(N') \ = \ \int_{N'}\, L(p,g^M)  \ - \ \eta\big(\B_{\text{trivial}}(g^M)\big),
\end{equation}
where $L(p,g^M)$ is the Hirzebruch $L$-polynomial in the Pontrjagin forms of any Riemannian metric on $N'$ which near $M$ is the product of $g^M$ and
the standard metric on the half-line. It follows from \refe{signth} that $\int_{N'}L(p,g^M)$ is independent of the choice of such a metric on $N'$.
Note that if $\dim{M}\equiv 1 \ (\MOD\, 4)$ then $L(p,g^M)$ does not have a term of degree $\dim{N'}$ and, hence, $\int_{N'}\,L(p,g^M)=0$.

Combining \refe{signth} with the metric independence of\, $\sign(N')$, \refe{etalam-eta}, and \refe{tileta}, we conclude that, modulo $\ZZ$,
\eq{eta-L}
    \eta_\lam(t) \ - \ (\rank E)\, \int_{N'}\, L(p,g^M_t)
\end{equation}
is independent of $t\in (t_0-\del,t_0+\del)$. It is also independent of the choice of $N'$ since for different choices of $N'$, the integral
$\int_{N'}L(p,g^M)$ differs by an integer.

\subsubsection{General case ($M$ does not necessarily bound an oriented manifold)}\Label{SS:Mnotbounds}
In general, there might be no smooth oriented manifold whose oriented boundary is diffeomorphic to $M$. However, since $\dim{}M$ is odd, there exists
an oriented manifold $N$ whose oriented boundary is the disjoint union of two copies of $M$ (with the same orientation), cf. \cite{Wall60},
\cite[Th.~IV.6.5]{Rudyak_book}. Then the same arguments as above show that modulo $\ZZ$
\eq{eta-L2}
    \eta_\lam(t) \ - \ \frac{\rank E}2\, \int_{N}\, L(p,g^M_t)
\end{equation}
is independent of $t\in (t_0-\del,t_0+\del)$. In particular, if $\dim{M}\equiv 1 \ (\MOD\, 4)$, then modulo $\ZZ$, $\eta_\lam(t)$ is independent of
$t$.

Note, however, that replacing $\eta_\lam(t)$ by \refe{eta-L2} removes the dependence of the metric, but creates a new dependence on the choice of the
manifold $N$. For different choices of $N$, the integrals $\int_{N}\, L(p,g^M_t)$ might differ by an integer. Thus the expression \refe{eta-L2} is
well defined as a function of $\n$ only modulo $\frac{\rank E}2\,\ZZ$.

\subsection{Removing the metric anomaly}\Label{E:removeanom}
We are now ready to state the main result of this section.

\th{metricindep}
Let $E$ be a flat vector bundle over a closed oriented odd-dimensional manifold $M$. Let $N$ be an oriented manifold whose oriented boundary is the
disjoint union of two copies of $M$. For each Riemannian metric $g^M$ on $M$ consider
\eq{metricindep}
    \rho(\n,g^M)\cdot e^{i\pi\,\frac{\rank E}2\,\int_N\,L(p,g^M)}\ \in \ \Det(H^\b(M,E)),
\end{equation}
where $\rho(\n,g^M)\in \Det(H^\b(M,E))$ is defined in \refe{rho} and $L(p,g^M)$ is the Hirzebruch $L$-polynomial in the Pontrjagin forms of any
Riemannian metric on $N$ which near $M$ is the product of $g^M$ and the standard metric on the half-line. Then the element \refe{metricindep} is
independent of $g^M$.

In particular, if\/ $\dim{}M\equiv 1 (\MOD 4)$, then $\int_NL(p,g^M)=0$ and, hence, $\rho(\n,g^M)$ is independent of $g^M$.
\eth
\prf
Let $g^M_t$, $(t\in \RR)$ be a family of Riemannian metrics on $M$. We shall use the notation of \refss{rho-eta}. From \refe{rho(t)=2} we obtain for
$t\in (t_0-\del,t_0+\del)$
\eq{rho(t)int}\notag
    \rho(t)\cdot e^{i\pi\,\frac{\rank E}2\,\int_N\,L(p,g^M_t)} \ = \ \pm\, e^{\xi_\lam(t)}\cdot e^{-i\pi\eta_\lam(t)}\cdot
    e^{-\frac{i\pi}2\,\sum_{j=0}^d\,(-1)^jj\,d_{j,\lam}}\cdot \rho_{{}_{\Gam_{\hskip-1pt{}t,[0,\lam]}}}\cdot e^{i\pi\,\frac{\rank
    E}2\,\int_N\,L(p,g^M_t)}.
\end{equation}
Combining this latter equality with \refl{xirho} and the metric independence of the reduction of \refe{eta-L2} modulo $\ZZ$ we conclude that for any
$t_1,t_2\in (t-\del,t+\del)$
\eq{rho(t)int2}
    \rho(t_1)\cdot e^{i\pi\,\frac{\rank E}2\,\int_N\,L(p,g^M_t)} \ = \ \pm\,\rho(t_2)\cdot e^{i\pi\,\frac{\rank E}2\,\int_N\,L(p,g^M_t)}
\end{equation}
Since the function $t\mapsto \rho(t)\cdot e^{i\pi\,\frac{\rank E}2\,\int_N\,L(p,g^M_t)}$ is continuous and non-zero the sign in the right hand side
of \refe{rho(t)int2} must be positive. Hence, the theorem is proven.
\eprf

\subsection{The definition of the refined analytic torsion}\Label{SS:defRAT}
We now arrive to the main definition of the paper.

\defe{refantor}
Let $E$ be a complex vector bundle over a closed oriented odd-dimensional manifold $M$ and let $\n$ be the flat connection on $E$. The {\em refined
analytic torsion} $\rat= \rat(\n)$ is the element of $\Det(H^\b(M,E))$ defined by \refe{metricindep}.
\edefe

\rem{TdependsonN}
The refined analytic torsion is independent of the angle $\tet\in (-\pi,0)$ and of the metric. But {\em it does depend on the choice of the manifold
$N$}. However, from the discussion at the end of \refss{Mnotbounds}, we conclude that $\rat(\n)$ is defined up to multiplication by
$i^{k\cdot\rank{E}}$ $(k\in \ZZ)$. If $\rank{}E$ is even then $\rat(\n)$ is well defined up to a sign, and if $\rank{}E$ is divisible by 4, then
$\rat(\n)$ is a well defined element of $\Det(H^\b(M,E))$.

(Here a quantity being well defined means that it depends only on $M$, $E$ and $\n$.)
\erem

\rem{Twhenbounds}
If $M$ is the oriented boundary of a smooth compact oriented manifold $N'$, one can define a version of the refined analytic torsion:
\begin{equation}\Label{E:refantor2}
    \rat'(\n)  \ := \
        \rho(\n,g^M)\cdot  \exp\Big(\,i\pi\cdot\rank{}E\,  \int_{N'}\, L(p,g^M)\,\Big).
\end{equation}
Note that the indeterminacy in the definition of $\rho'(\n)$ is smaller than the indeterminacy in the definition of $\rho(\n)$, cf. \refss{Mbounds}:
$\rho'(\n)$ is well defined up to a sign. If $\rank{}E$ is even, then $\rat'(\n)$ is a well defined element of $\Det(H^\b(M,E))$.
\erem

\section{A Duality Theorem for the Refined Analytic Torsion}\Label{S:duality}

In this section we establish a relationship between the refined analytic torsion corresponding to a flat connection and that of its dual. This result
is used in the next section in order to calculate the Ray-Singer norm of the refined analytic torsion, but it is also of independent interest.

\subsection{The dual connection}\Label{SS:n'}
Suppose $M$ is a closed oriented manifold of odd dimension $d=2r-1$. Let $E\to M$ be a complex vector bundle over $M$ and let $\n$ be a flat
connection on $E$. Fix a Hermitian metric $h^{E}$ on $E$. Denote by $\n'$ the {\em connection on $E$ dual to the connection $\n$}. It is defined by
the formula
\[
    dh^{E}(u,v) \ = \ h^{E}(\n u,v) \ + \ h^{E}(u,\n' v),
    \qquad u,v\in C^\infty(M,E).
\]

For $\ome_1, \ome_2\in \Ome^\b(M,E)$ of the form $\ome_i=s_i\otimes\chi_i$ with $s_i\in C^\infty(M,E),\ \chi_i\in \Ome^\b(M,\RR)$, define
\eq{OmeEOmeE'}
    h^E\big(\, (s_1\otimes\chi_1)\wedge(s_2\otimes\chi_2)\,\big) \ := \ h^E(s_1,s_2)\cdot\chi_1\wedge\chi_2.
\end{equation}
Then $h^E$ extends in a canonical way to a sesquilinear map
\eq{OmeEOmeE'2}
    h^E:\,\Ome^\b(M,E)\times \Ome^\b(M,E) \ \longrightarrow \ \Ome^\b(M,\CC).
\end{equation}
For each $j=0\nek d$, we then obtain a sesquilinear pairing
\eq{pairingEE'}
    \Ome^{j}(M,E)\times \Ome^{d-j}(M,E) \ \longrightarrow \ \CC, \quad (\ome_1,\ome_2) \ \mapsto \ \int_M\,h^E(\ome_1\wedge\ome_2).
\end{equation}

We denote by $E'$ the  flat vector bundle $(E,\n')$, referring to  $E'$ as the {\em dual of the flat vector bundle $E$}. The pairing
\refe{pairingEE'} induces a non-degenerate sesquilinear pairing
\eq{Poincarepair}\notag
    H^j(M,E')\otimes H^{d-j}(M,E) \ \longrightarrow \ \CC, \qquad j=0\nek d,
\end{equation}
and, hence, identifies $H^j(M,E')$ with the dual space of $H^{d-j}(M,E)$. Using the construction of \refss{dualgraded} (with $\tau:\CC\to \CC$ being
the complex conjugation) we thus obtain an anti-linear isomorphism
\eq{alpE-E'}
    \alp:\, \Det\big(H^\b(M,E)\big) \ \longrightarrow \ \Det\big(H^\b(M,E')\big).
\end{equation}

\subsection{The duality theorem}\Label{SS:dualityth}
Fix a compact oriented manifold $N$ whose boundary is diffeomorphic to two disjoint copies of $M$. From \refss{Mnotbounds}, we conclude that the
number
\eq{Reeta+int}
    \exp\Big(\,2i\pi\,\oeta(\n,g^M)-i\pi\,\rank E\,\int_{N}\, L(p,g^M)\,\Big),
\end{equation}
where $\oeta$ denotes the complex conjugate of $\eta$, is independent of the choice of the Riemannian metric $g^M$.

The main result of this section is the following {\em duality theorem}.
\th{dualityRAT}
Let\/ $E\to M$ be a complex vector bundle over a closed oriented odd-dimensional manifold $M$ and let\/ $\n$ be a flat connection on $E$. Let\/ $\n'$
denote the connection dual to $\n$ with respect to a Hermitian metric $h^{E}$ on $E$. Let\/ $N$ be a compact oriented manifold whose boundary is
diffeomorphic to two disjoint copies of $M$. Then
\eq{dualityRAT}
    \alp\big(\,\rat(\n)\,\big) \ = \ \rat(\n')\cdot e^{ 2i\pi\,\oeta(\n,g^M) - i\pi\,\rank E\,\int_{N}\, L(p,g^M)},
\end{equation}
where $\alp$ is the anti-linear isomorphism \refe{alpE-E'}, $g^M$ is any Riemannian metric on $M$, and $L(p,g^M)$ is the Hirzebruch $L$-polynomial in
the Pontrjagin forms of any Riemannian metric on $N$ which near $M$ is the product of\/ $g^M$ and the standard metric on the half-line.

In particular, if $\dim M\equiv 1 (\MOD 4)$, then
\eq{dualityRATodd}
    \alp\big(\,\rat(\n)\,\big) \ = \ \rat(\n')\cdot e^{2i\pi\,\oeta(\n,g^M)}.
\end{equation}
\eth

The rest of this section is concerned with the proof of \reft{dualityRAT}.

\subsection{Choices of the metric and the spectral cut}\Label{SS:choices}
Till the end of this section we fix a Riemannian metric $g^M$ on $M$ and set $\B= \B(\n,g^M)$ and $\B'= \B(\n',g^M)$. We also fix an Agmon angle
$\tet\in (-\pi/2,0)$ for the odd signature operator $\B= \B(\n,g^M)$ such that there are no eigenvalues of $\B$ in the solid angles
$L_{[-\tet-\pi,-\pi/2]}$, $L_{(-\pi/2,\tet]}$, $L_{[-\tet,\pi/2)}$, and $L_{(\pi/2,\tet+\pi]}$.

Let $\B'$ denote the odd signature operator associated to the connection $\n'$ and the metric $g^M$. One easily checks, cf.
\cite[Prop.~3.58]{BeGeVe}, that
\eq{n*=2}
    \n^* \ = \ \Gam\,\n'\,\Gam \qquad \text{and}\qquad
    (\n')^* \ = \ \Gam\,\n\,\Gam.
\end{equation}
Using \refe{oddsignGam}, \refe{n*=2}, and the equality $\Gam^*= \Gam$ (cf. Proposition~3.58 of \cite{BeGeVe}), one readily sees that the adjoint
$\B^*$ of $\B$ satisfies
\eq{Balp*-Balp'}
    \B^* \ = \ \B'.
\end{equation}
Our choice of the angle $\tet$ guarantees that $\pm2\tet$ are Agmon angles for the operator
\[
    (\Gam\n')^2 \ = \ \big(\,(\Gam\n)^2\,\big)^*.
\]
In particular, for each $\lam\ge0$, the number $\xi_\lam(\n',g^M,\tet)$ can be defined by  formula \refe{xi}, with the same angle $\tet$ and with
$\n$ replaced by $\n'$ everywhere.

\subsection{A choice of $\lam$}\Label{SS:choicelam}
Since the leading symbol of $\B$ is self-adjoint, there are at most finitely many purely imaginary eigenvalues of $\B$. Hence, there exists
$\lam\ge0$ such that there are no purely imaginary eigenvalues of $\B$ with absolute value $\ge\sqrt{\lam}$. In other words, the operator
$\B^{(\lam,\infty)}$ does not have purely imaginary eigenvalues. Moreover, our assumptions on $\tet$ in \refss{choices} imply that no eigenvalue of
$\B^{(\lam,\infty)}$ lies in the solid angles $L_{[-\tet-\pi,\tet]}$ and $L_{[-\tet,\tet+\pi]}$. It follows that no eigenvalue of the operator
$(\B^{(\lam,\infty)})^2$ lies in the solid angle $L_{[-2\tet,2\tet+2\pi]}$.

\lem{na-na'xi}
Let $\tet$ be as in \refss{choices} and let $\lam\ge0$ be big enough so that the operator $\B^{(\lam,\infty)}$ does not have purely imaginary
eigenvalues. Then
\eq{na-na'xi}
        \xi_\lam(\n',g^M,\tet) \ = \ \oxi_\lam(\n,g^M,\tet),
\end{equation}
and
\eq{na-na'eta}
    \eta_\lam(\n',g^M) \ = \ \oeta_\lam(\n,g^M),
\end{equation}
where $\oz$ denotes the complex conjugate of the number $z\in \CC$.
\elem
\prf
Let $\B^{\prime(\lam,\infty)}$ denote the restriction of $\B'$ to the span of the generalized eigenvectors of $(\B')^2$ corresponding to the
eigenvalues whose absolute values are greater than $\lam$. From \refe{Balp*-Balp'}, we obtain
\eq{B'*lam=Blam}
    \big(\B^{\prime(\lam,\infty)}\big)^*\ =\ \B^{(\lam,\infty)}.
\end{equation}
Hence, with our assumptions on $\tet$, we have, for any $j=0\nek d$,
\eq{LDetB-B'}
 \begin{aligned}
    \LD_{2\tet}\Big(\,\big(\,\B^{\prime(\lam,\infty)}\,\big)^2{}_{\big|_{\Ome^j_{(\lam,\infty)}(M,E)}}\,\Big) \ &= \
    \overline{  \LD_{-2\tet}\Big(\,\big(\,\big(\,\B^{\prime(\lam,\infty)}\,\big)^*\,\big)^2{}_{\big|_{\Ome^j_{(\lam,\infty)}(M,E)}}\,\Big)}
    \\ &= \
    \overline{  \LD_{-2\tet}\Big(\,\big(\,\B^{(\lam,\infty)}\,\big)^2{}_{\big|_{\Ome^j_{(\lam,\infty)}(M,E)}}\,\Big)}.
 \end{aligned}
\end{equation}
Since there are no eigenvalues of $\big(\B^{(\lam,\infty)}\big)^2$  in the solid angle $L_{[-2\tet,2\tet+2\pi]}$,
\[
    \LD_{-2\tet}\Big(\,\big(\,\B^{(\lam,\infty)}\,\big)^2{}_{\big|_{\Ome^j_{(\lam,\infty)}(M,E)}}\,\Big)
    \ = \ \LD_{2\tet}\Big(\,\big(\,\B^{(\lam,\infty)}\,\big)^2{}_{\big|_{\Ome^j_{(\lam,\infty)}(M,E)}}\,\Big).
\]
The equality \refe{na-na'xi} follows now from \refe{LDetB-B'} and the definition \refe{xi2} of $\xi_\lam(\n,g^M,\tet)$.

Let $\Pi_>$  (resp. $\Pi_<$) denote the spectral projection of $\B$ onto the span of all generalized eigenvectors of $\B$ corresponding to
eigenvalues with positive (resp. negative) real part. Let $Q$ denote the spectral projection of $\B$ onto the span of all generalized eigenvectors of
$\B$ corresponding to eigenvalues whose absolute value is larger than $\lam$. Let $\Pi'_>$, $\Pi'_<$, and $Q'$ be similarly defined spectral
projections of $\B'$. Then, since the operators $\B^{(\lam,\infty)}$ and $\B^{\prime(\lam,\infty)}$ have no purely imaginary eigenvalues, we conclude
from \refe{etasV} and \refe{etaV} that
\eq{etalamnn'}
 \begin{aligned}
    2\,\eta_\lam(\n,g^M) \ &= \ \zet_\tet(0,Q\Pi_>,\B_\even) \ - \ \zet_\tet(0,Q\Pi_<,\B_\even),\\
    2\,\eta_\lam(\n',g^M) \ &= \ \zet_\tet(0,Q'\Pi'_>,\B'_\even) \ - \ \zet_\tet(0,Q'\Pi'_<,\B'_\even).
 \end{aligned}
\end{equation}
From \refe{B'*lam=Blam} and our assumptions on $\tet$, we obtain
\eq{zetPi><}
\begin{aligned}
    \zet_\tet(s,Q'\Pi'_>,\B'_\even) \ & = \  \ozet_\tet(\os,Q\Pi_>,\B_\even),              \\
    \zet_\tet(s,Q'\Pi'_<,\B'_\even) \ &= \  e^{-2\pi i s}\cdot \ozet_\tet(\os,Q\Pi_<,\B_\even).
\end{aligned}
\end{equation}
The equality \refe{na-na'eta} follows immediately from \refe{etalamnn'} and \refe{zetPi><}.
\eprf

\subsection{Small eigenvalues of $\B$ and $\B'$}\Label{SS:small}
Let $\Ome^{\b}_{[0,\lam]}(M,E')\subset \Ome^\b(M,E)$ denote the span of the eigenvectors of $(\B')^2$ corresponding to the eigenvalues with absolute
value $\le\lam$. Then $\Ome^{\b}_{[0,\lam]}(M,E')$ is a subcomplex of $(\Ome^\b(M,E),\n')$ preserved by the chirality operator $\Gam$.

The pairing \refe{pairingEE'} defines a non-degenerate sesquilinear pairing
\eq{pairingOme'Omelam}
    \Ome^{j}_{[0,\lam]}(M,E')\,\times\,\Ome^{d-j}_{[0,\lam]}(M,E) \ \longrightarrow \ \CC,
\end{equation}
and, hence, identifies $\Ome^{\b}_{[0,\lam]}(M,E')$ with the dual complex of $\Ome^{\b}_{[0,\lam]}(M,E)$.

As in \refs{witheta}, set, for $j=0\nek d$,
\[
    d_{j,\lam} \ = \ \dim \Ome^j_{[0,\lam]}(M,E), \quad d'_{j,\lam} \ = \ \dim \Ome^{j}_{[0,\lam]}(M,E').
\]
Since $\Gam^2=1$, we obtain from \refe{oddsignGam}, $\B_j= \Gam\circ\B_{d-j}\circ\Gam$  and, hence,
\[
    \Gam\,\Big(\,\Ome^j_{[0,\lam]}(M,E)\,\Big)\ =\ \Ome^{d-j}_{[0,\lam]}(M,E), \qquad j=0\nek d.
\]
Therefore
\eq{dj=dN-j}
    d_{j,\lam} \ = \ d_{d-j,\lam}, \qquad j=0\nek d.
\end{equation}
Hence,
\eq{-1jdj}
    \sum_{j=0}^d\,(-1)^jj\,d_{j,\lam}\ = \ \sum_{p=0}^{r-1}\,\big(\,2p-(d-2p)\,\big)\,d_{2p,\lam}
    \ = \ 4\,\sum_{p=0}^{r-1}\,p\,d_{2p,\lam} \ - \ d\,\sum_{p=0}^{r-1}\,d_{2p,\lam}.
\end{equation}
In particular,
\eq{dj=Omeeven}
    \sum_{j=0}^d\,(-1)^jj\,d_{j,\lam}\ \equiv \ \sum_{p=0}^{r-1}\,d_{2p,\lam} \ = \ \dim \Ome^\even_{[0,\lam]}(M,E), \qquad \MOD\ 2\,\ZZ.
\end{equation}

From \refe{Balp*-Balp'}, we conclude that $d_{j,\lam}= d'_{d-j,\lam}$ $(j=0\nek d)$. Combining this equality with \refe{dj=dN-j}, we get $d_{j,\lam}=
d'_{j,\lam}$. Hence, by \refe{dj=Omeeven},
\eq{d'j=Omeeven}
    \sum_{j=0}^d\,(-1)^jj\,d'_{j,\lam}\ \equiv \  \dim \Ome^\even_{[0,\lam]}(M,E), \qquad \MOD\ 2\,\ZZ.
\end{equation}
From \refe{eta=eta}, \refe{eta0lam}, and \refe{d'j=Omeeven}, we obtain, modulo $2\ZZ$,
\eq{eta+d=eta}\notag
 \begin{aligned}
    2\,\eta(\B_\even(\n,g^M)) \ &= \ 2\,\eta(\B_\even^{(\lam,\infty)}(\n,g^M)) \ +\ 2\,\eta(\B^{[0,\lam]}_\even(\n,g^M))
    \\ &\equiv \ 2\,\eta_\lam(\n,g^M) \ +\ \sum_{j=0}^d\,(-1)^jj\,d_{j,\lam}.
 \end{aligned}
\end{equation}
Similarly,
\eq{etaprime+d=eta}
    2\,\eta(\B_\even(\n',g^M)) \ \equiv \ 2\,\eta_\lam(\n',g^M) \ +\ \sum_{j=0}^d\,(-1)^jj\,d_{j,\lam}, \qquad \MOD\ 2\,\ZZ.
\end{equation}

\subsection{Proof of \reft{dualityRAT}}\Label{SS:prdualyRAT}
Let $\rho'_{{}_{\Gam_{[0,\lam]}}}$ denote the refined torsion of the complex $\Ome^{\b}_{[0,\lam]}(M,E')$ associated to the restriction of $\Gam$ to
$\Ome^{\b}_{[0,\lam]}(M,E')$.

Since $\Gam^*=\Gam$ (cf. Proposition~3.58 of \cite{BeGeVe}), we obtain from \refl{reftordual} and the definition \refe{alpE-E'} of $\alp$,
\eq{rholam'=rholam}
    \rho'_{{}_{\Gam_{[0,\lam]}}} \ = \ \alp(\rho_{{}_{\Gam_{[0,\lam]}}}).
\end{equation}

From \refe{rho}, \refe{xietad}, and \refd{refantor}, we obtain
\meq{rhoan=rhoxi...}
    \rat(\n) \ = \ \rho_{{}_{\Gam_{[0,\lam]}}}\cdot \exp\Big(\,\xi_\lam(\n,g^M,\tet)-i\pi\eta_\lam(\n,g^M)
    \\-\frac{i\pi}2\,\sum_{j=0}^d\,(-1)^jj\,d_{j,\lam}+i\pi\,\frac{\rank E}2\,\int_N\,L(p,g^M)\,\Big).
\end{multline}
Since $\alp$ is an anti-linear isomorphism, $\alp(\rat\cdot{z})= \alp(\rat)\cdot\oz$ for any $z\in \CC$. Hence, from \refe{rholam'=rholam}
 and \refe{rhoan=rhoxi...} we get
\meq{alp(rat)}
    \alp\big(\rat(\n)\big) \ = \ \rho'_{{}_{\Gam_{[0,\lam]}}}\cdot \exp\Big(\,\oxi_\lam(\n,g^M,\tet)+i\pi\oeta_\lam(\n,g^M)
    \\ +\frac{i\pi}2\,\sum_{j=0}^d\,(-1)^jj\,d_{j,\lam}-i\pi\,\frac{\rank E}2\,\int_N\,L(p,g^M)\,\Big).
\end{multline}
Using \refl{na-na'xi} and the analogue of \refe{rhoan=rhoxi...} for $\rat(\n')$, we obtain from \refe{alp(rat)}
\meq{alp(rat)2}
    \alp\big(\rat(\n)\big) \ = \ \rat(\n')\cdot \exp\Big(\,2i\pi\oeta_\lam(\n,g^M)
     \\ + {i\pi}\,\sum_{j=0}^d\,(-1)^jj\,d_{j,\lam}-i\pi\,\rank E\,\int_N\,L(p,g^M)\,\Big).
\end{multline}
From \refe{alp(rat)2} and \refe{etaprime+d=eta} we obtain \refe{dualityRAT}.\hfill$\square$

\section{Comparison with the Ray-Singer Torsion}\Label{S:Ray-Singer}

In this section we calculate the Ray-Singer norm $\RSn{\rat}$ of the refined analytic torsion. In particular, we show that, if $\n$ is a Hermitian
connection, then $\RSn{\rat}= 1$.

\subsection{The Ray-Singer torsion}\Label{SS:RStorsion}
Let $E\to M$ be a complex vector bundle over a closed oriented manifold $M$ of odd dimension $d=2r-1$ and let $\n$ be a flat connection on $E$. Fix a
Riemannian metric $g^M$ on $M$ and a Hermitian metric $h^{E}$ on $E$. Let $\n^*$ denote the adjoint of\/ $\n$ with respect to the scalar product
$\<\cdot,\cdot\>$ on $\Ome^\b(M,E)$ defined by $h^{E}$ and the Riemannian metric $g^M$. Let
\eq{Laplacian}
    \Del\ =\ \n^*\,\n\ +\ \n\,\n^*
\end{equation}
be the Laplacian. We denote by $\Del_k$ the restriction of $\Del$ to $\Ome^k(M,E)$.

The {\em Ray-Singer torsion} $\TRS$ of $E$, \cite{RaSi1,BisZh92,BFK3}, is defined by
\footnote{Our sign convention is different from \cite{RaSi1}, \cite{BFK3}, and \cite{BrKappelerRAT} but is consistent with \cite{BisZh92}. In
our notations, the torsion defined in \cite{RaSi1,BFK3,BrKappelerRAT} is equal to $1/\TRS$.}
\eq{RaySingertor}
    \TRS \ = \TRS(\n) := \ \exp\,\Big(\,
     \frac12\,\sum_{k=0}^d\,(-1)^{k}\,k\, \LD_{-\pi}(\Del_k)\,\Big),
\end{equation}
Note that $\Del_k$ is a self-adjoint non-negative operator. Therefore, all its eigenvalues are non-negative and $\LD_{-\pi}(\Del_k)$ is well defined.

More generally, suppose $\calI$\/ is an interval of the form $[0,\lam],\ (\lam,\mu]$, or $(\lam,\infty]$ ($\mu\ge\lam\ge0$) and let $\Pi_{\Del_k,\calI}$ be the
spectral projection of $\Del$ corresponding to $\calI$, cf. \refss{spectralsubspace}.  Set
\eq{hatOmecalI}\notag
    \hatOme^k_{\calI}(M,E)\ := \ \Pi_{\Del_k,\calI}\big(\, \Ome^\b(M,E)\,\big)\ \subset\  \hatOme^\b(M,E).
\end{equation}
Let $\Del^\calI_k$ denote the restriction of $\Del_k$ to $\hatOme^k_\calI(M,E)$ and define
\eq{RaySingertorI}
    \TRS_\calI \ = \TRS_\calI(\n) := \ \exp\,\Big(\,
     \frac12\,\sum_{k=0}^d\,(-1)^{k}\,k\, \LD_{-\pi}(\Del^\calI_k)\,\big)\,\Big).
\end{equation}
In particular, by \refe{RaySingertor},  $\TRS= \TRS_{(0,\infty)}$. By \refe{Det=DetlamDetmu}, for any non-negative real numbers $\mu\ge \lam\ge 0$,
\eq{TRSIJ}
    \TRS_{(\lam,\infty)} \ = \ \TRS_{(\lam,\mu]}\cdot \TRS_{(\mu,\infty]}.
\end{equation}


\subsection{The Ray-Singer metric on the determinant line of cohomology}\Label{SS:RSnorm}
If the connection $\n$ is acyclic, i.e., if the cohomology $H^\b(M,E)$ vanishes, then the Ray-Singer torsion is independent of the Hermitian metric
$h^{E}$ and the Riemannian metric $g^M$, cf. \cite{RaSi1,BisZh92}. If the cohomology does not vanish, then $\TRS$, in general, depends on the choice
of the metrics. To construct a metric independent invariant of the flat vector bundle $E$ one needs to take into account the contribution of the
space of harmonic forms. An elegant way to do this, which was proposed by Quillen \cite{Quillen85}, is to construct a norm $\RSn{\cdot}$ on the
determinant line of $H^\b(M,E)$, called the {\em Ray-Singer metric}, which is independent of $g^M$ and $h^E$. We now briefly recall this
construction.

{}For each $\lam\ge 0$, the cohomology of the finite dimensional complex $\big(\hatOme^\b_{[0,\lam]}(M,E),\n\big)$ is naturally isomorphic to
$H^\b(M,E)$. Identifying these two cohomology spaces, we then obtain from \refe{isomorphism2} an isomorphism
\eq{philam}
    \phi_\lam \ = \ \phi_{{}_{\hatOme^\b_{[0,\lam]}(M,E)}}:\, \Det\big(\,\hatOme^\b_{[0,\lam]}(M,E)\,\big) \ \longrightarrow \
    \Det\big(\,H^\b(M,E)\,\big).
\end{equation}

The scalar product $\<\cdot,\cdot\>$ on $\hatOme^\b_{[0,\lam]}(M,E)\subset \Ome^\b(M,E)$ defined by $g^M$ and $h^E$ induces a metric on the
determinant line $\Det\big(\hatOme^\b_{[0,\lam]}(M,E)\big)$.  Let $\|\cdot\|_\lam$ denote the metric on $\Det\big(H^\b(M,E)\big)$ such that the
isomorphism \refe{philam} is an isometry. It is well known, cf., for example, Theorem~1.1 of \cite{BisZh92}, that for $0\le\lam\le\mu$
\eq{normlammu}
    \|\cdot\|_\mu \ = \ \|\cdot\|_\lam\cdot \TRS_{(\lam,\mu]}.
\end{equation}

The {\em Ray-Singer metric} on $\Det\big(H^\b(M,E)\big)$ is defined by the formula
\eq{RSnorm}
    \RSn{\cdot} \ := \  \|\cdot\|_\lam\cdot \TRS_{(\lam,\infty)}, \qquad \lam\ge 0.
\end{equation}
It follows immediately from \refe{TRSIJ} and \refe{normlammu} that $\RSn{\cdot}$ is independent of the choice of $\lam\ge0$.

\th{RSnormRAT}
Let $E$ be a complex vector bundle over a closed oriented odd-dimensional manifold $M$ and let $\n$ be a flat connection on $E$. Then
\eq{RSnormRAT}
    \RSn{\rat} \ = \ e^{\pi\,\IM \eta(\n,g^M)},
\end{equation}
where
\[
    \eta(\n,g^M)\ = \ \eta\big(\B_\even(\n,g^M)\big).
\]
In particular, if\/ $\n$ is a Hermitian connection, then $\eta(\n,g^M)\in \RR$ and
\eq{RSnormRATherm}
    \RSn{\rat} \ = \ 1.
\end{equation}
\eth
The rest of this section is concerned with the proof of \reft{RSnormRAT}.

\subsection{Choice of the spectral cut}\Label{SS:choicetet}
The determinants in \refe{RaySingertor} and \refe{RaySingertorI} are defined using the spectral cut $R_{-\pi}$ along the negative real axis. Since
the spectrum of the operator $\Del_k$ lies on the positive real axis, we can replace $R_{-\pi}$ with a spectral cut $R_\gam$ for any nonzero
$-\pi\le\gam<\pi$ without changing the value of $\Det'_\tet(\Del_k)$. In particular, we can take the spectral cut along $R_{2\tet}$, where $\tet\in
(-\pi/2,0)$ is any Agmon angle for the odd signature operator $\B= \B(\n,g^M)$ such that there are no eigenvalues of $\B$ in the solid angles
$L_{[-\tet-\pi,-\pi/2]}$, $L_{(-\pi/2,\tet]}$, $L_{[-\tet,\pi/2)}$, and $L_{(\pi/2,\tet+\pi]}$. We fix such an angle $\tet$ till the end of this
section.

\subsection{The Ray-Singer metric and the dual connection}\Label{SS:RSnormdual}
Let $\n'$ be the connection dual to $\n$ with respect to the Hermitian metric $h^E$, cf. \refss{n'}, and let $E'$ denote the flat bundle $(E,\n')$.
Let
\eq{Laplacianprime}\notag
    \Del'\ = \ (\n')^*\,\n' \ + \ \n'\,(\n')^*,
\end{equation}
denote the Laplacian of the connection $\n'$. For $\lam\ge0$, we denote by $\hatOme^{\b}_{[0,\lam]}(M,E')\subset \Ome^\b(M,E)$ the image of the
spectral projection $\Pi_{\Del',[0,\lam]}$, cf. \refss{spectralsubspace}. As in \refss{RSnorm}, we use the scalar product induced by $g^M$ and $h^E$
on $\hatOme^{\b}_{[0,\lam]}(M,E')$ to construct a metric $\|\cdot\|'_\lam$ on $\Det(H^\b(M,E'))$ and we define the Ray-Singer metric on
$\Det(H^\b(M,E'))$ by the formula
\eq{RSnormprime}
     \RSnp{\cdot} \ := \  \|\cdot\|'_\lam\cdot \TRS_{(\lam,\infty)}(\n').
\end{equation}

\subsection{Comparison between the Ray-Singer metrics associated to a connection and to its dual}\Label{SS:RS-RSprime}
From \refe{n*=2} we conclude that
\eq{Laplace-Laplacianprime}\notag
    \Del'\ = \ \Gam\circ\Del\circ\Gam.
\end{equation}
Hence, for each $\lam\ge0$, $j=0\nek d$,
\eq{GamhatOme=hatOme}
    \Gam\big(\,\hatOme^{j}_{[0,\lam]}(M,E')\,\big) \ = \ \hatOme^{d-j}_{[0,\lam]}(M,E).
\end{equation}
Recall that the notation $h^E(\alp\wedge\bet)$ was introduced in \refe{OmeEOmeE'}. It follows from \refe{GamhatOme=hatOme} and \refe{OmeEOmeE'} that
the map
\eq{Omelam-Omelam'}
    (\ome',\ome)\ \mapsto \ \int_M\,h^E(\Gam\ome'\wedge\ome), \qquad \ome\in \hatOme^{\b}_{[0,\lam]}(M,E), \ \ome'\in \hatOme^{\b}_{[0,\lam]}(M,E')
\end{equation}
defines a non-degenerate sesquilinear pairing between $\hatOme^{\b}_{[0,\lam]}(M,E')$ and $\hatOme^{\b}_{[0,\lam]}(M,E)$ and, hence, identifies
$\hatOme^{\b}_{[0,\lam]}(M,E')$ with the dual space of $\hatOme^{\b}_{[0,\lam]}(M,E)$. Moreover, this identification preserves the scalar products
induced by $g^M$ and $h^E$ on $\hatOme^{\b}_{[0,\lam]}(M,E')$ and on the dual to $\hatOme^{\b}_{[0,\lam]}(M,E)$. Hence, the anti-linear isomorphism
\refe{alpE-E'} is an isometry with respect to the metrics $\|\cdot\|_\lam$ and $\|\cdot\|'_\lam$. In particular,
\[
    \|\rat(\n)\|_\lam \ = \ \|\alp(\rat(\n))\|'_\lam.
\]
It follows now from \refe{dualityRAT} that
\eq{normrat-normprime}
    \|\rat(\n)\|_\lam \ = \
    \|\rat(\n')\|'_\lam\cdot e^{2\pi \IM\eta(\n,g^M)}.
\end{equation}

A verbatim repetition of the proof of Lemma~8.8 of \cite{BrKappelerRAT} yields that for each $\lam\ge0$
\eq{na-na'TRS}
        \TRS_{(\lam,\infty)}(\n') \ =\ \TRS_{(\lam,\infty)}(\n).
\end{equation}
Hence, from \refe{RSnorm}, \refe{RSnormprime}, and \refe{normrat-normprime}, we conclude that
\eq{RSnormrat-normprime}
    \RSn{\rat(\n)} \ = \ \RSnp{\rat(\n')}\cdot e^{2\pi \IM\eta(\n,g^M)}.
\end{equation}

\subsection{Direct sum of a connection and its dual}\Label{SS:n+n'}
Let
\[
   \tiln \ = \ \begin{pmatrix}
                        \n&0\\
                        0&\n'
                       \end{pmatrix},
\]
denote the flat connection on $E\oplus{E}$ obtained as a direct sum of the connections $\n$ and $\n'$. Clearly, for each $\lam\ge0$,
\eq{xietatil}
  \begin{aligned}
     \xi_\lam(\tiln,g^M,\tet) \ &= \ \xi_\lam(\n,g^M,\tet)\ + \ \xi_\lam(\n',g^M,\tet), \\
     \eta_\lam(\tiln,g^M) \ &= \ \eta_\lam(\n,g^M)\ + \ \eta_\lam(\n',g^M).
  \end{aligned}
\end{equation}

From \refl{rhodirectsum} we obtain for $\lam\ge0$
\eq{na-na'rho}
    \rho_{{}_{\Gam_{{}_{[0,\lam]}}}}(\tiln,g^M)  \ = \
    \mu_{H^\b(M,E),H^\b(M,E')}\big(\,\rho_{{}_{\Gam_{{}_{[0,\lam]}}}}(\n,g^M)\otimes\rho_{{}_{\Gam_{{}_{[0,\lam]}}}}(\n',g^M)\,\big).
\end{equation}

Hence, in view of \refe{rho}, \refe{xietad}, and \refd{refantor}, we obtain
\eq{na-na'rat}
    \rat(\tiln)  \ = \
    \mu_{H^\b(M,E),H^\b(M,E')}\big(\,\rat(\n)\otimes\rat(\n')\,\big).
\end{equation}

Let $\tilDel=\tiln^*\tiln+\tiln\tiln^*$ be the Laplacian of the connection $\tiln$ and let $\hatOme^\b_{[0,\lam]}(M,E\oplus E')$ denote the span of
the eigenvectors of $\tilDel$ corresponding to eigenvalues which are $\le\lam$. As in \refss{RSnorm}, we use the scalar product induced by $g^M$ and
$h^E$ on $\hatOme^{\b}_{[0,\lam]}(M,E\oplus E')$ to construct a metric $\|\cdot\|^{{}^\sim}_\lam$ on $\Det(H^\b(M,E\oplus{}E'))$ and we define the
Ray-Singer metric on $\Det(H^\b(M,E\oplus{}E'))$ by the formula
\eq{RSnormnormprime}
     \RSnnp{\cdot} \ := \  \|\cdot\|^{{}^\sim}_\lam\cdot \TRS_{(\lam,\infty)}(\tiln).
\end{equation}
Since
\[
    \hatOme^\b_{[0,\lam]}(M,E\oplus E')\ = \ \hatOme^\b_{[0,\lam]}(M,E)\,\oplus\,\hatOme^\b_{[0,\lam]}(M,E'),
\]
it follows from the definition \refe{directsum2} of the fusion isomorphism that, for any $h\in H^\b(M,E)$ and $h'\in H^\b(M,E')$,
\eq{normmu=product}
    \|\mu_{H^\b(M,E),H^\b(M,E')}(h\otimes h')\|^{{}^\sim}_\lam \ = \ \|h\|_\lam\cdot\|h'\|'_\lam.
\end{equation}
Therefore, we obtain from \refe{na-na'rat}
\eq{normlamrattil}
    \|\rat(\tiln)\|^{{}^\sim}_\lam \ = \ \|\rat(\n)\|_\lam\cdot\|\rat(\n')\|'_\lam.
\end{equation}
Since
\[
    \TRS_{(\lam,\infty)}(\tiln) \ = \ \TRS_{(\lam,\infty)}(\n)\cdot\TRS_{(\lam,\infty)}(\n'),
\]
we conclude from \refe{RSnormnormprime} and \refe{normlamrattil} that
\eq{normrattil0}
    \RSnnp{\rat(\tiln)} \ = \ \RSn{\rat(\n)}\cdot\RSnp{\rat(\n')}.
\end{equation}
Combining the later equality with \refe{RSnormrat-normprime}, we get
\eq{normrattil}
    \RSnnp{\rat(\tiln)} \ = \ \big(\,\RSn{\rat(\n)}\,\big)^2\cdot e^{-2\pi \IM\eta(\n,g^M)}.
\end{equation}
Hence, {\em \refe{RSnormRAT} is equivalent to the equality}
\eq{nomrattil=1}
    \RSnnp{\rat(\tiln)} \ = \ 1.
\end{equation}

\subsection{Deformation of the chirality operator}\Label{SS:deformationGam}
We will prove \refe{nomrattil=1} by a deformation argument. For $t\in [-\pi/2,\pi/2]$ introduce the rotation $U_t$ on
\[
    \Ome^\b \ := \ \Ome^\b(M,E)\,\oplus\,\Ome^\b(M,E),
\]
given by
\[
    U_t \ = \ \begin{pmatrix}
    \ \cos t&-\,\sin t \ \\
    \ \sin t& \ \ \ \cos t \
    \end{pmatrix}.
\]
Note that $U_t^{-1}= U_{-t}$. Denote by $\tilGam(t)$ the deformation of the chirality operator, defined by
\eq{Gam(t)}
        \tilGam(t) \ = \ U_t\circ
                    \begin{pmatrix}
                        \Gam&0\\ 0&-\Gam
                    \end{pmatrix}
                    \circ U_t^{-1} \ = \ \Gam\circ
                    \begin{pmatrix}
                        \cos 2t&\sin 2t\\ \sin 2t&-\cos2t
                    \end{pmatrix}.
\end{equation}
Then
\eq{Gam0pi4}
    \tilGam(0) \ = \ \begin{pmatrix}
                        \Gam&0\\ 0&-\Gam
                    \end{pmatrix}, \qquad
    \tilGam(\pi/4) \ = \ \begin{pmatrix}
                        0&\Gam\\ \Gam&0
                    \end{pmatrix}.
\end{equation}

\subsection{Deformation of the odd signature operator}\Label{SS:deformationB}
Consider a one-parameter family of operators $\tB(t):\Ome^\b\to \Ome^\b$ ($t\in [-\pi/2,\pi/2]$) defined by the formula
\eq{tB(t)}
    \tB(t) \ := \ \tilGam(t)\,\tiln\ + \ \tiln\,\tilGam(t).
\end{equation}
Then
\eq{tilDel0}
    \tB(0) \ = \ \begin{pmatrix}
                        \B&0\\ 0&-\B'
                    \end{pmatrix}
\end{equation}
and
\eq{tilBpi/4}
    \tB(\pi/4) \ = \ \begin{pmatrix}
                        0&\Gam\n'+\n\Gam\\ \Gam\n+\n'\Gam&0
                    \end{pmatrix}.
\end{equation}
Hence, using \refe{n*=2}, we obtain
\eq{tilDelpi}
                    \tB(\pi/4)^2 \ = \ \begin{pmatrix}
                        \Del&0\\ 0&\Del'
                    \end{pmatrix} \ = \ \tilDel.
\end{equation}
Set
\begin{gather}
    \Ome^\b_+(t) \ := \ \Ker\, \tiln\,\tilGam(t);\notag\\
        \Ome^\b_- \ := \ \Ker\tiln \ = \ \Ker \n\oplus \Ker\n'.\notag
\end{gather}
Note that $\Ome^\b_-$ is independent of $t$. Since the operators $\tiln$ and $\tilGam(t)$ commute with $\tB(t)$, the spaces $\Ome^\b_+(t)$ and
$\Ome^\b_-$ are invariant for $\tB(t)$.

Let $\calI$ be an interval of the form $[0,\lam],\ (\lam,\mu]$, or $(\lam,\infty]$ ($\mu\ge\lam\ge0$). Denote
\eq{Omecal(t)I}\notag
    \Ome^\b_{\calI}(t)\ := \ \Pi_{\tB(t)^2,\calI}\big(\, \Ome^\b(t)\,\big)\ \subset\  \Ome^\b(t),
\end{equation}
where $\Pi_{\tB(t)^2,\calI}$ is the spectral projection of $\tB(t)^2$ corresponding to $\calI$, cf. \refss{spectralsubspace}. For $j=0\nek d$, set
$\Ome^j_\calI(t)= \Ome^\b_\calI(t)\cap\Ome^j$ and
\eq{Ome+-tlam}
    \Ome^j_{\pm,\calI}(t) \ := \ \Ome^j_\pm(t)\cap \Ome^j_\calI(t).
\end{equation}

For each $\lam\ge0$, $t\in (-\pi/2,\pi/2)$, the space $\Ome^\b_{(\lam,\infty)}(t)$ is invariant by  $\tB^{(\lam,\infty)}(t)$ and the operator
$\tB^{(\lam,\infty)}(t):\Ome^\b_{(\lam,\infty)}(t)\to \Ome^\b_{(\lam,\infty)}(t)$ is bijective. Since the range of the restriction of\/
$\tilGam(t)\tiln$ to $\Ome^\b_{(\lam,\infty)}(t)$ is contained in $\Ome^\b_{+,(\lam,\infty)}(t)$ whereas the range of the restriction of\/
$\tiln\tilGam(t)$ to $\Ome^\b_{(\lam,\infty)}(t)$ is contained in $\Ome^\b_{-,(\lam,\infty)}(t)$, it follows from the surjectivity of\/
$\tB^{(\lam,\infty)}(t)$ that
\eq{Omet+-1}
    \Ome^\b_{+,(\lam,\infty)}(t)\, +\, \Ome^\b_{-,(\lam,\infty)}(t) \ = \ \Ome^\b_{(\lam,\infty)}(t),\qquad t\in [-\pi/2,\pi/2].
\end{equation}
Similarly, since the kernel of the restriction of\/ $\tiln\tilGam(t)$ to $\Ome^\b_{(\lam,\infty)}(t)$ is equal to $\Ome^\b_{+,(\lam,\infty)}(t)$
whereas the kernel of the restriction of\/ $\tilGam(t)\tiln$ to $\Ome^\b_{(\lam,\infty)}(t)$ is equal to $\Ome^\b_{-,(\lam,\infty)}(t)$, it follows
from the injectivity of\/ $\tB^{(\lam,\infty)}(t)$ that
\eq{Omet+-2}
    \Ome^\b_{+,(\lam,\infty)}(t)\cap \Ome^\b_{-,(\lam,\infty)}(t)  \ = \ \{0\},\qquad t\in [-\pi/2,\pi/2].
\end{equation}

Combining \refe{Omet+-1} and \refe{Omet+-2} we obtain
\eq{Omet+-}
    \Ome^\b_{(\lam,\infty)}(t) \ = \ \Ome^\b_{+,(\lam,\infty)}(t)\, \oplus\, \Ome^\b_{-,(\lam,\infty)}(t),\qquad t\in [-\pi/2,\pi/2].
\end{equation}

We define $\tB_j^\calI(t),\ \tB_\even^\calI(t),\ \tB_\odd^\calI(t),\ \tB_j^{\pm,\calI}(t),\ \tB_\even^{\pm,\calI}(t),\ \tB_\odd^{\pm,\calI}(t)$, etc.
in the same way as the corresponding maps were defined in \refss{drdetoddsign}.

\subsection{The graded determinant of the deformed odd signature operator}\Label{SS:grdettilB(t)}
By Definitions~\ref{D:detgrfd} and \ref{D:grdeterminantV}, for every $\calI$ of the form $[0,\lam],\ (\lam,\mu]$, or $(\lam,\infty]$
($\mu\ge\lam\ge0$), the graded determinant of $\tB^{\calI}_{\even}(t)$ is given by the formula
\eq{grdettilB(t)}
    \Detgrtet\big(\tB^{\calI}_{\even}(t)\big) \ := \ e^{\LDgrtet(\tB^{\calI}_{\even}(t))},
\end{equation}
where $\tet$ is an Agmon angle for $\tB^{\calI}_{\even}(t)$ and
\eq{grldettilB(t)}
    \LDgrtet\big(\tB^{\calI}_{\even}(t)\big) \ := \ \LD_\tet\big(\tB_{\even}^{+,{\calI}}(t)\big)
      \ - \ \LD_\tet\big(-\tB_{\even}^{-,{\calI}}(t)\big)
    \ \in \ \CC.
\end{equation}
Since $\tilGam(t)$ commutes with $\tB(t)$, we easily obtain
\eq{tilB-tilB+}\notag
    \tB^{+,\calI}_j(t) \ = \ \tilGam(t)\circ \tB^{-,\calI}_{d-j}(t)\circ \tilGam(t), \qquad j=0\nek d.
\end{equation}
Therefore, \refe{grldettilB(t)} can be rewritten as
\eq{grldettilB(t)2}
    \LDgrtet\big(\tB^{\calI}_{\even}(t)\big) \ := \ \sum_{j=0}^d\,(-1)^j\, \LD_\tet\big(\,(-1)^j\,\tB_{j}^{+,{\calI}}(t)\,\big).
\end{equation}

\lem{tB(pi/4)=TRS}
Suppose that $\tet\in (-\pi/2,0)$ is an Agmon angle for the operator $\tB^{(\lam,\infty)}_{\even}(\pi/4)$. Then, for each $\lam\ge0$,
\eq{tB(pi/4)=TRS}
    \big|\,\Detgrtetnp\big(\tB^{(\lam,\infty)}_{\even}(\pi/4)\big)\,\big| \ = \ \frac1{\TRS_{(\lam,\infty)}(\tiln)}.
\end{equation}
\elem
\prf
It follows from \refe{tilBpi/4} that the operator $\tB^{(\lam,\infty)}_{\even}(\pi/4)$ is self-adjoint. Hence,
$\eta\big(0,\tB^{(\lam,\infty)}_{\even}(\pi/4)\big)$ and $\zet_{2\tet}(0,\tB^{(\lam,\infty)}_{\even}(\pi/4)^2\big)$ are real, cf., for example,
Theorem~A.2 of \cite{BrKappelerRAT}. Thus, from \refe{tilDelpi}, \refe{grldettilB(t)2}, and \refe{grdet-etaV}, we conclude that
\eq{tB(pi/4)=xi}
 \begin{aligned}
    \RE\,\LDgrtetnp\big(\tB^{(\lam,\infty)}_{\even}(\pi/4)\big) \ &= \ \frac12\,\sum_{j=0}^{d}\,(-1)^{j}\,
    \LDnp_{2\tet}\big[\,\tilDel^{(\lam,\infty)}{}_{\big|_{\Ome^j_{+,{(\lam,\infty)}}(\pi/4)}}\,\big]
    \\ &= \ \frac12\,\sum_{j=0}^{d}\,(-1)^{j}\, \LDnp_{-\pi}\big[\,\tilDel^{(\lam,\infty)}_j{}_{\big|_{\Ome^j_{+,{(\lam,\infty)}}(\pi/4)}}\,\big].
 \end{aligned}
\end{equation}
As on page~340 of \cite{BFK3} (see also section~8.4 of \cite{BrKappelerRAT}), one shows that the right hand side of \refe{tB(pi/4)=xi} is equal to
\eq{tB(pi/4)=xi2}\notag
    \frac12\,\sum_{j=0}^{d}\,(-1)^{j+1}\,j\,
    \LDnp_{-\pi}\big[\,\tilDel^{(\lam,\infty)}_j\,\big].
\end{equation}
Hence, equality \refe{tB(pi/4)=TRS} follows from\refe{RaySingertorI} and \refe{grdettilB(t)}.
\eprf

\subsection{Deformation of the canonical element of the determinant line}\Label{SS:deformationrho}
Since the operators $\tiln$ and $\tB(t)^2$ commute, the space $\Ome^\b_{\calI}(t)$ is invariant under $\tiln$, i.e., it is a subcomplex of $\Ome^\b$.
The same arguments as in the proof of \refl{tcdoth} show that, for every $\lam\ge0$, the complex $\Ome^\b_{(\lam,\infty)}(t)$ is acyclic and, hence,
the cohomology of the finite dimensional complex $\Ome^\b_{[0,\lam]}(t)$ is naturally isomorphic to
\[
    H^\b(M,E\oplus{}E')\ \simeq\ H^\b(M,E)\oplus{}H^\b(M,E').
\]
Let $\tilGam_{[0,\lam]}(t)$ denote the restriction of $\tilGam(t)$ to $\Ome^\b_{[0,\lam]}(t)$. As $\tilGam(t)$ and $\tB(t)^2$ commute, it follows
that $\tilGam_{[0,\lam]}(t)$ maps $\Ome^\b_{[0,\lam]}(t)$ onto itself and, therefore, is a chirality operator for $\Ome^\b_{[0,\lam]}(t)$. Let
\eq{rhotilGam}
    \rho_{{}_{\tilGam_{\hskip-1pt[0,\lam]}(t)}}(t) \ \in \ \Det\big(H^\b(M,E\oplus E')\,\big)
\end{equation}
denote the refined torsion of the finite dimensional complex $\big(\Ome^\b_{[0,\lam]}(t),\tiln\big)$ corresponding to the chirality operator
$\tilGam_{[0,\lam]}(t)$, cf. \refd{refinedtorsion}.

{}For each $t\in (-\pi/2,\pi/2)$ fix an Agmon angle $\tet= \tet(t)\in (-\pi/2,0)$ for $\tB_{\even}(t)$  and define the element $\rho(t)\in
\Det\big(\,H^\b(M,E\oplus E')\,\big)$ by the formula
\eq{tilrho(t)}
    \rho(t) \ := \ \Detgrtetnp\big(\tB^{(\lam,\infty)}_{\even}(t)\big)\cdot \rho_{{}_{\tilGam_{\hskip-1pt{}_{[0,\lam]}}(t)}}(t),
\end{equation}
where $\lam$ is any non-negative real number. It follows from \refp{antor-grdetfd2} that $\rho(t)$ is independent of the choice of $\lam\ge0$.

\subsection{The Ray-Singer norm of $\rho(t)$}\Label{SS:normrho(t)}
For $t\in [-\pi/2,\pi/2]$, $\lam\ge0$, set
\begin{gather}
    \xi_\lam(t,\tet) \ := \ \frac12\,\sum_{j=0}^{d}\,(-1)^{j+1}\,j\,
    \LDnp_{2\tet}\big[\,{\tB^{(\lam,\infty)}_{\even}(t)^2}_{\big|_{\Ome^j_{{(\lam,\infty)}}(t)}}\,\big], \Label{E:xilam(t)}\\
    \eta_\lam(t) \ := \ \eta\big(\tB^{(\lam,\infty)}_{\even}(t)\big).\Label{E:etalam(t)}
\end{gather}
From \refe{zet-zet}, we see that $\RE\xi_\lam(t,\tet)$ is independent of the choice of the angle $\tet\in (-\pi/2,0)$ such that both $\tet$ and
$\tet+\pi$ are Agmon angles for $\tB(t)$. Hence, for any such angle $\tet\in (-\pi/2,0)$, we obtain from \refe{tilrho(t)} and \refe{xietad},
\eq{normrho(t)=eta}
    \RSnnp{\rho(t)} \ = \ \bigRSnnp{\, \rho_{{}_{\tilGam_{\hskip-1pt{}_{[0,\lam]}}(t)}}(t)\cdot e^{\xi_\lam(t,\tet)}\,}\cdot
    e^{\pi\,\IM\eta_\lam(t)}.
\end{equation}
Since the rank of the operator $\tB^{[0,\lam]}_{\even}(t)$ is finite, $\eta\big(\tB^{[0,\lam]}_{\even}(t)\big)\in \frac12\ZZ$. Hence,
\[
    \IM \eta_\lam(t) \ = \ \IM\big(\, \eta\big(\tB_{\even}(t)\big) \ - \ \eta\big(\tB^{[0,\lam]}_{\even}(t)\big)\,\big) \ = \
     \IM \eta\big(\tB_{\even}(t)\big)
\]
is independent of $\lam\ge0$. We conclude now from \refe{normrho(t)=eta} that
\[
    \bigRSnnp{\, \rho_{{}_{\tilGam_{\hskip-1pt{}_{[0,\lam]}}(t)}}(t)\cdot e^{\xi_\lam(t,\tet)}\,}
\]
is independent of $\lam\ge0$.

\subsection{The Ray-Singer norm of $\rho(0)$}\Label{SS:rho(0)}
By \refe{tilDel0},
\eq{tileta=eta-eta}\notag
    \eta_\lam(0)\ = \ \eta\big(\B^{(\lam,\infty)}_\even\big)\ +\ \eta\big(-\B^{\prime(\lam,\infty)}_\even\big).
\end{equation}
Since the operator $\B^{\prime(\lam,\infty)}_\even$ is invertible $m_0(\B^{\prime(\lam,\infty)}_\even)=0$ (see \refss{etainv} for the definition of
$m_0$). It then follows from the definition \refe{etaV} of the $\eta$-invariant that
\eq{eta(-B)=-eta(B)}
   \eta\big(-\B^{\prime(\lam,\infty)}_\even\big)\ =\ -\eta\big(\B^{\prime(\lam,\infty)}_\even\big).
\end{equation}
Hence, using definition \refe{eta=eta} of $\eta_\lam(\n,g^M)$ and  formula \refe{na-na'eta} for $\eta_\lam(\n',g^M)$, we get
\eq{tileta=IMeta}
    \eta_\lam(0) \ = \ \eta_\lam(\n,g^M) \ - \ \eta_\lam(\n',g^M) \ = \ 2\,i\,\IM \eta_\lam(\n,g^M).
\end{equation}

Using again the invertibility of $\B^{(\lam,\infty)}_\even$, we obtain from \refe{grdet-etaV} and \refe{eta(-B)=-eta(B)} (and \refe{eta=eta} and
\refe{na-na'eta}), that
\eq{Detgr(-B)}\notag
 \begin{aligned}
    \Detgrtetnp\big(-\B^{\prime(\lam,\infty)}_\even\,\big) \ &= \ \Detgrtetnp\big(\,\B^{\prime(\lam,\infty)}_\even\,\big)\cdot
            e^{2\pi i\,\eta(\B^{\prime(\lam,\infty)}_\even)}
    \\ &= \ \Detgrtetnp\big(\,\B^{\prime(\lam,\infty)}_\even\,\big)\cdot
            e^{2\pi i\,\oeta(\n,g^M)}.
 \end{aligned}
\end{equation}
Hence, from  \refe{tilDel0}, we get
\eq{DettilB=BB'}
    \Detgrtetnp\big(\tB^{(\lam,\infty)}_\even(0)\big) \ = \
    \Detgrtetnp\big(\B^{(\lam,\infty)}_\even\big)\cdot\Detgrtetnp\big(\B^{\prime(\lam,\infty)}_\even\big)\cdot
    e^{2\pi i\,\oeta_\lam(\n,g^M)}.
\end{equation}

{}From \refe{tilDel0} we conclude that $\Ome^\b_{[0,\lam]}(0)= \Ome^\b_{[0,\lam]}(M,E)\oplus\Ome^\b_{[0,\lam]}(M,E')$. It follows now from
\refe{Gam0pi4} and \refe{rhodirectsum} that
\eq{tilrho0lam=mu}
    \rho_{{}_{\tilGam_{\hskip-1pt{}_{[0,\lam]}}(0)}}(0) \ = \ \mu_{H^\b(M,E),H^\b(M,E')}\big(\,
            \rho_{{}_{\Gam_{\hskip-1pt{}_{[0,\lam]}}}}(\n,g^M)\otimes \rho_{{}_{-\Gam_{\hskip-1pt{}_{[0,\lam]}}}}(\n',g^M)\,\big).
\end{equation}
{}From \refe{Gamc} and definition \refe{refinedtor} of the element $\rho$, we conclude that
\eq{rho-Gam}
    \rho_{{}_{-\Gam_{\hskip-1pt{}_{[0,\lam]}}}}(\n',g^M) \ = \ (-1)^{\sum_{j=0}^{r-1}\dim \Ome^j_{[0,\lam]}(M,E')}\cdot
    \rho_{{}_{\Gam_{\hskip-1pt{}_{[0,\lam]}}}}(\n',g^M).
\end{equation}
Hence, by \refe{normmu=product} and \refe{tilrho0lam=mu}, we have
\meq{normrhotilgam0}
    \bigRSnnp{\rho_{{}_{\tilGam_{\hskip-1pt{}_{[0,\lam]}}(0)}}(0)} \ = \
    \bigRSnp{ \rho_{{}_{\Gam_{\hskip-1pt{}_{[0,\lam]}}}}(\n,g^M) }\cdot \RSnp{ \rho_{{}_{\Gam_{\hskip-1pt{}_{[0,\lam]}}}}(\n',g^M) }.
\end{multline}

Combining \refe{tilrho(t)}, \refe{DettilB=BB'}, and \refe{normrhotilgam0} with \refe{rho} and \refe{metricindep} (definition of the refined analytic
torsion), we conclude that
\eq{normrho(0)}
    \RSnnp{\rho(0)} \ = \ \RSnp{\rat(\n)}\cdot \RSnp{\rat(\n')}\cdot e^{2\pi \IM\eta_\lam(\n,g^M)}.
\end{equation}
Hence, by \refe{normrattil0},
\eq{normrho(0)=normtilrho}
    \RSnnp{\rho(0)} \ = \ \RSnnp{\rat(\tiln)}\cdot e^{2\pi \IM\eta_\lam(\n,g^M)}.
\end{equation}

\subsection{The Ray-Singer norm of $\rho(\pi/4)$}\Label{SS:rho(pi/4)}
Recall that the norm $\|\cdot\|_\lam^{{}^\sim}$ was defined in \refss{n+n'}. From \refe{tilDelpi} and \refl{normrho}, we obtain
\eq{normrholam(pi)}\notag
    \big\|\,\rho_{{}_{\tilGam_{\hskip-1pt[0,\lam]}(\pi/4)}}(\pi/4)\,\big\|^{{}^\sim}_\lam \ = \ 1.
\end{equation}
Hence, from \refl{tB(pi/4)=TRS}, \refe{tilrho(t)}, and \refe{RSnormnormprime}, we obtain
\eq{normrho(pi/4)}
    \bigRSnnp{\,\rho(\pi/4)\,} \ = \ 1.
\end{equation}

\subsection{Proof of \reft{RSnormRAT}}\Label{SS:prRSnormRAT}
From \refe{normrho(t)=eta}, \refe{tileta=IMeta}, and \refe{normrho(0)=normtilrho}, we obtain
\eq{rhoeta0=tilrho}
    \bigRSnnp{\, \rho_{{}_{\tilGam_{\hskip-1pt{}_{[0,\lam]}}(0)}}(0)\cdot e^{\xi_\lam(0,\tet)}\,} \ = \ \RSnnp{\rat(\tiln)}.
\end{equation}

By \refe{tilBpi/4}, the operator $\tB^{(\lam,\infty)}_{\even}(\pi/4)$ is self-adjoint. Hence, $\eta_\lam(\pi/4)\in \RR$. Therefore, we get from
\refe{normrho(t)=eta} and \refe{normrho(pi/4)} that
\eq{rhoeta0pi/4=1}
    \bigRSnnp{\, \rho_{{}_{\tilGam_{\hskip-1pt{}_{[0,\lam]}}(\pi/4)}}(\pi/4)\cdot e^{\xi_\lam(\pi/4,\tet)}\,} \ = \ 1.
\end{equation}
From \refe{rhoeta0=tilrho} and \refe{rhoeta0pi/4=1} we conclude that {\em in order to prove \refe{nomrattil=1} (and, hence, \refe{RSnormRAT}) it
suffices to show that
\eq{normrholam(t)xilam(t)}
    \bigRSnnp{\, \rho_{{}_{\tilGam_{\hskip-1pt{}_{[0,\lam]}}(t)}}(t)\cdot e^{\xi_\lam(t,\tet)}\,}
\end{equation}
is independent of\/ $t$. }

{}Fix $t_0\in[-\pi/2,\pi/2]$ and let $\lam\ge0$ be such that the operator $\tB_{\even}(t_0)^2$ has no eigenvalues with absolute value $\lam$. Choose
an angle $\tet\in (-\pi/2,0)$ such that both $\tet$ and $\tet+\pi$ are Agmon angles for $\tB(t_0)$. Then there exists $\del>0$ such that for all
$t\in (t_0-\del,t_0+\del)\cap[-\pi/2,\pi/2]$, the operator $\tB_{\even}(t)^2$ has no eigenvalues with absolute value $\lam$ and both $\tet$ and
$\tet+\pi$ are Agmon angles for $\tB(t)$.

A verbatim repetition of the proof of \refl{xirho} shows that
\eq{ddttilxilam}
    \frac{d}{dt}\, \rho_{{}_{\tilGam_{\hskip-1pt{}_{[0,\lam]}}(t)}}(t)\cdot e^{\xi_\lam(t,\tet)} \ = \ 0.
\end{equation}
Hence, \refe{normrholam(t)xilam(t)} is independent of $t$. \hfill$\square$

\providecommand{\bysame}{\leavevmode\hbox to3em{\hrulefill}\thinspace} \providecommand{\MR}{\relax\ifhmode\unskip\space\fi MR }
\providecommand{\MRhref}[2]{%
  \href{http://www.ams.org/mathscinet-getitem?mr=#1}{#2}
} \providecommand{\href}[2]{#2}

\end{document}